\documentclass[10pt]{article}
\usepackage[in]{fullpage}

\usepackage{amsthm,amsfonts,amssymb,mathrsfs,amsmath,latexsym,mathtools}
\usepackage{amsmath, amssymb, amsthm}
\usepackage[msc-links,initials]{amsrefs}
\usepackage{epsfig}
\usepackage{float}
\usepackage{graphicx}
\usepackage{comment}
\usepackage{enumerate}
\usepackage[shortlabels]{enumitem}
\usepackage{enumitem}
\usepackage{caption}
\usepackage{subcaption}
\usepackage{overpic}
\usepackage[super]{nth}
\usepackage{bbm}
\usepackage{xcolor}
\usepackage{hyperref}
\hypersetup{
 colorlinks   = true, 
 urlcolor     = blue, 
 linkcolor    = blue, 
 citecolor   = red 
}
\usepackage{tikz}

\mathtoolsset{showonlyrefs}

\numberwithin{equation}{section}

\newtheorem{thm}{Theorem}[section]
\newtheorem{prop}[thm]{Proposition}
\newtheorem{lem}[thm]{Lemma}
\newtheorem{cor}[thm]{Corollary}

\theoremstyle{definition}
\newtheorem{definition}[thm]{Definition}

\theoremstyle{remark}

\newcommand{\beq}{ \begin{equation} }
\newcommand{\eeq}{ \end{equation} }
\newcommand{\beqq}{ \begin{equation*} }
\newcommand{\eeqq}{ \end{equation*} }

\DeclareMathOperator{\ai}{Ai}
\DeclareMathOperator{\re}{Re}
\DeclareMathOperator{\im}{Im}
\DeclareMathOperator{\Tr}{Tr}
\DeclareMathOperator*{\res}{Res}

\DeclareMathOperator{\diag}{diag}
\newcommand*{\deff}{\mathrel{\vcenter{\baselineskip0.5ex \lineskiplimit0pt
                     \hbox{\scriptsize.}\hbox{\scriptsize.}}}%
                     =}
\newcommand*{\revdeff}{=\mathrel{\vcenter{\baselineskip0.5ex \lineskiplimit0pt
                     \hbox{\scriptsize.}\hbox{\scriptsize.}}}%
                     }

\newcommand{\Boh}{\mathcal{O}}

\newcommand{\mb}{\mathbb}
\newcommand{\bm}{\mathbf}
\newcommand{\bs}{\boldsymbol}

\newcommand{\wt}{\widetilde}
\newcommand{\wh}{\widehat}

\newcommand{\R}{{\mathbb R}}
\newcommand{\C}{{\mathbb C}}

\newcommand{\Z}{{\mathbb Z}}

\renewcommand{\mid}{\, | \,}

\newcommand{\msf}{\mathsf}

\newcommand{\dd}{\mathrm{d}}
\newcommand{\ii}{\mathrm{i}}

\newcommand{\oz}{\zeta} 

\newcommand{\itR}{\mathsf{R}}  
\newcommand{\itM}{\mathsf{M}}  
\newcommand{\itA}{\mathsf{A}}  

\newcommand{\itt}{\mathsf{t}} 
\newcommand{\itl}{\mathsf{y}} 
\newcommand{\ith}{\mathsf{x}}  
\newcommand{\ittm}{\boldsymbol{\mathsf{t}}} 
\newcommand{\itlm}{\boldsymbol{\mathsf{y}}} 
\newcommand{\ithm}{\boldsymbol{\mathsf{x}}}  
\newcommand{\itF}{\mathsf{F}}  
\newcommand{\itD}{\mathfrak{D}^{\mathrm{(KPZ)}}} 
\newcommand{\itP}{\mathsf{P}}  
\newcommand{\itQ}{\mathsf{Q}}  
\newcommand{\itK}{\mathsf{K}}  
\newcommand{\itid}{\mathsf{1}}  
\newcommand{\itG}{\mathsf{G}}  
\newcommand{\ite}{\mathsf{e}}  
\newcommand{\itU}{\mathsf{U}}  
\newcommand{\itLambda}{\mathsf{\Lambda}}  
\newcommand{\itH}{\mathsf{H}}  

\newcommand{\itJ}{\mathsf{J}}  
\newcommand{\itI}{\mathsf{I}}  
\newcommand{\itzero}{\mathsf{0}}  
\newcommand{\itperm}{\mathsf{\Pi}}  
  
\newcommand{\itconjb}{\mathsf{\Delta}} 
\newcommand{\itYb}{\mathsf{\Phi}}  
\newcommand{\itC}{\mathsf{C}}  

\newcommand{\itE}{\mathsf{E}}  

\newcommand{\itss}{\mathsf{S}}  
\newcommand{\itpermss}{\mathsf{\Pi}^{\itss}}  
\newcommand{\itEss}{\mathsf{E}^{\itss}}  
\newcommand{\itqss}{\mathsf{q}} 
\newcommand{\itpss}{\mathsf{p}} 
\newcommand{\itrss}{\mathsf{r}} 
\newcommand{\itsss}{\mathsf{s}} 
\newcommand{\ituss}{\mathsf{u}} 

\newcommand{\itDsnt}{\partial_{\itt}} 
\newcommand{\itDsnl}{\partial_{\itl}} 
\newcommand{\itDsnh}{\partial_{\ith}} 
\newcommand{\itEsn}{\mathsf{E}}  
\newcommand{\itAsn}{\mathsf{C}}  
\newcommand{\itqsn}{\mathsf{q}}  
\newcommand{\itpsn}{\mathsf{p}}  
\newcommand{\itrsn}{\mathsf{r}}  
\newcommand{\itssn}{\mathsf{s}}  

\newcommand{\prob}{\mathbbm{P}}

\newcommand{\kpt}{\tau} 
\newcommand{\kpz}{\mathrm{(KPZ)}}  
\newcommand{\kpga}{\gamma} 
\newcommand{\kpf}{\mathsf{f}} 

\newcommand{\kpheight}{\mathcal{H}} 
\newcommand{\kph}{\mathsf{h}}  
\newcommand{\kpc}{\mathsf{c}}  
\newcommand{\kpqq}{\mathsf{a}}  
\newcommand{\kprr}{\mathsf{b}}  

\newcommand{\chim}{\boldsymbol{\mathsf{X}}} 
  
\newcommand{\itGo}{\bm M_{\mathsf{o}}}  
\newcommand{\itGe}{\bm M_{\mathsf{e}}}  
\newcommand{\Am}{\boldsymbol{\mathsf{A}}} 
\newcommand{\Bm}{\boldsymbol{\mathsf{B}}} 
\newcommand{\Mm}{\boldsymbol{\mathsf{M}}} 
\newcommand{\Fm}{\boldsymbol{\mathsf{U}}} 
\newcommand{\Gm}{\boldsymbol{\mathsf{V}}} 
\newcommand{\Gmo}{\bm M_{\mathsf{o}}}  
\newcommand{\Gme}{\bm M_{\mathsf{e}}}  
\newcommand{\Sm}{\bm S}

\newcommand{\kpzFdist}{\mathfrak{F}^{\mathrm{(KPZ)}}} 
\newcommand{\kpL}{\mathsf{L}}  
\newcommand{\kpR}{\mathsf{R}}  

\newcommand{\pkmath}{\mathsf}  
\newcommand{\per}{\mathrm{(per)}}  
\newcommand{\pkt}{\tau} 
\newcommand{\pkga}{\gamma} 
\newcommand{\pkh}{\pkmath{h}}  
\newcommand{\pkS}{\pkmath{S}}  
\newcommand{\pkQ}{\pkmath{Q}}  
\newcommand{\pkH}{\pkmath{V}}  
\newcommand{\pkA}{\pkmath{A}}  
\newcommand{\pkB}{\pkmath{B}}  

\newcommand{\pkmchim}{\bs{\pkmath{X}}} 

\newcommand{\pkFdist}{\mathfrak{F}^{\mathrm{(per)}}} 
\newcommand{\pkDdisc}{\mathfrak{D}^{\mathrm{(per)}}}  
\newcommand{\pkCdisc}{\mathfrak{C}^{\mathrm{(per)}}}  
\newcommand{\pkL}{\pkmath{L}}  
\newcommand{\pkR}{\pkmath{R}}  

\newcommand{\ail}{y} 
\newcommand{\aih}{x}  
\newcommand{\aitau}{t} 
\newcommand{\aiE}{E} 
\newcommand{\aiW}{W}  

\newcommand{\laxD}{\msf D} 
\newcommand{\laxC}{\msf C} 

\newcommand{\laxKP}{\msf B} 

\newcommand{\genchi}{\bs{\mathsf{X}}} 

\newcommand{\gena}{\mathsf{a}}  
\newcommand{\genb}{\mathsf{b}}  
\newcommand{\genc}{\mathsf{c}}  
\newcommand{\genf}{\mathsf{f}}  
\newcommand{\geng}{\mathsf{g}}  
\newcommand{\genF}{\mathsf{F}}  
\newcommand{\genM}{\mathsf{M}}  
\newcommand{\genG}{\mathsf{G}}  
\newcommand{\genK}{\mathsf{K}}  
\newcommand{\genH}{\mathsf{H}} 
\newcommand{\genT}{\mathsf{T}} 
\newcommand{\genL}{\mathsf{L}} 
\newcommand{\genm}{\mathsf{m}}  
\newcommand{\genFunc}{\mathfrak{D}}  

\newcommand{\genq}{\mathsf{q}} 
\newcommand{\genp}{\mathsf{p}} 
\newcommand{\genr}{\mathsf{r}} 
\newcommand{\gens}{\mathsf{s}} 
\newcommand{\genu}{\mathsf{u}} 
\newcommand{\genv}{\mathsf{v}} 

\newcommand{\genD}{\bs{\mathsf{D}}} 

\newcommand{\genHilb}{L^2(\Omega,\mu)} 

\newcommand{\gent}{\mathsf{t}}  
\newcommand{\geny}{\mathsf{y}}  
\newcommand{\genx}{\mathsf{x}}  

\newcommand{\gendetD}{\mathfrak{D}}  

\DeclareMathOperator{\Fgue}{F_{GUE}}

\newcommand{\itB}{\mathsf{B}}  
\newcommand{\itHe}{\mathsf{X}}  
\newcommand{\itPo}{\mathsf{P}}  
\newcommand{\itZ}{\mathsf{\Phi}}  
\newcommand{\itZsn}{\mathsf{\Phi}}  
\newcommand{\itW}{\mathsf{\Psi}}  
\newcommand{\itWsn}{\mathsf{\Psi}}  
\newcommand{\itEn}{\mathsf{E}}  

\newcommand{\kpYd}{\mathsf{L}}  
\newcommand{\kpYo}{\mathsf{O}}  
\newcommand{\kpV}{\mathsf{V}}  
\newcommand{\kpU}{\mathsf{U}}  
\newcommand{\itvss}{\mathsf{v}} 
\newcommand{\whkptm}{\itM_\mathsf{t}}  

\newcommand{\whkpym}{\itM_\mathsf{y}}  

\newcommand{\whkpxm}{\itM_\mathsf{x}}  

\newcommand{\kpMz}{\mathsf{M}_{\mathsf{z}}}  

\usepackage{upgreek}

\newcommand{\veca}{a}  

\title{Differential equations for the KPZ and periodic KPZ fixed points}

\author{Jinho Baik\footnote{Department of Mathematics, University of Michigan,
Ann Arbor, MI, 48109, USA, \texttt{baik@umich.edu}} \and Andrei Prokhorov\footnote{Department of Mathematics, University of Michigan,
Ann Arbor, MI, 48109, USA, \texttt{andreip@umich.edu}}\,\,\footnote{St. Petersburg State University, Universitetskaya emb. 7/9, 199034, St. Petersburg, Russia}
\and Guilherme L.~F.~Silva\footnote{Instituto de Ci\^encias Matem\'aticas e de Computa\c{c}\~ao, Universidade de S\~ao Paulo, S\~ao Carlos, SP, 13566-590, Brazil, \texttt{silvag@usp.br}}}

\date{\today}

\begin{document}

\maketitle

\begin{abstract}
The KPZ fixed point is a 2d random field, conjectured to be the universal limiting fluctuation field for the height function of models in the KPZ universality class.
Similarly, the periodic KPZ fixed point is a conjectured universal field for spatially periodic models. 
For both fields, their multi-point distributions in the space-time domain have been computed recently. 
We show that for the case of the narrow-wedge initial condition, these multi-point distributions can be expressed in terms of so-called integrable operators. 
We then consider a class of operators that include the ones arising from the KPZ and the periodic KPZ fixed points, and find that they are related to various matrix integrable differential equations such as  coupled matrix mKdV equations, coupled matrix NLS equations with complex time, and matrix KP-II equations. 
When applied to the KPZ fixed points, our results extend previously known differential equations for one-point distributions and equal-time, multi-position distributions to multi-time, multi-position setup. 
\end{abstract}

\setcounter{tocdepth}{1}
\tableofcontents



\section{Introduction}

The KPZ fixed point 
\beqq
	\kpheight(\kpga, \kpt), \qquad (\kpga, \kpt)\in \R\times \R_+,
\eeqq
is a 1+1 dimensional random field that is conjectured to be the universal limit for the height fluctuations of the random growth models  belonging to the KPZ universality class. 
Many remarkable properties of the KPZ fixed point were obtained over the last two decades, and 
several models are proved to converge to it in various senses. 
The field itself was constructed relatively recently by Matetski, Quastel and Remenik \cite{Matetski-Quastel-Remenik17}. 
Another construction was also obtained by Dauvergne, Ortmann and Virág \cite{Dauvergne-Ortmann-Virag18}. 
In this paper we focus on the multi-point distribution functions and study their connections to deterministic  (integrable) differential equations. 
Such connections may depend on the initial condition of the KPZ fixed point. 
We consider the most well-studied initial condition, the narrow wedge initial condition given by 
$\kpheight(\kpga,\kpt=0)= 0$ for $\gamma=0$ and $\kpheight(\kpga,\kpt=0)= -\infty$ for $\gamma \neq 0$. 
In this case,  the $1:2:3$ invariance property,
\begin{equation}\label{eq:KPZscalinginv}
	\left( \epsilon^{-1} \kpheight(\epsilon^2 \kpga, \epsilon^3\kpt) \right)_{\kpga, \kpt} \overset{d}{=} \left( \kpheight(\kpga, \kpt) \right)_{\kpga, \kpt},
\end{equation}
holds for all $\epsilon>0$. 

Many differential equations are known for the one-time distributions and the equal-time, multi-position distributions. 
Some of them (for the narrow-wedge initial condition) are the following.
\begin{itemize}
\item The one-point distribution reduces to the GUE Tracy-Widom distribution $\Fgue$. For each fixed $\gamma, \tau$, 
\begin{equation}\label{eq:Fstepgue1}
	\prob\left( \kpheight(\kpga,\kpt)\leq \kph  \right)=\Fgue\left( \frac{\kph}{\tau^{1/3}}+\frac{\gamma^2}{\tau^{4/3}} \right),
\end{equation}
where $\Fgue$ is the GUE Tracy-Widom distribution. It  satisfies  \cite{Tracy-Widom94}
\begin{equation}\label{eq:FstepTW1}
	\frac{\dd^2}{\dd \xi^2} \log \Fgue(\xi) = -u(\xi)^2,
\end{equation}
where $u$ solve the Painlevé II equation 
\begin{equation}\label{eq:FstepTW2}
	u''= \xi u^2+2 u^3
\end{equation}
with the boundary condition $u(\xi)\sim \ai(\xi)$ as $\xi\to +\infty$. 

\item The \underline{equal-time}, multi-position distribution is a scaled and shifted version of the Airy$_2$ process\footnote{The Airy$_2$ process is often regarded as an evolution. The time in the Airy$_2$ process corresponds to the position of the KPZ fixed point.}  introduced by Pr{\"a}hofer and Spohn \cite{Prahofer-Spohn02}. 
In 2005, two differential equations were obtained for the multi-point distribution of the Airy$_2$ process. 
\begin{itemize}
\item Tracy and Widom \cite{Tracy-Widom05} found a system of matrix ordinary differential equations (ODEs)  with respect to a combination of the height variables. 
\item On the other hand, Adler and van Moerbeke \cite{Adler-van_Moerbeke05} obtained a nonlinear third order partial differential equation (PDE) for the case of the 2-point distribution with respect to the (shifted and scaled) height and spatial variables. Their result was extended to multi-point distribution by Wang \cite{Wang09a}. 

\end{itemize}
In these two results, the time variable was kept as a constant and did not appear in the equations.
\item Quastel and Remenik \cite{Quastel-Remenik19b} also considered  the equal-time, multi-point distribution, but included the time $\tau$ as an evolution variable instead of a fixed parameter. 
Thus, they regarded the \underline{equal-time} $m$-point multi-position distribution as a function on $2m+1$ variables. 
They derived a matrix version of the second  Kadomtsev-Petviashvili (KP-II) equation with respect to these $2m+1$ variables. 
If we specialize their result to the $m=1$ case, it implies that \eqref{eq:Fstepgue1}, viewed as a function of three scalar variables $x, \gamma, \tau$, is related to the scalar KP-II equation. 
Remarkably, the result of Quastel and Remenik applies to general initial conditions far beyond the narrow wedge initial condition. 
\end{itemize}

One of the goals of this paper is to extent the above results to \underline{multi-time}, multi-position distributions of the KPZ fixed point with the narrow-wedge initial condition. 
Fix a positive integer $m$ and consider the function of $3m$ variables 
\begin{equation} \label{eq:Fstepdefnintro}
	\kpzFdist(\kph, \kpga, \kpt)\deff	\prob\left( \bigcap_{i=1}^m \{\kpheight(\kpga_i, \kpt_i) \le \kph_i \} \right),
\end{equation}
where $\kpt=(\kpt_1, \cdots ,\kpt_m)^T \in \R_+^m$, $\kpga=(\kpga_1, \cdots, \kpga_m)^T\in \R^m$, and $\kph=(\kph_1, \cdots, \kph_m)^T \in \R^m$
represent time, position, and height, respectively. 
We prove that this function is related to the following differential equations. 
\begin{enumerate}
\item As a function of $\kpga$ and $\kph$, we obtain a system of coupled $m\times m$ matrix nonlinear Schr\"odinger (NLS) equations with complex time $\ii \gamma$.

\item As a function of $\kpt$ and $\kph$, we obtain a system of coupled matrix modified Korteweg-de Vries (mKdV) equations. 

\item As a function of all $3m$ variables, we obtain matrix Kadomtsev-Petviashvili (KP)  equations that generalizes the matrix KP-II equation obtained by Quastel and Remenik \cite{Quastel-Remenik19b} for the case $\tau_1=\cdots= \tau_m$. 
We also show a connection to the multi-component KP hierarchy.

\item As a function of $\kph$, we obtain a matrix ODE system that generalizes the matrix ODE system of Tracy and Widom \cite{Tracy-Widom05} for the Airy$_2$ process. 

\end{enumerate}

For the case $m=2$, we also derive several PDEs with the aid of symbolic computations. If we take the special case $\tau_1=\tau_2$, then a combination of two of them becomes the PDE of Adler and van Moerbeke \cite{Adler-van_Moerbeke05} for the Airy$_2$ process. 

\bigskip

We also study a periodic version of the KPZ fixed point. 
If we consider random growth models in the KPZ universality class on a periodic domain (instead of the infinite line) and take a large time, large period limit in a certain critical way, then a new (1+1) dimensional random field emerges. 
This field, which we call the periodic KPZ fixed point, is expected to interpolate between the Brownian motion and the KPZ fixed point \cites{Baik-Liu19, Baik22}. 
This interpolation property was proved for the one-point distribution in \cite{Baik-Liu-Silva}. Even though the periodic KPZ fixed point is not yet constructed as a field, its multi-point distributions were obtained recently as limiting distributions for the totally asymmetric simple exclusion process on a ring \cite{Baik-Liu19}. 
It was shown in \cite{Baik-Liu-Silva} that the one-point distribution of the periodic KPZ fixed point with the narrow step initial condition is related to a coupled  NLS (with complex time), a coupled  mKdV, and a KP-II equation.
In this paper, we extend this result to the multi-point distributions and show that they are related to the same equations as the KPZ fixed point mentioned in 1--3 above. 
However, the Tracy-Widom type ODE system is not expected to hold for the periodic KPZ fixed point. In particular, the one-point distribution is not the Tracy-Widom distribution; it depends non-trivially on the time. 

\medskip

The multi-time, multi-position distributions of both the KPZ fixed point and the periodic KPZ fixed point have been computed explicitly and are expressible in terms of integrals of Fredholm determinants \cites{Johansson-Rahman19, Liu19, Baik-Liu19}. 
We will show that the Fredholm determinants for both fields are equal to the Fredholm determinants of so-called IIKS integrable operators \cites{IIKS, Deift99b}. 
These integrable operators have a certain common structure which we call cubic integrable. We derive the differential equations 1--3, and also 4 in some cases, for general cubic integrable operators.  
Such cubic integrable operators can be either defined on contours, which is the case of the KPZ fixed point, or on discrete sets in the complex plane, which is the case of the periodic fixed point.

IIKS integrable operators have a special algebraic structure that associates them to Riemann-Hilbert problems in a canonical way.
In turn, Riemann-Hilbert problems are deeply integrated in the theory of integrable differential equations, providing inverse transform methods for the latter \cites{Ablowitz-Clarkson91, Its03, Fokas-Its-Kapaev06}. 
In particular, if a Riemann-Hilbert problem depends on parameters and the jump/residue matrix for the Riemann-Hilbert problem has a certain decomposition  structure, then one can derive linear matrix differential equations, called Lax equations, for its solution. 
The compatibility of a pair of Lax equations often yields nonlinear differential equations in those parameters. 

For the case of cubic integrable operators, we find Lax equations for each of the parameters, $\kph, \kpga, \kpt$. 
The compatibility between the Lax equations in the variables $\kph$ and $\kpga$ yields the coupled NLS equation mentioned in 1 above, the compatibility between $\kpt$ and $\kph$ yields the coupled mKdV mentioned in 2, and the compatibility between the equations on the three variables yield the KP-II equation and multi-component KP hierarchy from 3.
If we impose an additional property on the cubic integrable operators, that we call \underline{strong} cubic integrability, we obtain an additional Lax equation with respect to the so-called spectral variable. This additional equation allows us to obtain the equation in 4. 
We will see that the cubic integrable operator for the KPZ fixed point is strongly cubic integrable hence the equation 4 holds.

\bigskip

In the next section we introduce the cubic integrable operators formally, and state our main results. 
See Subsection~\ref{sec:organization} for the organization of the rest of the paper.

\subsection*{Acknowledgments}

The work of Baik was supported in part by NSF grant DMS-1954790.
Silva acknowledges his current support by São Paulo Research Foundation under grants \# 2019/16062-1 and \# 2020/02506-2, and by Brazilian National Council for Scientific and Technological Development (CNPq) under grant \# 315256/2020-6. Prokhorov was supported by NSF MSPRF grant DMS-2103354 and RSF grant 22-11-00070.

\section{Statement of Main Results}\label{sec:mainresults}

In the first subsection, we define cubic integrable operators and state various differential equations associated to them. 
In the next subsection we introduce the Fredholm determinants of cubic integrable operators. 
The following subsection then discusses how the multi-time, multi-position distributions of the KPZ and the periodic KPZ fixed points are related to cubic integrable operators. 
We finish this section with a roadmap of the organization of the rest of the paper. 

\subsection{Cubic integrable operators and differential equations} \label{sec:DEsection}

Let $\Omega\subset \C$ be either a finite union of disjoint simple contours or a discrete set without accumulations points. 
If it is a union of contours, we further assume that each contour is either closed or extends to infinity, so in particular none of these contours have finite endpoints. 
Let 
\beq \label{eq:mudeff}
	\text{$\mu$ be the  counting measure if $\Omega$ is discrete and $\dd\mu=\dd z$ if $\Omega$ consists of contours,} 
\eeq
so that 
\begin{equation}\label{deff:genHilb}
	\genHilb =  
	\begin{cases}
	\ell^2(\Omega) & \text{if } \Omega \text{ is discrete}, \\
	L^2(\Omega, \dd z) & \text{if } \Omega \text{ is a union of contours}. 
\end{cases}
\end{equation}
Fix an integer $m\geq 1$ and consider $3m$ parameters 
\beq \label{eq:tyxde}
	\gent=(\gent_1,\hdots,\gent_m)\in \R^m,\qquad \geny=(\geny_1,\hdots,\geny_m)\in \R^m, \qquad \genx=(\genx_1,\hdots,\genx_m)\in \R^m. 
\eeq
In the (periodic) KPZ fixed point, they will correspond to time, spatial location, and height variables.
Define the cubic exponential functions 
\begin{equation}\label{def:mjexpon}
	\genm_j(z) 
	\deff \exp \left( \gent_j z^3 + \geny_j z^2 + \genx_j z\right), \qquad z\in \C, 
\end{equation}
and introduce the $(m+1)\times (m+1)$ diagonal matrix-valued function
\begin{equation}\label{eq:itconjb}
\itconjb(z)=	\itconjb(z \mid \genx,\geny,\gent)\deff \diag\left(\genm_1(z),\genm_2(z),\hdots,\genm_m(z),1\right). 
\end{equation}
Let, for $i\in \{1. \cdots, m+1\}$, 
\beq \label{eq:Eidefn}
	\text{$\itE_{i}$ be the $(m+1)\times(m+1)$ matrix whose $(i,i)$ entry is $1$ and all other entries are $0$.}
\eeq

An integral operator $\genH: \genHilb\to \genHilb$ is called {\it integrable} in the sense of Its, Izergin, Korepin and Slavnov (or, shortly, an IIKS operator, or an integrable operator) if its kernel is of the form\footnote{The condition $\genH(u,u)=0$ is not strictly necessarily. We included it here to simplify technical details in the theory, especially when $\Omega$ is a discrete set. This simplification is enough for our purposes regarding the (periodic) KPZ fixed points.} 
\begin{equation}\label{eq:admisskernel}
	\genH(u,v)=\frac{\genf(u)^T\geng(v)}{u-v} \quad \text{for }u\neq v\qquad \text{and}\qquad \genH(u,u)=0 \quad \text{for }u,v\in \Omega, 
\end{equation}
for some vector-valued functions $\genf$ and $\geng$ (see \cites{IIKS, Deift99b}).  
The operators of relevance to us possess an additional structure on the vector functions $\genf$ and $\geng$. 

\begin{definition}[Cubic integrable operator] \label{def:cubicintegrabl} 
A bounded integral operator $\genH:\genHilb\to\genHilb$ is called {\it cubic integrable} if it is integrable in the sense of \eqref{eq:admisskernel} and the vector functions are of the form 
\begin{equation}\label{eq:fg}
	\genf(u)=\genf(u\mid \genx,\geny,\gent) \deff \genc(u) \itconjb(u\mid \genx,\geny,\gent) \Fm(u) , 	 
	\qquad  
	\geng(u)=\geng(u\mid \genx,\geny,\gent) \deff \frac1{\genc(u)} \itconjb(u\mid \genx,\geny,\gent)^{-1} \Gm(u)
\end{equation}
for $u\in \Omega$, where $\itconjb$ is as in \eqref{eq:itconjb},
with the following additional structures:
\begin{itemize}
    \item $\Fm$ and $\Gm$ are  $(m+1)$-dimensional column vector functions which do not depend on $\genx,\geny,\gent$ and satisfy
\begin{equation} \label{eq:withEzero}	
	\Fm(u)^T \itE_{i}  \Gm(u)=0 \qquad \text{ for every } u\in \Omega \text{ and }i\in \{1, \cdots, m+1\};
\end{equation}
	\item $\genc$ is a non-vanishing scalar function  on $\Omega$ (which may depend on $\genx,\geny,\gent$); 
    \item $\genf,\geng \in L^2(\Omega,\mu)\cap L^\infty(\Omega,\mu) $. 
\end{itemize}
If, in addition, the set $\Omega$ is a union of disjoint contours and the vectors $\Fm,\Gm$ are constants on each connected component of $ \Omega$, then we say that $\genH$ is a {\it strongly} cubic integrable operator.
\end{definition}

If we change the scalar $\msf c$, the operator is a conjugate of the original operator $\genH$. 
Since we will be mostly interested in properties of the operators that are invariant under such conjugations, such as the Fredholm determinant and the trace, the scalar $\msf c$ does not play a major role in our results.
However, including it in the definition as a free parameter is useful, as it provides an extra layer of flexibility that we explore in examples, to prove that $\genH$ is not only bounded but also a trace class operator.

The term {\it cubic} comes from the cubic dependence of the kernel of $\genH$ on the parameters $\genx,\geny,\gent$ through the factor $\itconjb$ in \eqref{eq:fg}. 

\bigskip

For a cubic integrable operator, we define the following matrix function. 

\begin{definition} \label{defn:Phi1d}
Let $\genH$ be a cubic integrable operator on $\genHilb$ and assume that $\itid-\genH$ is invertible. 
Define the $(m+1)\times (m+1)$ matrix 
\begin{equation}\label{defn:itYb1}
	\itYb_1=\itYb_1(\genx,\geny,\gent)\deff \int_\Omega ((\itid-\genH)^{-1} \genf)(s)\geng(s)^T \dd\mu(s) 
\end{equation}
where recall \eqref{eq:mudeff} for $\mu$. 
\end{definition}

The subscript $1$ is used in \eqref{defn:itYb1}  since $\itYb_1$ is the first element of a sequence of matrix functions; see \eqref{eq:expYk} below. 
Observe that the integral \eqref{defn:itYb1} is unchanged even if we change the scalar $\genc$ in the definition of $\genf, \geng$ to a different scalar.

A non-empty set $\itss\subset \{1,\cdots,m\}$ induces a permutation $\pi$ on $\{1, \cdots, m+1\}$ that fixes the element $m+1$; this permutation is given by $\pi(i)=k_i$ for $1\le i\le m$ and $\pi(m+1)=m+1$, where $k_1< \cdots<k_{|\itss|}$ are elements of $\itss$ and $k_{|\itss|+1}< \cdots<k_{m}$ are elements of $\{1, \cdots, m\}\setminus\itss$. 
For a non-empty set $\itss\subset \{1,\cdots,m\}$, let  $\itpermss=(\delta_{\pi(i),j})_{i, j=1}^{m+1}$ denote the  $(m+1)\times (m+1)$ permutation matrix for $\pi$. Note that the $(m+1, m+1)$ entry of the matrix $\itpermss$ is $1$. 

\begin{definition} \label{defn:Phioneblock}
Let $\itZ_1$ be the matrix in \eqref{defn:itYb1}. 
For a non-empty set $\itss\subset \{1,\cdots,m\}$, introduce matrix-valued functions $\itqss,\itpss,\itrss,\itsss$ through the formula 
\beq \label{eq:L1O1not}
	\itpermss \itZ_1 (\itpermss)^T  = 
	\begin{pmatrix} \itqss & \itpss \\ \itrss & \itsss \end{pmatrix},
\eeq 
with $\itqss$ and $\itsss$ being of sizes $|\itss|\times |\itss|$ and  $(m+1-|\itss|)\times (m+1-|\itss|)$, respectively. 
Also, introduce the differential operators
\begin{equation}\label{eq:diffoperitss}
	\partial_\itt=\partial^\itss_{\itt}\deff \sum \limits_{k\in \itss}\partial_{\itt_{k}}  ,\qquad \partial_\itl=\partial^\itss_{\itl}\deff \sum_{k\in\itss} \partial_{\itl_{k}} \qquad \text{and}\qquad \partial_\ith=\partial^\itss_{\ith}\deff \sum_{k\in \itss}\partial_{\ith_k}. 
\end{equation}
\end{definition}

The matrices $\itqss,\itpss,\itrss,\itsss$ depend on the permutation $\itss$, but we do not display this dependence explicitly to simplify the notations. 
We note that the functions $\itqss,\itpss,\itrss,\itsss$ are complex-valued matrices. 
The following are basic relations between them. 

\begin{lem} \label{lem:basic}
Under the notations of Definition~\ref{defn:Phi1d},
$$
	\Tr(\genq) + \Tr(\gens)=0, \qquad 
	\itDsnh \genq=-\genp\genr, \qquad 
	\itDsnh \itsss = \genr\genp. 
$$
\end{lem}

\bigskip

We now state three main theorems on differential equations. The first theorem concerns the pair $\genp$ and $\genr$, which are matrices of sizes $|\itss|\times (m+1-|\itss|)$ and $(m+1-|\itss|)\times |\itss|$, respectively. 

\begin{thm}[Coupled matrix NLS with complex time and mKdV] \label{thm:KdVheat}
Let $\genH$ be a cubic integrable operator and assume that $\itid-\genH$ is invertible for the parameters $(\genx,\geny,\gent)$ in an open set of $\R^{3m}$. 
Then, the following results hold. 
\begin{enumerate}[(a)]
\item As functions of $\ith$ and $\itl$, the matrices $\genp$ and $\genr$ satisfy a system of coupled matrix nonlinear Schr\"odinger (NLS) equations  with complex time $\itl\mapsto \ii \itl$, 
\beq\label{eq:heat}
	\partial_{\itl}\genp=\partial^2_{\ith}\genp+2\genp\genr\genp,
	\qquad 
	\partial_{\itl}\genr=-\partial^2_\ith\genr-2\genr\genp\genr .
\eeq 

\item As functions of $\itt$ and $\ith$, the matrices $\genp$ and $\genr$ satisfy a system of coupled matrix modified Korteweg-de Vries (mKdV) equations 
\beq\label{eq:mkdv}
	\partial_\itt\genp=\partial^3_\ith\genp+3(\partial_\ith\genp)\genr\genp+3\genp\genr(\partial_\ith \genp), 
	\qquad 
	\partial_\itt\genr=\partial^3_\ith\genr+3(\partial_\ith\genr)\genp\genr+3\genr\genp(\partial_\ith \genr) . 
\eeq
\end{enumerate}
\end{thm} 

These equations were obtained for the (periodic) KPZ fixed point when $m=1$ in  \cite{Baik-Liu-Silva}. 
The coupled NLS equations with complex time with $m=1$ also appeared in the recent work \cite{Krajenbrink_Le_Doussal} on the weak noise theory of the KPZ equation. 

\medskip

The next differential equations involve all $3m$ variables $\itt$, $\itl$ and $\ith$. 

\begin{thm}[Matrix KP] \label{thm:KPeqtion}
Let $\genH$ be a cubic integrable operator and assume that $\itid-\genH$ is invertible for the parameters $(\genx,\geny,\gent)$ in an open set of $\R^{3m}$.
Then, the following result hols. 
\begin{enumerate}[(a)]
\item The $|\itss|\times |\itss|$ matrices $\genu\deff\genp\genr$ and $\genq$ satisfy the matrix KP-II equation\footnote{The usual KP-II equation is 
$ \partial_x \left( \partial_t w + 6w \partial_x w +  \partial_{xxx} w \right)+ 3\partial_{yy} w=0$ for a scalar function $w$. 
When $m=1$, the $x$-derivative of Equation~\eqref{eq:kpu} becomes this equation if we set $w(t,y,x)=2u(-4t,y,x)$.
}  
\begin{equation}\label{eq:kpu}
	-4\partial_\itt\genu+\partial_\ith^3\genu+6\partial_\ith(\genu^2) - 3\partial_\itl^2\genq  +6[\genu,\partial_\itl \genq] =0,
	\qquad 
	\partial_\ith \genq=-\genu.
\end{equation} 

\item The $(m+1-|\itss|)\times(m+1- |\itss|)$ matrices $\itvss \deff  \itrss \itpss$ and $\gens$ 
satisfy the matrix KP-II equation 
\beq\label{eq:kpv}
	- 4 \itDsnt \itvss + \itDsnh^3 \itvss  + 6 \itDsnh ( \itvss^2) + 3  \itDsnl^2 \itsss  + 6 [ \itvss,   \itDsnl \itsss ]   =0, 
	\qquad 
	\itDsnh \itsss = \itvss.
\eeq
\end{enumerate}
\end{thm}

If we insert $\genu=- \partial_\ith \genq$, \eqref{eq:kpu} becomes a fourth order differential equation for the matrix $\genq$ only. 
If we take $\itss=\{1, \cdots, m\}$ and set $\itt_1=\cdots=\itt_m$, Equation~\eqref{eq:kpu} becomes the matrix KP-II equation obtained in  \cite{Quastel-Remenik19b} for the KPZ fixed point; see Section~\ref{sec:previousdes} below for further discussion on this end.

We also obtain a connection to the so-called multi-component KP hierarchy in Proposition~\ref{prop:KPhi} below. 

\medskip

For the next result, define $(m+1)\times (m+1)$ matrices 
\beq \label{eq:multibytyx}
	\whkptm \deff \diag( \itt_1, \cdots, \itt_m, 0) , \quad
	\whkpym \deff \diag( \itl_1, \cdots, \itl_m, 0) , \quad
	\whkpxm\deff \diag( \ith_1, \cdots, \ith_m, 0),
\eeq 
and the differential operator 
\begin{equation}\label{def:diffoperPform}
	\partial  \deff \sum_{k=1}^m \itt_k \partial_{\ith_k}.
\end{equation}
We take  $\itss=\{1, \cdots, m\}$ in Definition~\ref{defn:itYb1} so that $\itpermss$ is the identity matrix, and 
$\genq$, $\genp$, $\genr$, and $\gens$ are matrices of sizes $m\times m$, $m\times 1$, $1\times m$, and $1\times 1$, respectively.

\begin{thm}[ODE system] \label{thm:matrixODE}
Let $\genH$ be a \underline{strongly} cubic integrable operator and assume that $\itid-\genH$ is invertible for the parameters $(\genx,\geny,\gent)$ in an open set.
Then,  $\itZ_1$ in \eqref{defn:itYb1} satisfies
\beq \label{eq:matrxp2ii}
	3\partial^2\itZ_1 = 2[\partial\itZ_1, \whkpym] + [[\itZ_1, \whkptm], 3\partial\itZ_1-2[\itZ_1, \whkpym] -\whkpxm] .
\eeq 
Furthermore, the additional relation
\beq \label{eq:parsinrmp}
	\partial\itssn=\itrsn \ittm \itpsn
\eeq
also holds, where $\ittm$ is the matrix defined in \eqref{def:boldtyx} below. 
\end{thm} 

If we write the block form \eqref{eq:matrxp2ii} (recall that $\itpermss=\itI$ in this case) and use \eqref{eq:parsinrmp}, then we find the following matrix ODE system for $\itqsn, \itpsn, \itrsn$.
Introduce the $m\times m$ matrices
\beq \label{def:boldtyx}
	\ittm\deff \diag(\itt_1,\hdots,\itt_m),\qquad 
	\itlm\deff \diag(\itl_1,\hdots,\itl_m),\qquad 
	\ithm\deff \diag(\ith_1,\hdots,\ith_m). 
\eeq 

\begin{cor} \label{thm:matrixp2}
Under the same conditions of Theorem~\ref{thm:matrixODE}, $\itqsn, \itpsn, \itrsn$ satisfy the matrix ODE system 
\begin{equation}\label{eq:thPformSystem}
\begin{split}
	& 3\partial^2\itqsn-2[\partial \itqsn, \itlm] -3 [[\itqsn, \ittm], \partial\itqsn] +3 \ittm \itpsn \partial\itrsn+3 ( \partial\itpsn ) \itrsn\ittm  +2 [[\itqsn, \ittm], [\itqsn, \itlm]] -2\ittm \itpsn \itrsn\itlm +2 \itlm \itpsn \itrsn\ittm + [[\itqsn, \ittm], \ithm] =0, \\
	& 3\partial^2\itpsn+2\itlm \partial\itpsn  - 3 ( \partial\itqsn ) \ittm \itpsn - 3[\itqsn, \ittm]\partial\itpsn  +3\ittm \itpsn\itrsn\ittm \itpsn   -2 (\itlm\itqsn\ittm-\ittm\itqsn\itlm)\itpsn +\ithm\ittm \itpsn =0, \\
	& 3\partial^2\itrsn-2 \partial\itrsn\itlm  -3\itrsn\ittm \partial\itqsn+3\partial\itrsn[\itqsn, \ittm] + 3\itrsn\ittm \itpsn\itrsn\ittm -2\itrsn (\itlm\itqsn\ittm-\ittm\itqsn\itlm)+ \itrsn\ittm \ithm =0.
\end{split}
\end{equation}
In the particular case when $\itt_1=\cdots=\itt_m\revdeff \itt$, they satisfy the reduced ODE system 
\begin{equation}\label{eq:thPformSystemred}
\begin{split}
	& \partial_{\ith} \itqsn+\itpsn\itrsn=0, \\
	& 3\itt \partial_{\ith}^2\itpsn+ 2\itlm  \partial_{\ith} \itpsn +6\itt  \itpsn\itrsn \itpsn  - 2[\itlm,\itqsn]  \itpsn  +\ithm \itpsn =0, \\
	& 3\itt \partial^2_{\ith}\itrsn- 2  ( \partial_{\ith}\itrsn ) \itlm +6\itt \itrsn\itpsn\itrsn -  2\itrsn [\itlm,\itqsn] + \itrsn \ithm =0. \\
\end{split}
\end{equation}
where $\partial_{\ith}\deff \sum_{i=1}^m \partial_{\ith_i}$.
\end{cor}

As we discuss in Subsection~\ref{sec:multipointairydiscussion} below, the system \eqref{eq:thPformSystemred} is equivalent to the ODE system found by Tracy and Widom for the Airy$_2$ process \cite{Tracy-Widom03}. 

\bigskip

We also derive a scalar PDE of Adler-van Moerbeke  \cite{Adler-van_Moerbeke05} when $m=2$ for strongly cubic integrable operators (see Subsection~\ref{sec:AVPDE} and Section~\ref{sec:symbolic}), and find a connection with the multi-component KP hierarchy (see Subsection~\ref{sec:KPhierarchy}).

\subsection{Fredholm determinants of cubic integrable operators}\label{sec:Freddetcio}

When the operator $\genH$ is trace class, 
the Fredholm determinant of $\genH$ is related to the matrix $\itZ_1$ in \eqref{defn:itYb1}, and thus to the differential equations of the last subsection. 

\begin{definition}[Cubic admissible determinant] \label{def:admissible} 
For a (strongly) cubic integrable operator $\genH$ that is trace class, we call its Fredholm determinant 
\beq 
	\genFunc(\genx,\geny,\gent)\deff\det(\itid-  \itH ) 
\eeq 
a {\it (strongly) cubic admissible} determinant.
\end{definition}

Recall the definition \eqref{eq:Eidefn} of $\itE_i$. 

\begin{prop} \label{lem:basic2}
Let $\itH$ be a cubic admissible operator. 
Let $\genFunc(\genx,\geny,\gent)=\det(\itid-  \itH )$ be the associated cubic admissible determinant and $\itZ_1$ be the matrix from \eqref{defn:itYb1}. Then, 
\beqq
	\partial_{\ith_i} \log \genFunc(\genx,\geny,\gent) =- \Tr \left( \itZ_1 \itE_{i} \right) , \qquad i\in \{ 1, \cdots, m\}. 
\eeqq
\end{prop}

\medskip

Recalling Definition~\ref{defn:Phi1d} and Lemma~\ref{lem:basic}, this result implies that 
\beqq
	\partial_\ith \log \gendetD=-\Tr(\genq) = \Tr(\gens) \quad \text{and} \quad 
	\partial_\ith^2 \log \gendetD= \Tr(\genp\genr) = \Tr(\genr\genp), 
\eeqq
where $\partial_\ith=\sum \limits_{k\in \itss}\partial_{\ith_{k}}$ is as in \eqref{eq:diffoperitss}. 
Thus, for example, we find from Theorem \ref{thm:KPeqtion} that if we take $\itss=\{1, \cdots, m\}$, then $\gens= \partial_\ith \log \gendetD$ and $\itvss=\partial_\ith^2 \log \gendetD$ satisfy the matrix KP-II equation \eqref{eq:kpv}, see Subsection~\ref{sec:PDEfordetad}.

In \eqref{eq:pjinfg} and Proposition~\ref{prop:additionaldeformation} below, we also find the derivatives of $\log \gendetD$ with respect to $\itl$ and $\itt$.

\subsection{The (periodic) KPZ fixed point and cubic integrable operators}\label{sec:DEKPZresult}

We are now discuss how the multi-point distributions of the KPZ and the periodic KPZ fixed points with the narrow wedge initial condition are related to cubic admissible determinants. 

Johansson and Rahman  \cite{Johansson-Rahman19} and also independently Liu \cite{Liu19} obtained explicit formulas for the multi-point distributions of the KPZ fixed point. 
These results were preceded by two-time distribution formulas by Johansson  \cites{JohanssonTwoTime, Johansson18}. 
The formula of Liu is the following. 
Let $\pkga=(\pkga_1,\cdots,\pkga_m)$ be position parameters, $\kpt=(\kpt_1,\cdots,\kpt_m)$ time parameters, and $\pkh=(\pkh_1,\cdots,\pkh_m)$ height parameters.  
Arrange the time parameters so that $\kpt_1\le \cdots\le \kpt_m$ and arrange the location parameters further so that $\kpga_i<\kpga_{i+1}$ whenever $\kpt_i=\kpt_{i+1}$. 
Then, \cite[Theorem 2.20 and Definition 2.23]{Liu19} states that the multi-point distribution in \eqref{eq:Fstepdefnintro} for the KPZ fixed point with the narrow wedge initial condition is given by 
\begin{equation} 
\label{eq:Fstepdefn}
\kpzFdist(\kph, \kpga, \kpt)
	= \frac1{(2\pi \ii)^{m-1}}  \oint \cdots \oint \itD (\kph,\kpga,\kpt\mid\oz) 
	\prod_{i=1}^{m-1} \frac{\dd \oz_{i}}{(1-\oz_i)\oz_{i}}. 
\end{equation}
Here, the integrals are over disjoint circles centered at the origin and of radii smaller than $1$, and the term $\itD (\kph,\kpga,\kpt\mid\oz) $ is a Fredholm determinant that depends on the parameters $\kph,\kpga$ and $\kpt$, and also on $\oz\deff (\oz_1, \cdots, \oz_{m-1})$.
We will recall its precise definition in Subsection~\ref{sec:Liuformula}.

In \cite{Baik-Liu19}, the authors computed the so-called relaxation time limit of the multi-point distributions of the periodic totally asymmetric simple 
exclusion process, extending their earlier work \cite{Baik-Liu16} on the one-point distribution.
The limit obtained in \cite[Theorem 2.1, Definition 2.5]{Baik-Liu19} is expected to be the multi-point distribution of the conjectured periodic KPZ fixed point. 
Due to the spatial periodicity, we may assume that the location variables satisfy $\pkga=(\pkga_1, \cdots, \pkga_m)^T\in [0,1)^m$.
As before, arrange the parameters so that $\pkt_1\leq \cdots\leq \pkt_m$. 
If $\kpt_j=\kpt_{j+1}$, we assume that $\pkh_j<\pkh_{j+1}$. 
Then, for the step initial condition, it was shown that 
\begin{equation} \label{eq:PKPZFstepdefn}
	\pkFdist(\pkh, \pkga, \pkt)
	\deff \frac1{(2\pi \ii)^{m}}  \oint \cdots \oint  \pkCdisc(\oz) \pkDdisc (\pkh,\pkt,\pkga\mid \oz)
	\prod_{i=1}^{m} \frac{\dd \oz_{i}}{\oz_{i}} 
\end{equation}
where the integration contours are circles about the origin satisfying $0<|\zeta_1|<\cdots<|\zeta_m|<1$, the function $\pkCdisc(\oz)$ is simple and explicit, and $\pkDdisc (\pkh,\pkt,\pkga\mid \oz) $ is again a Fredholm determinant. 
We discuss the precise formula in Subsection~\ref{sec:BLformula}. 

One main difference between $\itD$ and $\pkDdisc$ is that the operator for the former acts on an $L^2$ space of contours while the one for the latter acts on an $\ell^2$ space of a discrete set. 
The underlying operators, as originally obtained, share some similarities but the formulas are somewhat involved. 
It turns out that we can recast them to cubic admissible determinants, as we state as our next result.

\begin{thm}\label{thm:genintoperred}\hfill
\begin{enumerate}[(i)] 
\item (KPZ fixed point) For each $\zeta=(\zeta_1, \cdots, \zeta_{m-1})\in \C^{m-1}$ satisfying $0<|\zeta_1|, \cdots, |\zeta_{m-1}|<1$, the function  $(\ith,\itl,\itt)\mapsto \itD(\kph,\kpga,\kpt\mid \oz)$ is a \underline{strongly} cubic admissible determinant with the identification
\begin{equation}\label{def:xyttohgammataukpz}
	\itt_i = - \kpt_i/3, \qquad \itl_i=\kpga_i, \qquad \ith_i= \kph_i.
\end{equation}
Furthermore, the associated matrix $\itYb_1^\kpz=\itZ_1$ satisfies the symmetry relation
\begin{equation}\label{eq:symmetryprkpz}
\genL\itYb_1^\kpz(\ith,\itl,\itt) \genL^{-1}=\itYb_1^\kpz(\ith,-\itl,\itt)^T
\end{equation}
where  $\genL\deff\diag(\genL_1,\cdots, \genL_m,1)$ with $\genL_1= -(1-\zeta_1)$ and for $j\ge 2$, 
\begin{equation}\label{eq:symmetryLkpz} \begin{split}
	\genL_j\deff   (-1)^j \frac{(1-\zeta_1)(1-\zeta_2)\cdots (1-\zeta_j)}{(1-\frac1{\zeta_1}) (1-\frac1{\zeta_2}) \cdots (1-\frac1{\zeta_{j-1}})} . 
\end{split} \end{equation}

\item (Periodic KPZ fixed point) For each $\zeta=(\zeta_1, \cdots, \zeta_{m})\in \C^{m}$ satisfying  $0<|\zeta_1|<\cdots<|\zeta_m|<1$, the function $(\ith,\itl,\itt)\mapsto \pkDdisc(\pkh,\pkga,\pkt\mid \oz)$ is  a cubic admissible determinant with the identification
\begin{equation}\label{def:xyttohgammatau}
	\itt_i = - \pkt_i/3, \qquad \itl_i= \pkga_i/2, \qquad \ith_i= \pkh_i.
\end{equation}
\end{enumerate}
\end{thm}

We will describe the cubic integrable operators explicitly in Section~\ref{sec:cubicKPZ}.
This theorem implies that Proposition~\ref{lem:basic2}, Lemma~\ref{lem:basic}, Theorem \ref{thm:KdVheat}, and Theorem \ref{thm:KPeqtion} apply to (the cubic admissible operators $\genH^\kpz$ and $\genH^\per$ for) $\itD$ and $\pkDdisc$. 
The operator $\genH^\kpz$ also satisfies Theorem \ref{thm:matrixODE} since it is strongly admissible. 
These differential equation results assume that the Fredholm determinant is non-zero. 
Since the Fredholm determinant for the (periodic) KPZ fixed point is analytic in the parameters, this condition holds in an open set of the space of parameters.
For the case of the equal-time distributions of the KPZ fixed point, Quastel and Remenik \cite{Quastel-Remenik19b} obtained the same matrix KP-II equation as in Theorem \ref{thm:KPeqtion} although the quantities that solve the equation are different than ours\footnote{The paper also obtained the result for general initial conditions.}. 
When $m=1$, Theorems~\ref{thm:KdVheat}, \ref{thm:KPeqtion} and \ref{thm:KPeqtion} for 
the periodic KPZ fixed point and the KPZ fixed point were obtained in \cite{Baik-Liu-Silva}\footnote{In the paper \cite{Baik-Liu-Silva}, the roles of $\itrss$ and $\itpss$ are switched.}.  

The equal-time slice of the KPZ fixed point is a simple change of the Airy$_2$ process. 
However, due to the integrals in the formula \eqref{eq:Fstepdefn}, the function $\itD$ is not directly related to the multi-point distribution of the Airy$_2$ process. 
Nonetheless, a different Fredholm determinant formula for the Airy$_2$ process is known and it can also be turned into the cubic admissible determinant formalism. 
Thus, the matrix NLS system with complex time in Theorem \ref{thm:KdVheat} will hold for the Airy$_2$ process. 
This discussion will be carried out in a separate paper. 

For the periodic KPZ fixed point, when $m=1$ it was shown in \cite[Theorem~1.2]{Baik-Liu-Silva} that 
$$
	\begin{pmatrix} -1 & 0 \\ 0 & 1 \end{pmatrix} 
	\itYb_1^\per (\ith,\itl,\itt) 
	\begin{pmatrix} -1 & 0 \\ 0 & 1 \end{pmatrix} 
	= \itYb_1^\per (\ith,-\itl,\itt)^T. 
$$
However, it is not clear if a symmetry relation holds for $m>1$. 

\bigskip

Proposition \ref{lem:basic2} relates log derivatives of $\itD$ and $\pkDdisc$ with $\itZ_1$ which, in turn, satisfies differential equations. 
Integrating the log derivatives, we obtain the next formula.

\begin{prop}\label{thm:kpformula}
Set
\beq \label{deff:ithR}
	\veca\deff (1, 2,\cdots, m). 
\eeq
Let $\gendetD$ denote either $\itD$ or $\pkDdisc$. 
Assume that the parameters $(\ith,\itl,\itt)\in \R^{3m}$ satisfy $\gendetD(\ith+ \xi \veca,\itl,\itt)\neq 0$ for every $\xi\ge 0$.  
Then, 
\begin{equation}\label{eq:KPTWformula}
	\gendetD( \ith,  \itl, \itt)=\exp\left[    -\int_{0}^\infty \sum_{i=1}^m\Tr\left( \itqsn_i( \ith+ \xi \veca,  \itl, \itt )  \right)\dd \xi \right]
	=\exp\left[   \int_{0}^\infty \sum_{i=1}^m\Tr\left( \itssn_i( \ith+ \xi \veca,  \itl, \itt )  \right)\dd \xi \right]
\end{equation}
where $\itqsn_i$ and $\itssn_i$ denote the matrices $\itqsn$ and $\itssn$ in \eqref{eq:L1O1not} with the choice of the index set $\itss=\itss_i\deff \{1,\hdots,i\}$ for each $i=1,\cdots,m$.
\end{prop}

\subsection{Organization of the rest of the paper} \label{sec:organization}

In Section~\ref{sec:previousdes} we discuss a few special cases of the differential equations and compare with previous works in the literature. 
The proofs of the differential equations are given in Section~\ref{sec:IIKSRHP} using a connection between integrable operators and Riemann-Hilbert problems, applied to cubic integrable operators.  
In Section~\ref{sec:aspectsFreddet}, we introduce operators of a certain structure, which is not restricted to but does arise in the multi-point distributions for the KPZ and periodic KPZ fixed points. We then show that the corresponding Fredholm determinants can be recast as cubic admissible determinants. A special case of this result recasts Theorem~\ref{thm:genintoperred} for the KPZ and periodic KPZ fixed points, and Section~\ref{sec:cubicKPZ}  is devoted to this application. 
The remaining two sections are short.  
Section~\ref{sec:TWformulas} is about simple asymptotic properties of cubic integrable operators and a proof of Proposition~\ref{thm:kpformula}.
In Section~\ref{sec:symbolic}, we discuss  scalar PDEs that extend a work of Adler and van Moerbeke, which we first state in Subsection~\ref{sec:AVPDE}.

\section{Discussions on the differential equations and comparison with existing works} \label{sec:previousdes}

We compare the differential equations in Section~\ref{sec:DEsection} with various results in the literature.

\subsection{The case $m=1$: self-similar solutions}

When $m=1$, the equations \eqref{eq:heat}, \eqref{eq:mkdv}, \eqref{eq:kpu} and \eqref{eq:thPformSystemred} become equations for scalar-valued functions. We show that they admit self-similar solutions constructed out of Painlevé II transcendents.

Suppose that $w(\xi)$ is a single-variable function that satisfies the Painlevé II equation \eqref{eq:FstepTW2}
\beqq
	w''=\xi w+2w^3.  
\eeqq
Then, it is straightforward to check that for any non-zero constant $\alpha$, the functions
\beqq
	\itpsn (\itt, \itl, \ith) = \frac{\alpha}{(-3\itt)^{1/3}} \exp\left(-\frac{1}{3\itt}\ith\itl+\frac{2}{27\itt^2}\itl^3\right) w(\xi), \quad 
	\itrsn (\itt, \itl, \ith)= - \frac1{\alpha (-3\itt)^{1/3}} \exp\left(\frac{1}{3\itt}\ith\itl-\frac{2}{27\itt^2}\itl^3\right) w(\xi) 
\eeqq
with 
\beqq
	\xi=\frac{\ith}{(-3\itt)^{1/3}}+\frac{\itl^2}{(-3\itt)^{4/3}}
\eeqq
satisfy the NLS system with complex time \eqref{eq:heat} and the mKdV system \eqref{eq:mkdv}. 
They also satisfy the Tracy-Widom ODE system \eqref{eq:thPformSystemred}. 
It is also straightforward to check that the function
\beqq
 	\itqsn(\itt, \itl, \ith)=-(-3\itt)^{-1/3}Q(\xi) \quad \text{and} \quad 
	\genu(\itt, \itl, \ith)= (-3\itt)^{-2/3}w(\xi)^2
\eeqq
solve the scalar KP equation \eqref{eq:kpu} (observe that the commutator drops out $[ \genu, \partial_{\itl} \itqsn ]=0$, since the functions are scalar-valued). 
The functions $\itssn(\itt, \itl, \ith)=- \itqsn(\itt, \itl, \ith)$ and $\genv(\itt, \itl, \ith)=\genu(\itt, \itl, \ith)$ also solve the equation \eqref{eq:kpv}. 
Thus, these equations have solutions constructed out of a Painlevé II transcendent. . 

For the KPZ fixed point, we can check that the functions $\itpsn, \itrsn, \itqsn, \itssn$ for $m=1$ are indeed given by above forms with $w$ chosen to be the Hastings-McLeod solution to the Painlevé II equation. 
However, the functions for the periodic KPZ fixed point are not expected to be given by the above forms associated to the Painlev\'e II equation due to the fact that the space for the operator is discrete, which then implies that the relevant Riemann-Hilbert problem is of discrete type instead of of continuous type (see Section~\ref{sec:IIKSRHP}), the solutions are associated to solitons, and no self-similarity can be recast.

\subsection{The case of equal times: reduction of KP-II}

If we take $\itss=\{1, \cdots, m\}$, set $\kpt_1=\cdots=\kpt_m= - t/3$, and change the notations 
$\itl_i, \ith_i, \genu, \genq$ to $x_i, r_i, q, - Q$, then \eqref{eq:kpu} becomes
\beqq
	12\partial_t q+\partial_r^3q+6\partial_r(q^2)+ 3\partial_x^2 Q+  6[q,\partial_x Q] =0,\qquad \partial_r Q=q,
\eeqq 
where $\partial_r=\sum_{i=1}^m \partial_{r_i}$ and $\partial_x= \sum_{i=1}^m \partial_{x_i}$. 
This is the same matrix KP equation obtained in \cite[(1.6)]{Quastel-Remenik19b} by Quastel and Remenik. 
However,  the functions $q, Q$ considered in \cite{Quastel-Remenik19b} are different from what we considered in this paper. 
Namely, the functions in \cite{Quastel-Remenik19b} are directly related to the multi-point distributions $\kpzFdist(\kph, \kpga, \kpt)$ while in this paper, the functions are related to $\itD (\kph,\kpga,\kpt\mid\oz)$ which is a part of the formula \eqref{eq:Fstepdefn} for the multi-point distribution. 
It is puzzling that they satisfy the same KP-II equation. 
If we compare the result of \cite{Quastel-Remenik19b} and Equation \eqref{eq:KPTWformula} with formula \eqref{eq:Fstepdefn} in mind, we obtain the identity 
\begin{equation}\label{eq:KPintKPtr}
\Tr(q(\kph,\kpga,t))=\partial^2_{\kph}\log\left[ \frac1{(2\pi \ii)^{m-1}}  \oint \cdots \oint \exp\left[ \int_{0}^\infty \sum_{i=1}^m\Tr\left( \itssn_i( \ith+ \xi \veca,  \itl, \itt )  \right)\dd \xi  \right]  
	\prod_{i=1}^{m-1} \frac{\dd \oz_{i}}{(1-\oz_i)\oz_{i}} \right]
\end{equation}
between a single solution $q$ to the KP equation on the left-hand side and a superposition on the right-hand side of families of solutions $\gens_i=\gens_i^\kpz$ to different KP-II equations corresponding to different permutations $\itss=\itss_i=\{1,\hdots, i\}$. 
It is unclear whether this identity holds because the KP-II solutions related to the KPZ fixed point are very special, or if such an identity holds for a broad class of solutions. 
We note that $q$ on the left-hand side is a real solution while $\gens_i$ are complex-valued solutions.

\subsection{The case of equal times: matrix ODE system}\label{sec:multipointairydiscussion}

Tracy and Widom obtained in  \cite{Tracy-Widom03} a matrix ODE system for the Airy$_2$ process, which is a simple transformation of the equal-time slice of the KPZ fixed point. 
We show here that the ODE system \eqref{eq:thPformSystemred} for the equal-time case of the cubic integrable operator is the same as the matrix ODE system of Tracy and Widom. 

The ODE system \eqref{eq:thPformSystemred} involves matrix functions $\genp,\genr,\genq$, and derivatives with respect to the variables $\ith_k$. We rename $\itt$ to $-\frac{\itt}{3}$ so that
$\itt_1=\cdots=\itt_m\revdeff -\frac{\itt}{3}$. 
The ODE system contains $\itt$ and $\itl_k$ as parameters. 
Change the variables $\ith_k$ to $\xi_k$ by 
$$
	\xi_k\deff \frac{\ith_k}{\itt^{1/3}}+\frac{\itl_k^2}{\itt^{4/3}}. 
$$
Set $\bs\xi\deff\diag(\xi_1,\cdots,\xi_m)=\frac{1}{\itt^{1/3}}\bs\ith+\frac{1}{\itt^{4/3}}\bs\itl^2$, where $\itlm\deff \diag(\itl_1,\hdots,\itl_m)$ and $\ithm\deff \diag(\ith_1,\hdots,\ith_m)$, and also $\xi\deff(\xi_1,\cdots,\xi_m)$, $\partial_\xi\deff \sum_{j=1}^m \partial_{\xi_j}$.
Define $\msf P$, $\msf R$, $\msf Q$ from $\itpsn$, $\itrsn$, $\itqsn$ by setting 
\begin{align*}
	& \itpsn (\itt, \itl, \ith) = \ii\itt^{-1/3}\exp\left(\frac{1}{\itt}\bs\ith\bs\itl+\frac{2}{3\itt^2}\bs\itl^3\right) \msf P(\xi), \\
	&\itrsn (\itt, \itl, \ith) = \ii\itt^{-1/3}\msf R(\xi)\exp\left(-\frac{1}{\itt}\bs\ith\bs\itl-\frac{2}{3\itt^2}\bs\itl^3\right), \\
& \itqsn (\itt, \itl, \ith) =-\ii \itt^{-1/3}\exp\left(\frac{1}{\itt}\bs\ith\bs\itl+\frac{2}{3\itt^2}\bs\itl^3\right) \msf Q(\xi)\exp\left(-\frac{1}{\itt}\bs\ith\bs\itl-\frac{2}{3\itt^2}\bs\itl^3\right), 
\end{align*}
Then, with $\wt{\bs\itl}\deff \itt^{-2/3}\bs\itl$, the system \eqref{eq:thPformSystemred} becomes
\begin{equation}\label{eq:TWsysteminourform}
\partial_\xi \msf Q=-\left[\wt{\bs\itl},\msf Q\right]-\msf P\msf R, \quad \partial_\xi^2\msf P=\bs\xi \msf P+2\msf P\msf R\msf P+2\left[\wt{\bs\itl},\msf Q\right]\msf P,\quad \partial_\xi^2\msf R=\msf R\bs\xi +2\msf R\msf P\msf R+2\msf R\left[\wt{\bs\itl},\msf Q\right].
\end{equation}

Equations (1)--(3) of \cite{Tracy-Widom03}\footnote{The published version of the system from \cite{Tracy-Widom03} has a typo which was later fixed in the ArXiv version (ArXiv:0302033v4, see Equations~(1)--(3) and the last two displayed equations in page 3 therein).  Our discussion here uses the latest Arxiv version.} form an ODE system for three $m\times m$ matrices $r, q, \tilde{q}$. 
For the $m$-dimensional vector 
$\ite\deff (1, \cdots, 1)^T$
consider the $m\times m$, $m\times 1$, and $1\times m$ matrices $r, q \ite$, and $\ite^T \tilde q$. 
Then, equations (1)--(3) of \cite{Tracy-Widom03} changed for these functions are precisely \eqref{eq:TWsysteminourform} with the identification 
\beqq
	\msf Q= r, \qquad \msf P= q\ite, \qquad \msf R= \ite^T\tilde q\qquad \text{and}\qquad \tau_k=\itt^{-2/3}\itl_k=\wt\itl_k.
\eeqq

\subsection{Differential equations related to the Adler and van Moerbeke PDE}\label{sec:AVPDE}

For the equal-time, 2-position distributions of the KPZ fixed point (i.e. the Airy$_2$ process), Adler and van Moerbeke also obtained a PDE \cite{Adler-van_Moerbeke05}. We discuss how this PDE arises from our analysis.

\subsubsection{A scalar PDE for cubic admissible determinants} \label{sec:PDEfordetad}

Assume that the cubic integrable operator $\genH$ is trace class and consider the log determinant
\beqq
	M \deff \log \genFunc(\genx,\geny,\gent)= \log \det(\itid-\itH).
\eeqq 
The KP-II equation in Theorem \ref{thm:KPeqtion}, which is for matrix functions $\itpsn, \itrsn, \itqsn, \itssn$, yields a scalar partial differential equation for $M$ in the following way. 
Choose the set $\itss= \{1,\hdots,m\}$ so that $\itsss$ and $\itvss=\itrss\itpss$ become scalar functions. 
Then, the commutator in  \eqref{eq:kpv} disappears, and from Proposition \ref{lem:basic2} and Lemma \ref{lem:basic}, we find that $\itDsnh M= \itsss$ and $\itDsnh^2 M= \itvss$. 
Thus, the equation \eqref{eq:kpv}  becomes
\begin{equation}\label{eq:MpreKP}
	\itDsnh(- 4 \itDsnh\itDsnt M + \itDsnh^4 M  + 6  ( (\itDsnh^2 M )^2) + 3  \itDsnl^2 M)=0.
\end{equation}
If $M$ decays fast enough as $\ith_i\to \infty$, which is the case for the (periodic) KPZ fixed points, integrating the above equation we obtain a fourth order nonlinear differential equation for $M$, 
\beq \label{eq:kpfordet}
	- 4 \itDsnh\itDsnt M + \itDsnh^4 M  + 6  ( (\itDsnh^2 M )^2) + 3  \itDsnl^2 M=0.
\eeq

\subsubsection{Invariance of strongly cubic admissible determinants}

Let $\Omega$ be a union of contours, and suppose that $\genH$ is a strongly cubic integrable operator that is trace class with parameters $\ith, \itl, \itt$. Let $\Omega'$ be a scaled and translated contour of $\Omega$ given by 
\beqq
	\Omega'= (-3\itt_1)^{\frac{1}{3}}   \Omega -  \frac{\itl_1}{(-3\itt_1)^{2/3}}. 
\eeqq
Using the fact that the functions $\Fm,\Gm$ are constants on each connected component of $ \Omega$ (since $\genH$ is strongly cubic integrable), it is straightforward to check that 
\beq \label{eq:Hstinv}
	\frac1{(-3\itt_1)^{\frac{1}{3}}} \genH \bigg( \frac1{ (-3\itt_1)^{ \frac{1}{3}}} u- \frac{\itl_1}{3\itt_1}, \frac1{ (-3\itt_1)^{ \frac{1}{3}}} v- \frac{\itl_1}{3\itt_1} \mid  \ith, \itl, \itt\bigg) 
	=  \genH ( u, v \mid x, y, t), 
	\qquad u, v\in \Omega', 
\eeq 
for the parameters $x=(x_1, \cdots, x_m)$, $y=(y_1, \cdots, y_m)$, $t=(t_1, \cdots, t_m)$ given by 
\beq \label{eq:pararel} \begin{split}
	&x_i= \frac{1}{(-3\itt_1)^{ \frac{1}{3}}}\ith_i+\frac{2\itl_1}{(-3\itt_1)^{ \frac{4}{3}}} \itl_i - \frac{\itl_1^2}{(-3\itt_1)^{ \frac{4}{3}}\itt_1} \itt_i, \qquad 
	y_i=\frac{1}{(-3\itt_1)^{ \frac{2}{3}}}\itl_i- \frac{\itl_1}{{(-3\itt_1)^{ \frac{2}{3}}}\itt_1} \itt_i, \qquad
	t_i= \frac{1}{-3\itt_1} \itt_i 
\end{split} \eeq
for $i=1, \cdots, m$. Note that for $i=1$, 
\beqq
	y_1=0, \qquad t_1=-\frac13.
\eeqq

The left-hand side of \eqref{eq:Hstinv} is a change of variables of the kernel  $\genH(u,v \mid \ith, \itl, \itt)$, $u,v \in \Omega$, and hence its Fredholm determinant is equal to $\det(\itid-\genH(  \ith, \itl, \itt))_{L^2(\Omega)}$. 
On the other hand, the Fredholm determinant of the right-hand side of \eqref{eq:Hstinv}, $\det(\itid-\genH ( x, y, t))_{L^2(\Omega')}$, is equal to $\det(\itid-\genH ( x, y, t))_{L^2(\Omega)}$ by the Cauchy's theorem using the fact that restrictions of $\genH ( u, v \mid x, y, t)$ for $u$ and $v$ on connected components of $\Omega'$ extends to entire functions in $u$ and $v$. Thus, we find a relationship
\beqq
	M= \det(\itid-\genH ( \ith, \itl, \itt)) = \det(\itid-\genH ( x,y,t)) .
\eeqq
Under the change of variables given by \eqref{eq:pararel}, the equation \eqref{eq:kpfordet} evaluated at $\itt_1=-1/3$ and $\itl_1=0$ becomes the following equation in $x, y, t$:
\begin{equation} \label{eq:kpfordetfin}
\begin{split}
    &
    12\sum_{i=2}^m{(1-\aitau_i)}{} \partial_{\aih}\partial_{\aitau_i}M
    -8\sum_{i=2}^m{}{}\ail_i\partial_{\aih}\partial_{\ail_i}M
    -4\sum_{i=1}^m\aitau_i^{-2}\left(\aih_i\aitau_i^2+3\ail_i^2\right)\partial_{\aih}\partial_{\aih_i}M
    -4\partial_{\aih}M +\partial_{\aih}^4M +6(\partial_{\aih}^2M)^2\\&
    +3\sum_{j=2}^m\sum_{i=2}^m(1-\aitau_i)(1-\aitau_j)\partial_{\ail_j}\partial_{\ail_i}M
    +12\sum_{j=2}^m\sum_{i=1}^m\ail_i\aitau_i^{-1}(1-\aitau_j)\partial_{\ail_j}\partial_{\aih_i}M \\
    &+12\sum_{i,j=1}^m\ail_i\ail_j\aitau_i^{-1}\aitau_j^{-1}\partial_{\aih_j}\partial_{\aih_i}M
    +6\sum_{i=1}^m\aitau_i^{-1}\partial_{\aih_i}M =0.
\end{split}
\end{equation}

Note that if $\itt_1=-1/3$ and $\itl_1=0$, then $(\ith, \itl, \itt)=(x,y,t)$. 
Thus, we find that $M(\ith, \itl, \itt)$ with $\itt_1=-1/3$ and $\itl_1=0$ satisfies the differential equation \eqref{eq:kpfordetfin} with respect to the original variables $(\ith, \itl, \itt)$. 

\subsubsection{Scalar PDEs for the two-time and one-time cases - Adler and van Moerbeke PDE}

Let $\genH$ be a strongly cubic integrable operator that is trace class, and consider the case when $m=2$ and $\itt_1=-1/3$ and $\itl_1=0$. 
Then, $M=\det(\itid-\genH)$ satisfies the equation \eqref{eq:kpfordetfin} with $m=2$ with respect to the variables $(\ith, \itl, \itt)$. 
Motivated by the work of \cite{Adler-van_Moerbeke05} on the Airy$_2$ process, introduce the change of variables 
\beq \label{eq:twoptpara}
	\itt_2= -\frac{\aitau}{3},\qquad  \ith_1=\frac{\aiE+\aiW}{2},\qquad\ith_2=\frac{\aiE-\aiW}{2}-\frac{\ail^2}{\aitau},  \qquad   \itl_2 =\ail. 
\eeq
Then \eqref{eq:kpfordetfin} with $m=2$ becomes
\begin{equation} \label{eq:maple5}
\begin{split}
    &
    24{(1-\aitau)}{}\aitau^2 \partial_{\aiE}\partial_{\aitau}M+16\aitau^2\partial_{\aiE}^4M+96\aitau^2(\partial_{\aiE}^2M)^2
    +3(1-\aitau)^2\aitau^2\partial_{\ail}^2M+4\ail\aitau(3-7\aitau)\partial_{\aiE}\partial_{\ail}M\\&+12\ail\aitau(\aitau-1)\partial_{\ail}\partial_{\aiW}M-4(2\aiE\aitau^2+3\ail^2)\partial_{\aiE}^2M-8\aiW\aitau^2\partial_\aiE\partial_\aiW M+12\ail^2\partial_{\aiW}^2M+2\aitau(3-\aitau)\partial_\aiE M\\&+6\aitau(\aitau-1)\partial_\aiW M=0.
\end{split}
\end{equation}
In the further special case when $t=1$ so that $\itt_1=\itt_2$, the above equation becomes simpler. If we consider 
$\frac14 \partial_\aiE   \eqref{eq:maple5}$, then we arrive at 
\beq \label{eq:QR}
	\left( -4\partial_{\aiE}^5+2{\aiW}\partial_{\aiE}^2\partial_{\aiW}+(3\ail^2+2{\aiE})\partial_{\aiE}^3-3\ail^2\partial_{\aiE}\partial_{\aiW}^2-\partial_{\aiE}^2 +4\ail\partial_{\aiE}^2\partial_{\ail}\right )  M
	-48\partial_{\aiE}^2M\partial_{\aiE}^3M=0.
\eeq
This PDE is the same as the one derived by Quastel and  Remenik \cite[Theorem~2.5]{Quastel-Remenik19b} for the $2$-point distribution of the KPZ fixed point with equal time\footnote{The correspondence of terms between \cite[Theorem~2.5]{Quastel-Remenik19b} and here is $\psi=\partial_{\aih}M$, $\bm\Psi=M$, $r=x$, and the variables $\ail$ coincide. }, except for a flipped sign in the term $\partial^2_\aiE$.

\medskip

In Section~\ref{sec:symbolic}, we derive the equation \eqref{eq:maple5} in a different way using arguments on the underlying Lax pairs, supported by symbolic calculations in linear algebra. 
As we also indicate therein, such arguments allow us to obtain three additional PDEs. One of them is 
\begin{equation} \label{eq:maple6}
 \begin{split}
 &3(\aitau-1)\partial_\aiE^4M-8(\aitau+1)\partial_\aiE^3\partial_\aiW M+6(\aitau-1)\partial_\aiE^2\partial_\aiW^2M+(1-\aitau)\partial_\aiW^4M+12\aiW\partial_\aiE^2M+8\aiE \partial_\aiE\partial_\aiW M\\&
 -4\aiW\partial_\aiW^2M+3\aitau(1-\aitau)\partial_\ail^2M-8\ail\partial_{\aiW}\partial_\ail M-4\partial_\aiW M+18(\aitau-1)(\partial_{\aiE}^2M)^2+24(\aitau-1)(\partial_\aiE\partial_\aiW M)^2\\&
 +6(1-\aitau)(\partial_{\aiW}^2M)^2+12(\aitau-1)\partial_\aiE^2 M\partial_\aiW^2 M-48(\aitau+1)\partial_\aiE^2 M\partial_\aiE\partial_\aiW M=0.
 \end{split}  
 \end{equation}
 If we set $t=1$ in this equation, the equation again simplifies, and $\frac1{12}\partial_\aiW   \eqref{eq:maple5} +\frac1{12} \partial_\aiE  \eqref{eq:maple6}$ becomes  
 \beq \label{eq:AV}
	\left( \ail^2(\partial_{\aiW}^3-\partial_{\aiE}^2\partial_{\aiW})-{\aiW}(\partial_{\aiE}\partial_{\aiW}^2-\partial_{\aiE}^3 )-2\ail \partial_{\aiE}\partial_{\aiW}\partial_{\ail }\right) M 
	-8(\partial_{\aiE}\partial_{\aiW}M)(\partial_{\aiE}^3M)+8(\partial_{\aiE}^2M)(\partial_{\aiE}^2\partial_{\aiW}M)=0. 
\eeq
This is the same PDE obtained by Adler and van Moerbeke \cite[Corollary~1.3]{Adler-van_Moerbeke05} for the $2$-point distribution of the Airy$_2$ process. 
Thus, Equation \eqref{eq:maple6} is a non-equal time extension of the Adler-van Moerbeke PDE. 
We note that the equation \eqref{eq:AV} was also re-derived in \cite{Bertola-Cafasso-2012b} using a Riemann-Hilbert approach, and the Adler-van Moerbeke PDE was extended to general $m>2$ (single time, multi-point) in \cite{Wang09a}.


\bigskip


As we just discussed, the multi-time multi-location differential equations we obtained reduce in the one-time multi-location case to several different equations obtained previously. We stress that although the differential equations admit such reductions, the solutions themselves are not directly comparable: all the previous works we mentioned obtain equations for the (log derivatives of the) distributions themselves, whereas in our case the solutions to the equations enter into the distribution through the integration in the $\zeta$-variables in \eqref{eq:Fstepdefn} (recall for instance the discussion on \eqref{eq:KPintKPtr}). It is interesting that the underlying equations are nevertheless the same.

\section{Derivation of differential equations}\label{sec:IIKSRHP}

In this section we prove the results stated in Subsection~\ref{sec:DEsection} and~\ref{sec:Freddetcio}. In order to derive differential equations, we use a known connection between IIKS integrable operators and Riemann-Hilbert problems (RHPs), and then exploit the structure of the Riemann-Hilbert problem for the case of cubic integrable operators. 
We follow the general methodology, often called the dressing method, which derives so-called Lax equations from a Riemann-Hilbert problem when the jump/residue matrix does not depend on parameters. 
The Lax equations are then combined appropriately to derive our claimed differential equations.

\subsection{A review of integrable operators and RHP}\label{sec:IIKSdiscussion}

An operator $\genH$ of the form \eqref{eq:admisskernel} for general vectors $\genf$ and $\geng$ is called an (IIKS) integrable operator, after its introduction by Its, Izergin, Korepin and Slavnov in \cite{IIKS}. 
Their general theory is discussed in \cites{IIKS, Deift99b} for integrable operators acting on continuous contours and \cite{Borodin00} for those on discrete sets. In this subsection we review some of their properties that we will use.

Let $\Omega$ be a union of contours or a discrete set as described in Subsection \ref{sec:DEsection},  and let $\mu$ and $\genHilb$ be as in \eqref{eq:mudeff} and \eqref{deff:genHilb}.

Let $\genf$ and $\geng$ be $(m+1)$-dimensional column vector-valued functions (not necessarily of the form in Definition \ref{def:cubicintegrabl}) in $L^2(\Omega,\mu)\cap L^\infty(\Omega,\mu)$ satisfying 
$$
	\genf(u)^T\geng(u)=0, \; u\in\Omega,\qquad \text{and set}\qquad 
\genH(u,v)\deff \frac{\genf(u)^T \geng(v)}{u-v}, \; u,v\in \Omega.
$$ 
Observe that for $\genf$ and $\geng$ of the form in Definition \ref{def:cubicintegrabl}, the assumption \eqref{eq:withEzero} implies this orthogonality. 
%
Assume that $\itid-\genH$ is invertible. With $\itI_k$ being the identity matrix of size $k$, define $(m+1)$-dimensional column vectors 
\begin{equation}\label{eq:ResolventKernel}
	\genF\deff (\itid-\genH)^{-1}\genf, \qquad \genG\deff (\itid-\genH^T)^{-1}\geng,
\end{equation}
and the $(m+1)\times (m+1)$ matrix-valued function
\beq \label{eq:phiinFgfGfm}
	\itYb(z)=\itI_{m+1}-\int_{\Omega} \frac{\genF(u)\geng(u)^T}{u-z}\dd\mu(u), \qquad z\in \C\setminus\Omega. 
\eeq
An important property of integrable operators is that $\itYb$ is the unique solution to the following Riemann-Hilbert problem (RHP). 
\begin{enumerate}[label={\bfseries RHP (\alph*):}]
\item $\itYb: \mb C\setminus \Omega \to\mb C^{{(m+1)}\times {(m+1)}}$ is analytic. 
\item In the case when $\Omega$ is a union of contours, $\itYb$ satisfies the jump condition 
\beq \label{eq:Jcon}
	\itYb_+(u)=\itYb_-(u) \itJ(u), \quad u\in \Omega, \qquad \text{where}\quad\itJ(u) \deff\itI_{m+1}-2\pi \ii \genf(u) \geng(u)^T .
\eeq 
In the case when $\Omega$ is a discrete set, the matrix $\itYb$ has simple poles at the points $u\in \Omega$, and its residues satisfy
\beq\label{eq:Jdis}
	\res_{z=u} \itYb(z) = \lim_{z\to u} \itYb(z) \itR(u) \qquad \text{with}\qquad \itR(u)\deff \genf(u) \geng(u)^T.
\eeq
The convergence of the limit is a part of the condition. 

\item $\itYb$ has an asymptotic series 
\beq \label{eq:expYk}
	\itYb(z) \sim \itI_{m+1} + \sum_{n=1}^\infty \frac{\itYb_n}{z^n}
\eeq
as $z\to\infty$ uniformly away from $\Omega$. 
\end{enumerate}

It can be shown that $\det \itYb(z) =1$ for all $z$ and 
\beq \label{eq:phiinFgfGfm2}
	\itYb(z)^{-1}=\itI_{m+1}+\int_{\Omega} \frac{\genf(u)\genG(u)^T}{u-z}\dd\mu(u) , \qquad z\in \C\setminus\Omega. 
\eeq
Note that \eqref{eq:phiinFgfGfm} and \eqref{eq:phiinFgfGfm2} imply that $\itYb(z) \sim \itI_{m+1} + \sum_{n=1}^\infty \frac{\itYb_n}{z^n}$ and 
$\itYb(z)^{-1}\sim \itI_{m+1}+ \sum_{n=1}^\infty \frac{(\itYb^{-1})_n}{z^n}$ as $z\to \infty$ where 
\beq \label{eq:pjinfg}
	\itYb_n= \int_{\Omega}  u^{n-1} \itF(u) \geng(u)^T \dd \mu(u)  , 
	\qquad
	(\itYb^{-1})_n = - \int_{\Omega} u^{n-1} \genf(u)  \itG(u)^T  \dd \mu(u) , \qquad n=1,2, \cdots. 
\eeq
The notation $\itYb_1$ is consistent with Definition~\ref{defn:Phi1d}.

\subsection{Lax equations}\label{sec:Lax}
The discussion in the last subsection is valid for general integrable operators. For the rest of Section~\ref{sec:IIKSRHP}, we further assume that $\genH$ is a cubic integrable operator as introduced in Definition~\ref{def:cubicintegrabl}. 

From \eqref{eq:fg}, $\genf=\genc\Delta \Fm$ and $\geng=\genc^{-1}\Delta^{-1}\Gm$ for vector functions $\Fm$ and $\Gm$ that do not depend on $\genx, \geny, \gent$.
Hence, the jump matrix and the residue matrix for the RHP in 
\eqref{eq:Jcon} and \eqref{eq:Jdis} take the form 
\beq \label{eq:rhpjupp}
	\itJ(u)=\itconjb(u)  \itJ_0(u) \itconjb(u)^{-1}
	\qquad \text{and}\qquad
	\itR(u)= \itconjb(u)\itR_0(u) \itconjb(u)^{-1},
\eeq 
where the matrices $\itJ_0(u) = \itI_{m+1}-2\pi \ii \Fm(u) \Gm(u)^T$ and $\itR_0(u)= \Fm(u) \Gm(u)^T$ 
do not depend on $\genx, \geny, \gent$. 
We follow the standard dressing method of deriving Lax equations for RHPs when the jump or residue matrix is a conjugation of a parameter-independent matrix. We think of $\itJ$ and $\itR$ as dressed-up versions of the matrices $\itJ_0$ and $\itR_0$ which are constant in the parameters.

Recall that $\itE_i$ denotes the $(m+1)\times (m+1)$ matrix whose  $(i,i)$-entry is $1$ and all other entries are $0$. The linear ODEs \eqref{eq:laxeq} below are called Lax equations for $\itZsn$. 

\begin{lem}[Lax equations]
Let $\itZsn(z)$ be the solution, which is given by \eqref{eq:phiinFgfGfm}, of the Riemann-Hilbert problem associated to the cubic integrable operator $\genH$. Then, the matrix 
\beq \label{eq:itWsndef}
	\itWsn(z) \deff \itZsn(z) \itconjb (z)
\eeq
satisfies 
\begin{equation} \label{eq:laxeq} \begin{split}
	\partial_{\genx_i}  \itWsn(z) = \laxD_1^{(i)}(z) \itWsn(z), \qquad
	 \partial_{\geny_i}  \itWsn(z) = \laxD_2^{(i)}(z)  \itWsn(z) , \qquad
	 \partial_{\gent_i}  \itWsn(z) = \laxD_3^{(i)}(z)  \itWsn(z) , \\
\end{split}
\end{equation}
for $i\in \{1, \cdots, m\}$, where
\beq \label{eq:Bdefn}
	\msf D_1^{(i)}(z)=z \itE_i+\laxC_{1}^{(i)} , \qquad 
	\msf D_2^{(i)}(z)=z^2 \itE_i+ z \laxC_{1}^{(i)} +\laxC_{2}^{(i)}, \qquad 
	\msf D_3^{(i)}(z)= z^3 \itE_i+ z^2 \laxC_{1}^{(i)} + z \laxC_{2}^{(i)} +\laxC_{3}^{(i)},  
\eeq
with 
\beq \label{eq:Aiinmulhi}
\begin{split} 
	\laxC_1^{(i)}  \deff [\itZ_1, \itE_i], \qquad
	\laxC_2^{(i)}  \deff [\itZ_2, \itE_i] - \laxC_1^{(i)} \itZ_{1}, \qquad 
	\laxC_3^{(i)}  \deff [\itZ_3, \itE_i] - \laxC_1^{(i)} \itZ_{2} - \laxC_2^{(i)} \itZ_{1} . 
\end{split}	
\eeq
\end{lem}

\begin{proof}
From \eqref{eq:rhpjupp}, the jump/residue condition for the new matrix becomes 
$$
\itWsn_+(u)= \itWsn_-(u)\itJ_0(u)\quad \text{or}\quad \res_{z=u} \itWsn(z) = \lim_{z\to u} \itWsn(z) \itR_0(u).
$$ 
Since $\itJ_0(u)$ and $\itR_0(u)$ do not depend on $\genx, \geny, \gent$, we find that a partial derivative $\partial \itWsn$ with respect any one of the parameters satisfies the same jump/residue condition as $\itWsn$ itself. 
Hence, $(\partial \itWsn (z))\itWsn(z)^{-1}$ is continuous across $\Omega$, and thus it is an entire matrix function. 
From the definition \eqref{eq:itWsndef}, 
\begin{equation}\label{eq:itWsndiffeqgen}
	(\partial \itWsn )\itWsn^{-1} = (\partial \itZsn )\itZsn^{-1} + \itZsn (\partial \itconjb ) \itconjb^{-1} \itZsn^{-1}.
\end{equation}
Also, from the formula of $\itconjb$ in \eqref{eq:itconjb}, 
\beq \label{eq:deltadde}
	\partial_{\genx_i} \itconjb (z) = z \itE_i \itconjb(z), \qquad
	\partial_{\geny_i} \itconjb (z) = z^2 \itE_i \itconjb(z), \qquad
	\partial_{\gent_i} \itconjb (z) = z^3 \itE_i \itconjb(z). 
\eeq
Hence, using the fact that $\itZsn \sim \itI_{m+1}$ as $z\to \infty$, Liouville's Theorem implies that the matrix functions
$(\partial_{\genx_i}  \itWsn(z)) \itWsn(z)^{-1}$ , $(\partial_{\geny_i}  \itWsn(z)) \itWsn(z)^{-1}$, and $(\partial_{\gent_i}  \itWsn(z)) \itWsn(z)^{-1}$ 
are polynomials in $z$ of degrees $1$, $2$, and $3$, respectively. 
Inserting \eqref{eq:expYk} into \eqref{eq:itWsndiffeqgen}, we obtain the coefficients of these polynomials in terms of $\itZsn_1, \itZsn_2$ and $\itZsn_3$, and find \eqref{eq:laxeq}. 
\end{proof}

\subsection{NLS and mKdV: proof of Lemma~\ref{lem:basic} and Theorem~\ref{thm:KdVheat}}\label{sec:proofcoupledsystem}

Recall that in Definition~\ref{defn:Phi1d}  $\itss\subset \{1,\hdots,m\}$ is a given set, and the associated differential operators $\itDsnt$, $\itDsnl$, and $\itDsnh$ are given in \eqref{eq:diffoperitss} such as $\partial_{\genx}= \sum_{i\in \itss} \partial_{\genx_i}$.  
Set
\begin{equation}\label{def:itEnitEss}
	\itEn = \itEss\deff \sum_{i\in \itss} \itE_i . 
\end{equation}

We first derive the following differential equations for  $\itZsn_1, \itZsn_2, \itZsn_3$  (recall \eqref{eq:pjinfg}).

\begin{prop} \label{lem:DEforphi1}
The matrix $\itZsn_1$ satisfies the differential equations
\beq \label{eq:kpzmateqdf} \begin{split}
	&   \itDsnh^2 \itZsn_1 + [\itDsnl \itZsn_1, \itEn] +  [\itDsnh \itZsn_1, [\itZsn_1, \itEn]] =0, \\
	&   \itDsnh \itDsnl \itZsn_1 +[\itDsnt \itZsn_1, \itEn] +   [\itDsnl \itZsn_1, [\itZsn_1, \itEn]] =0, \\
	&    \itDsnl^2 \itZsn_1 - \itDsnt \itDsnh \itZsn_1 - [\itDsnl \itZsn_1, \itDsnh \itZsn_1] =0. 
\end{split} \eeq
Furthermore, the three matrices $\itZsn_1, \itZsn_2, \itZsn_3$ are related by
\beq \label{eq:kpzmateqdf123} \begin{split}
	&\itDsnh \itZsn_1= {- [\itZsn_{2}, \itEsn]+[\itZsn_1, \itEn] \itZsn_{1}}, \\
	&\itDsnl \itZsn_1= {- [\itZsn_{3}, \itEsn]+[\itZsn_1, \itEn]\itZsn_{2}  -(  \itDsnh \itZsn_1) \itZsn_{1}} .
\end{split} \eeq
\end{prop}

\begin{proof}
The matrix $\itWsn(z) \deff \itZsn(z) \itconjb (z)$ of \eqref{eq:itWsndef} satisfies the Lax equations \eqref{eq:laxeq}. 
Summing over $i\in \itss$, 
\begin{equation} \label{eq:tyxeqforZ} \begin{split}
	& \itDsnh \itWsn(z) = ( z \itEn+  \itAsn_1 ) \itWsn(z), \\
	& \itDsnl \itWsn(z) = ( z^2 \itEn+  z \itAsn_1 + \itAsn_2 ) \itWsn(z) , \\
	& \itDsnt \itWsn(z) = ( z^3 \itEn+ z^2 \itAsn_1 + z \itAsn_2 + \itAsn_3 ) \itWsn(z) , \\
\end{split}
\end{equation}
with
\begin{equation} \label{eq:itAssdef1}
	\itAsn_1 \deff [\itZsn_1, \itEn], \qquad
	\itAsn_2 \deff [\itZsn_2, \itEn] - \itAsn_1 \itZsn_{1}, \qquad
	\itAsn_3 \deff [\itZsn_3, \itEn] - \itAsn_1 \itZsn_{2} - \itAsn_2 \itZsn_{1}. 
\end{equation}
Inserting $\itWsn(z)= \itZsn(z) \itconjb (z)$, these equations  become equations for $\itZsn(z)$. We then insert the asymptotic series \eqref{eq:expYk} and consider the coefficients of $\Boh(z^{-k})$. 
Collecting the coefficients, we obtain a sequence of equations relating the matrices $\itZsn_n$. 
In particular, the $\Boh(z^{-1})$ terms of the first two equations of \eqref{eq:tyxeqforZ} become  
\beq \label{eq:AiinZizm} 
	\itDsnh \itZsn_1= {- [\itZsn_{2}, \itEsn]+\itAsn_{1}\itZsn_{1}} ,\qquad 	
	\itDsnl \itZsn_1= {- [\itZsn_{3}, \itEsn]+\itAsn_{1}\itZsn_{2}+\itAsn_{2}\itZsn_{1}}.
\eeq
Using these relations, we find that the three matrices in \eqref{eq:itAssdef1} can all be expressed in terms of $\itZsn_1$ as
\beq \label{eq:eqfromlargez}
	\itAsn_1 = [\itZsn_1, \itEn], \qquad \itAsn_2= -\itDsnh \itZsn_1,\qquad \itAsn_3= -\itDsnl \itZsn_1.
\eeq
Inserting the formulas of $\itAsn_1$ and $\itAsn_2$ into \eqref{eq:AiinZizm}, we obtain \eqref{eq:kpzmateqdf123}.

To obtain the differential equations~\eqref{eq:kpzmateqdf}, we look at the so-called zero curvature equations for each pair of equations from \eqref{eq:tyxeqforZ}. 
These conditions are obtained inserting \eqref{eq:tyxeqforZ} in the identities $\partial_\geny \partial_\gent \itWsn= \partial_\gent \partial_\geny  \itWsn$, $\partial_\genx \partial_\gent \itWsn= \partial_\gent \partial_\genx  \itWsn$, and $\partial_\genx \partial_\geny \itWsn= \partial_\geny \partial_\genx  \itWsn$. 
These three equations become
\beqq \begin{split}
	& \itDsnl (z^2 \itAsn_1 + z \itAsn_2 + \itAsn_3) - \itDsnt (z\itAsn_1 + \itAsn_2) + z^2 [\itAsn_3, \itEsn] + z [\itAsn_3, \itAsn_1] + [\itAsn_3, \itAsn_2] =0 , 	\\
	& \itDsnh (z^2 \itAsn_1 + z \itAsn_2 + \itAsn_3) - \itDsnt \itAsn_1 + z^2 [\itAsn_2, \itEsn] + z [\itAsn_3, \itEsn] + z [\itAsn_2, \itAsn_1] + [\itAsn_3, \itAsn_1] =0 , 	\\
	& \itDsnh (z \itAsn_1 +  \itAsn_2)  - \itDsnl \itAsn_1 + z [\itAsn_2, \itEsn] + [\itAsn_2, \itAsn_1] =0 . 
\end{split} \eeqq
These equations hold for all $z$. 
Considering the coefficients of $z^i$, we obtain 8 identities. 
The $2$ identities from the last equation follow from the remaining equations. 
We also observe from \eqref{eq:eqfromlargez} that $\itDsnl \itAsn_2= \itDsnh \itAsn_3$. 
Using this, we see that the identity from the coefficients of $z$ of the first equation and the identity from the constant terms of the second equation are equivalent.
Thus, we are left with $5$ identities. Two of them are $\itDsnh \itAsn_1 + [\itAsn_2, \itEsn] =0$ and $\itDsnl \itAsn_1 + [\itAsn_3, \itEsn] =0$, 
but these follow from \eqref{eq:eqfromlargez}.  
Summarizing, the zero curvature conditions yield $3$ new identities, which are
\beq \label{eq:3eqnsfromzc} \begin{split}
	& \itDsnh \itAsn_2 - \itDsnl \itAsn_1 + [\itAsn_2, \itAsn_1] =0, \qquad 
	 \itDsnh \itAsn_3 - \itDsnt \itAsn_1 + [\itAsn_3, \itAsn_1] =0,\qquad 
	  \itDsnl  \itAsn_3 - \itDsnt \itAsn_2 + [\itAsn_3, \itAsn_2] =0.
\end{split} \eeq
Inserting \eqref{eq:eqfromlargez}, these equations become \eqref{eq:kpzmateqdf}.
\end{proof}

We now show how the last result implies  Lemma~\ref{lem:basic} and Theorem~\ref{thm:KdVheat}. 
Let $\itpermss$ be the $(m+1)\times (m+1)$ permutation matrix as introduced before Definition \ref{defn:Phi1d} and $\itI_j$ the identity matrix of size $j$.
Note that 
\begin{equation}
\label{eq:perm1onul}
	\itpermss \itEn (\itpermss)^T =
	\diag(\itI_{|\itss|},0,\cdots,0)
\end{equation}
with $\itEn$ as in \eqref{def:itEnitEss}.

\begin{cor}\label{nlsandkdv00}
Decompose $\itpermss \itZ_1 (\itpermss)^T$ as the sum of diagonal and off-diagonal block matrices,  
\begin{equation} \label{eq:decompitZ1block} 
	\itpermss \itZ_1 (\itpermss)^T  =  \kpYd+ \kpYo, 
	\qquad 
	\kpYd \deff \begin{pmatrix} \itqss & \itzero \\ \itzero & \itsss \end{pmatrix}, 
	\quad
	\kpYo \deff  \begin{pmatrix} \itzero & \itpss \\ \itrss & \itzero \end{pmatrix}, 
\end{equation}
where we used the notations of \eqref{eq:L1O1not}. Define the $(m+1)\times (m+1)$ diagonal matrix 
\begin{equation*}
	\kpV= \diag(\itI_{|\itss|},-\itI_{m+1-|\itss|}).
\end{equation*}
Then,  
\beq \label{eq:kpeqpfte} \begin{aligned}
	& \itDsnh \kpYd =  - \kpV \kpYo^2, \qquad
	&&\itDsnl \kpYd=    \kpYo ( \itDsnh \kpYo )  -(  \itDsnh \kpYo) \kpYo, \\
	& \itDsnl \kpYo  =  \kpV \big( \itDsnh^2 \kpYo + 2  \kpYo^3 \big) , \qquad 
	&&\itDsnt \kpYo =  \itDsnh^3 \kpYo + 3 ( \itDsnh \kpYo) \kpYo^2  + 3 \kpYo^2 ( \itDsnh \kpYo ). 
\end{aligned} \eeq
In terms of the matrices $\itqss, \itpss, \itrss, \itsss$ from \eqref{eq:L1O1not}, the first two equations become 
\beqq
	\itDsnh \itqss=-\itpss \itrss,  \qquad \itDsnh \itsss = \itrss \itpss , \qquad 
	\itDsnl \itqss= \itpss ( \itDsnh \itrss ) - ( \itDsnh \itpss ) \itrss,  
	\qquad \itDsnl \itsss =  \itrss ( \itDsnh \itpss ) - ( \itDsnh \itrss ) \itpss , 
\eeqq
the third equation becomes the coupled matrix NLS equation with complex time \eqref{eq:heat}, and the fourth equation becomes the coupled matrix mKdV equations \eqref{eq:mkdv}. Finally, the additional relation
$$
\Tr(\itsss)=-\Tr(\itqss)
$$
also holds.
\end{cor}

\begin{proof}
For each $n=1,2,3$, write in diagonal/off-diagonal block form,
\beqq
	 \itpermss \itZsn_n (\itpermss)^T=  \kpYd_n+ \kpYo_n,\qquad 
	 \kpYd_n = \begin{pmatrix} \ast & \itzero \\ \itzero & \ast \end{pmatrix},
	 \quad \kpYo_n =\begin{pmatrix} \itzero & \ast \\ \ast & \itzero \end{pmatrix}. 
\eeqq
In particular, $\kpYd_1= \kpYd$ and $\kpYo_1=\kpYo$. In the calculations that come, we make use of the identities
\begin{equation*}
	[   \itpermss \itZsn_n (\itpermss)^T , \itpermss \itEn (\itpermss)^T] = \kpYo_n \kpV = - \kpV \kpYo_n , \qquad 
	\kpYd_n  \kpV = \kpV \kpYd_n , \qquad 
	\kpV^2=\itI_{m+1}, 
\end{equation*}
which are straightforward to check. 

Conjugate the equations  \eqref{eq:kpzmateqdf} and \eqref{eq:kpzmateqdf123} by $\itpermss$ and consider the diagonal and off-diagonal blocks. The diagonal blocks of the first equation of  \eqref{eq:kpzmateqdf123}  give the identity  
\begin{equation}\label{eq:LxO200} 
	 \itDsnh \kpYd_1 =  - \kpV \kpYo_1^2 . 
\end{equation}
This is the first equation in \eqref{eq:kpeqpfte}. 

The off-diagonal blocks of the first equation of \eqref{eq:kpzmateqdf123}  and the diagonal blocks of the second equation yield 
\begin{equation*}
	\itDsnh \kpYo_1=  \kpV  \kpYo_{2} -  \kpV \kpYo_1   \kpYd_{1},  \qquad \itDsnl \kpYd_1=  -\kpV \kpYo_1  \kpYo_{2}  -(  \itDsnh \kpYd_1) \kpYd_1 -(  \itDsnh \kpYo_1) \kpYo_1 .
\end{equation*}
Using $\kpV^2=\itI_{m+1}$, the first equation implies $\kpYo_{2} = \kpV \itDsnh \kpYo_1 +  \kpYo_1 \kpYd_{1}$. 
Inserting this into the second equation and also using \eqref{eq:LxO200}, we find
$\itDsnl \kpYd_1=    \kpYo_1 ( \itDsnh \kpYo_1 )  -(  \itDsnh \kpYo_1) \kpYo_1$. 
This is the second equation in \eqref{eq:kpeqpfte}. 

The off-diagonal blocks of the first equation in \eqref{eq:kpzmateqdf} give us  
\begin{equation*}
 - \kpV \itDsnl \kpYo_1   +  \itDsnh^2 \kpYo_1 - [\itDsnh \kpYd_1, \kpV \kpYo_1  ] =0.
\end{equation*}
Inserting \eqref{eq:LxO200} and multiplying $\kpV$ on the left, we find 
$\itDsnl \kpYo_1  =  \kpV \itDsnh^2 \kpYo_1 + 2 \kpV \kpYo_1^3$,
which is the third equation in \eqref{eq:kpeqpfte}. 

The off-diagonal blocks of the second equation of \eqref{eq:kpzmateqdf} imply 
\begin{equation*}
	  - \kpV \itDsnt \kpYo_1 +  \itDsnh \itDsnl \kpYo_1 - [\itDsnl \kpYd_1, \kpV \kpYo_1 ]  =0  .
\end{equation*}
We remove $\itDsnl \kpYd_1$  and $\itDsnl \kpYo_1$ using the second and the third equations in \eqref{eq:kpeqpfte}, and find 
$\itDsnt \kpYo_1 =  \itDsnh^3 \kpYo_1 + 3 ( \itDsnh \kpYo_1) \kpYo_1^2  + 3 \kpYo_1^2 ( \itDsnh \kpYo_1 )$.
This is the fourth equation in \eqref{eq:kpeqpfte}. 

The equations for $\itqss, \itpss, \itrss, \itsss$ follow by inserting the formulas of $\kpYd$ and $\kpYo$. 

Finally, since $1=\det \itZsn(z) = \det (\itI_{m+1}+ \sum_{n=1}^\infty \frac{\itYb_n}{z^n})$ as $z\to \infty$, we have 
$\Tr \itZsn_1=0$, implying the identity $\Tr(\genq) + \Tr(\gens)=0$. 
\end{proof}

Thus, we proved Lemma~\ref{lem:basic} and Theorem~\ref{thm:KdVheat}.

\subsection{The matrix KP equation: proof of Theorem~\ref{thm:KPeqtion}}\label{sec:KPequationproof}

We derive the matrix KP equation from the equations \eqref{eq:kpeqpfte}. 

\begin{proof}[Proof of Theorem~\ref{thm:KPeqtion}]
We keep using the matrices $\kpYo,\kpYd,\kpV$ from Corollary~\ref{nlsandkdv00}, which satisfy $\kpYo \kpV = - \kpV \kpYo$, $\kpYd  \kpV = \kpV \kpYd$, and $\kpV^2=\itI_{m+1}$.  
Let $\ituss=\itpss\itrss$ and $\itvss=\itrss\itpss$.
Denote
\begin{equation*}
	\kpU\deff \kpYo^2= \begin{pmatrix} \itpss\itrss & \itzero \\ \itzero & \itrss\itpss \end{pmatrix} 
	= \begin{pmatrix} \ituss & \itzero \\ \itzero & \itvss \end{pmatrix}.
\end{equation*}
We first derive four equations for $\kpU$, $\kpYo$, and $\kpYd$. 
The first claimed equation follows from the definition of $\kpU$, 
\begin{equation}\label{eq:kpprfmtrx1}
	\itDsnh^3 \kpU =    ( \itDsnh^3 \kpYo) \kpYo+ \kpYo ( \itDsnh^3 \kpYo)  +  3 ( \itDsnh^2 \kpYo) (\itDsnh \kpYo)  + 	3 (\itDsnh \kpYo)  ( \itDsnh^2 \kpYo)  .
\end{equation}
To obtain the second equation, we use the second equation of \eqref{eq:kpeqpfte} to get 
\beqq \begin{split}
	\itDsnl^2 \kpYd
	&=    (\itDsnl\kpYo) ( \itDsnh \kpYo ) +  \kpYo ( \itDsnl\itDsnh \kpYo ) 
	- (  \itDsnl\itDsnh \kpYo) \kpYo  -(  \itDsnh \kpYo) (\itDsnl\kpYo ). 
\end{split} \eeqq
Using the third equation of \eqref{eq:kpeqpfte}, this implies our second equation, 
\beq \label{eq:kpprfmtrx2}
	-\kpV \itDsnl^2 \kpYd 
	=   ( \itDsnh^3 \kpYo) \kpYo+ \kpYo ( \itDsnh^3 \kpYo)  -   ( \itDsnh^2 \kpYo) (\itDsnh \kpYo) -  (\itDsnh \kpYo)  ( \itDsnh^2 \kpYo)  + 4 \kpYo ( \itDsnh \kpU)) \kpYo.
\eeq
Next, we multiply the fourth equation of \eqref{eq:kpeqpfte} by $\kpYo$ to the right, and separately to the left, and add both equations. The result is our third equation,
\beq \label{eq:kpprfmtrx3}
	\itDsnt \kpU =  ( \itDsnh^3 \kpYo) \kpYo+ \kpYo ( \itDsnh^3 \kpYo) + 3 \itDsnh (\kpU^2).
\eeq
Finally, for the fourth and last one, we compute $[\kpU,\itDsnl \kpYd]$ with the second equation of \eqref{eq:kpeqpfte}, obtaining
\beq \label{eq:kpprfmtrx4}
	- [ \kpU, \itDsnl \kpYd] + \itDsnh (\kpU^2) 
	= 2 \kpYo ( \itDsnh \kpU) \kpYo .
\eeq

We combine the above four equations to eliminate $\kpYo$ and obtain an equation of $\kpU$ and $\kpYd$ as follows: 
Multiplying the equation \eqref{eq:kpprfmtrx2} by $3$ and adding it to the equation \eqref{eq:kpprfmtrx1}, and then subtracting the  equation \eqref{eq:kpprfmtrx3} multiplied by $4$ and also the equation \eqref{eq:kpprfmtrx4} multiplied by $6$, we obtain 
\beqq
	\itDsnh^3 \kpU - 3 \kpV \itDsnl^2 \kpYd - 4 \itDsnt \kpU   + 6 [ \kpU,  \itDsnl \kpYd] 
	+ 6 \itDsnh (\kpU^2) =0  .
\eeqq
If we insert $\kpU = \big( \begin{smallmatrix} \ituss & \itzero \\ \itzero & \itvss \end{smallmatrix}\big)$ and 
$\kpYd=  \big( \begin{smallmatrix} \itqss & \itzero \\ \itzero & \itsss \end{smallmatrix}\big)$, we obtain the first equations in \eqref{eq:kpu} and \eqref{eq:kpv}. The second equations of them follow from $\itDsnh \kpYd =  - \kpV \kpU$, which is the first equation of \eqref{eq:kpeqpfte}. Thus the result is proved. 
\end{proof}

\subsection{Multi-component KP hierarchy}\label{sec:KPhierarchy}

In this subsection, we show that the three Lax equations \eqref{eq:laxeq} with respect to the parameters are related to multi-component KP hierarchy theory. 
This will give us another proof of Theorem~\ref{thm:KPeqtion}. The article \cite{kac_van_de_leur} gives a thorough exposition on the multi-component KP hierarchy theory, see also \cite{Teo11}. For the scalar KP hierarchy theory, a helpful overview is found in \cite{BaloghHarnad2021}. 
Multi-component KP hierarchy theory are already known to be connected with models related to the Airy$_2$ process and random matrix theory, see  \cites{AdlerVanMoerbekeVanhaecke2009,AdlerCafassoVanMoerbeke2012} and references therein. We can also notice that cubic admissible determinant is multi-component KP tau function, but we don't use this fact in this paper.

The cubic admissible operator $\genH$ contains $3m$ parameters $\ith_1, \cdots, \ith_m, \itl_1, \cdots, \itl_m, \itt_1, \cdots, \itt_m$.  
We relabel them as
\beqq
\begin{split} 
    x^{(i)}_1=\ith_i,\qquad x^{(i)}_2=\itl_i,\qquad x^{(i)}_3=\itt_i
,\qquad i\in \{1,\ldots , m\}. 
\end{split}
\eeqq 
Recall the matrices $\itZ_k$ for $k\ge 1$ from \eqref{eq:pjinfg}. 
Fix $n\in\{1, \cdots, m\}$ and write in block form\footnote{We can also consider $\itpermss \itZ_k (\itpermss)^T$ similar to \eqref{eq:L1O1not}. Then the discussions below still hold if we define the parameters $x^{(i)}_\ell$, $i=1, \cdots, |\itss|$, $\ell=1, 2,3$ appropriately. For the convenience of presentation, we state results only when $\itpermss=\{1, \cdots, n\}$ here.}
\beq \label{eq:phikblok}
    {\itZ}_k=\begin{pmatrix} \itqss_k & \itpss_k \\ \itrss_k & \itsss_k \end{pmatrix},
\eeq
where $\itqss_k$ is a matrix of size $n\times n$. 
For $k\ge 1$, let $\itE_{j}^{(k)}$ be the $k\times k$ matrix with $1$ in the $j$-th diagonal entry and zero in the other entries. 

Now define the formal $n\times n$ matrix-valued pseudo-differential operators on the variables $x_1^{(1)},\cdots, x_1^{(n)}$, 
\beq\label{eq:Pdef11}
    \partial\deff \sum_{j=1}^{n}{\partial_{x^{(j)}_1}} \qquad 
    \text{and} \qquad 
    \itPo\deff \itI_n+\sum_{k=1}^{\infty}\itqss_k\partial^{-k},
\eeq
and for $j \in \{ 1,\cdots, n\}$ and $\ell \in \{1,2, 3\}$ also introduce
\beq \label{eq:Dkjdenf}
    \laxKP_{\ell}^{(j)}\deff(\itPo\itE_{j}^{(n)}\partial^\ell\itPo^{-1})_+.
\eeq 
We stress that this convention for $\partial$ is only valid for this subsection.
The subindex $+$ in \eqref{eq:Dkjdenf} means that we take the part of the term between brackets that contains only non-negative powers of $\partial$, and the inverse $\itPo^{-1}$ is in the sense of pseudo-differential operators, so that 
$$
\itPo^{-1}=\itI_{n}-\itqss_1\partial^{-1}+(\itqss_1^2-\itqss_2)\partial^{-2}
+(-\itqss_1^3-\itqss_1(\partial\itqss_1) + \itqss_1\itqss_2+ \itqss_2\itqss_1 - \itqss_3 )\partial^{-3}+\cdots. 
$$
From the definition, we have\footnote{To compute it, note, for example, that $\partial^3 u \partial^{-2}= u \partial+ 3 (\partial u)  + 3(\partial^2 u) \partial^{-1} + (\partial^3 u)\partial^{-2}$  and $\partial^{-1} u \partial^2=u\partial - (\partial u)+ (\partial^2 u) \partial^{-1}-\cdots$ from the product rule of derivatives, and thus,  $(\partial^3 u \partial^{-1})_+= u \partial^2+ 3 (\partial u) \partial$  and $(\partial^{-1} u \partial^2)+=u\partial - (\partial u)$ for any function $u$.}
\beq \label{eq:Dj}
\begin{split} 
    &\laxKP_{1}^{(j)}=\itE_{j}^{(n)}\partial+[\itqss_1,\itE_{j}^{(n)}],\\ 
    &\laxKP_{2}^{(j)}=\itE_{j}^{(n)}\partial^2+[\itqss_1,\itE_{j}^{(n)}]\partial+[\itqss_2,\itE_{j}^{(n)}]-[\itqss_1,\itE_{j}^{(n)}]\itqss_1-2\itE_{j}^{(n)}(\partial \itqss_1)\\
    &\laxKP_{3}^{(j)}= \itE_{j}^{(n)}\partial^3 +[\itqss_1,\itE_{j}^{(n)}]\partial^2 + \big( [\itqss_2,\itE_{j}^{(n)}]-[\itqss_1,\itE_{j}^{(n)}]\itqss_1
    -3\itE_{j}^{(n)}(\partial \itqss_1) \big) \partial \\
    &\qquad\qquad +[\itqss_3,\itE_{j}^{(n)}] +[\itqss_1,\itE_{j}^{(n)}]\itqss_1^2 
    -[\itqss_1,\itE_{j}^{(n)}]\itqss_2-[\itqss_2,\itE_{j}^{(n)}]\itqss_1 \\
    &\qquad\qquad  -3\itE_{j}^{(n)}(\partial \itqss_2) -3\itE_{j}^{(n)}(\partial^2 \itqss_1)-2[\itqss_1,\itE_{j}^{(n)}](\partial \itqss_1)+3\itE_{j}^{(n)}(\partial \itqss_1)\itqss_1 .
\end{split}
\eeq 
These are particular combinations of $\itqss_1, \itqss_2, \itqss_3$ and their derivatives. 

There are several equivalent formulations for multi-component KP hierarchy. Remark 4.4 of \cite{kac_van_de_leur} gives a good summary of their relations. 
Here, we prove the following version.

\begin{prop}[Sato equations for multi-component KP hierarchy] \label{prop:Sato}
Fix $1\le n\le m$. 
The $n\times n$ matrix-valued pseudo-differential operators $\itPo$ and  $\laxKP_{\ell}^{(j)}$ from \eqref{eq:Pdef11} and \eqref{eq:Dj} satisfy the following equations, known as the Sato equations\footnote{This equation is the special case when $\alpha=0$ in equation (113) of \cite{kac_van_de_leur}. The full Sato equations involve pseudo-differential operators parametrized by $\alpha\in \Z^{n-1}$. The multi-component KP hierarchy may also involve infinitely many parameters $x^{(i)}_\ell$, $i\in \{1, \cdots, n\}$, $\ell\in \{1, 2, \cdots\}$. Here we have the special case when $\ell\in \{1,2, 3\}$.} for the $n$-component KP hierarchy: 
\beq \label{Sato}
    \dfrac{\partial\itPo}{\partial{x^{(i)}_\ell}}=\laxKP_\ell^{(i)}\itPo-\itPo\itE_{i}^{(n)}\partial^\ell
    \qquad
    \text{for } i\in \{ 1, \cdots, n \} \text{ and } \ell\in \{1,2,3\}.
\eeq
\end{prop}

\begin{proof}
Recall the Lax equations \eqref{eq:laxeq}. 
Inserting $\itWsn = \itZsn \itconjb$ into \eqref{eq:laxeq} and collecting coefficients of $\Boh(z^{-k})$ we obtain the identities
\beq \label{eq:zeq}\begin{split}
&\partial_{x^{(i)}_1} {\itZ}_k = - [{\itZ}_{k+1}, \itE_{i}^{(n+1)}]+{\laxC}_{1}^{(i)}{\itZ}_{k} ,\\
&	\partial_{x^{(i)}_2} {\itZ}_k = - [{\itZ}_{k+2}, \itE_{i}^{(n+1)}]+{\laxC}_{1}^{(i)}{\itZ}_{k+1}+{\laxC}_{2}^{(i)}{\itZ}_{k},\\
&\partial_{x^{(i)}_3}{\itZ}_k = - [{\itZ}_{k+3},\itE_{i}^{(n+1)}]+{\laxC}_{1}^{(i)}{\itZ}_{k+2}+{\laxC}_{2}^{(i)}{\itZ}_{k+1}+{\laxC}_{3}^{(i)}{\itZ}_{k},
\end{split} \eeq 
for $k\in \{1, 2, \cdots\}$ and $i\in\{1, \cdots, n\}$, where $\laxC_{1}^{(i)}$ are given by \eqref{eq:Aiinmulhi}. 
Using the formula of $\laxC_{1}^{(i)}$ and inserting the block form \eqref{eq:phikblok}, the first equation of \eqref{eq:zeq} yields 
\beq \label{eq:zeqblock}\begin{split}
    &\partial_{x^{(i)}_1} \itqss_k = - [\itqss_{k+1}, \itE_{i}^{(n)}]+[\itqss_{1}, \itE_{i}^{(n)}]\itqss_k-\itE_{i}^{(n)}\itpss_1\itrss_k 
\end{split} \eeq 
and
\beq \label{eq:zeqblock2}\begin{split}
    &\partial_{x^{(i)}_1} \itpss_k = \itE_{i}^{(n)}\itpss_{k+1}+[\itqss_{1}, \itE_{i}^{(n)}]\itpss_k-\itE_{i}^{(n)}\itpss_1\itsss_k ,\\
    &\partial_{x^{(i)}_1} \itrss_k = -\itrss_{k+1}\itE_{i}^{(n)}+\itrss_1\itE_{i}^{(n)}\itqss_k ,\\
    &\partial_{x^{(i)}_1} \itsss_k = \itrss_1\itE_{i}^{(n)}\itpss_k ,\\
\end{split} \eeq 
Summing over $i=1, \cdots, n$ and noting $\itE_{1}^{(n)}+\cdots+ \itE_{n}^{(n)}=\itI_n$, we find 
\begin{equation} \label{eq:zeqsum}\begin{split}
    &\partial \itqss_k = -\itpss_1\itrss_k 
\end{split} \end{equation} 
and
\begin{equation} \label{eq:zeqsumprs}\begin{split}
    \partial \itpss_k =  \itpss_{k+1}-\itpss_1\itsss_k , \qquad
    \partial \itrss_k = -\itrss_{k+1}+\itrss_1\itqss_k ,\qquad
    \partial \itsss_k = \itrss_1\itpss_k .
\end{split} \end{equation}
Using these four equations we can also compute higher derivatives of $\itqss_k$: 
\beq \label{eq:zeqhighersum} \begin{split}
    &\partial^2 \itqss_k=-\itpss_2\itrss_k+\itpss_1\itsss_1\itrss_k+\itpss_1\itrss_{k+1}-\itpss_1\itrss_1\itqss_k
,\\
    &\partial^3 \itqss_k=2 \itpss_1\itsss_1\itrss_1\itqss_k+\itpss_1\itrss_1\itqss_{k+1}+\itpss_1\itrss_{2}\itqss_k-2 \itpss_{2}\itrss_1\itqss_k-\itpss_1\itrss_1\itqss_1\itqss_k -2 \itpss_1\itsss_1\itrss_{k+1} \\
    &\qquad \qquad +\itpss_1\itsss_{2}\itrss_k+\itpss_{2}\itsss_1\itrss_k -\itpss_1\itsss_1^2\itrss_k -\itpss_1\itrss_{k+2}+2 \itpss_{2}\itrss_{k+1}-\itpss_{3}\itrss_k+2 \itpss_1\itrss_1\itpss_1\itrss_k .
\end{split} \eeq 

In addition, the top left blocks of the second and third equation of \eqref{eq:zeq} are
\beq\label{eq:zeqotherblock} \begin{split} 
    \partial_{x^{(i)}_2} \itqss_k 
    & = - [\itqss_{k+2}, \itE_{i}^{(n)}]+[\itqss_{1}, \itE_{i}^{(n)}]\itqss_{k+1}-\itE_{i}^{(n)}\itpss_1\itrss_{k+1}+[\itqss_{2}, \itE_{i}^{(n)}]\itqss_{k} \\
    & \qquad -\itE_{i}^{(n)}\itpss_2\itrss_{k}-[\itqss_{1}, \itE_{i}^{(n)}]\itqss_{1}\itqss_{k}+\itE_{i}^{(n)}\itpss_1\itrss_{1}\itqss_{k} -[\itqss_{1}, \itE_{i}^{(n)}]\itpss_{1}\itrss_{k}+\itE_{i}^{(n)}\itpss_{1}\itsss_1\itrss_{k}
\end{split} \eeq 
and
\beq\label{eq:zeqotherblock2} \begin{split} 
     \partial_{x^{(i)}_3}\itqss_k & = - [\itqss_{k+3},\itE_{i}^{(n)}]+[\itqss_{1}, \itE_{i}^{(n)}]\itqss_{k+2}-\itE_{i}^{(n)}\itpss_1\itrss_{k+2}+[\itqss_{2}, \itE_{i}^{(n)}]\itqss_{k+1} \\
     &\qquad -\itE_{i}^{(n)}\itpss_2\itrss_{k+1}-[\itqss_{1}, \itE_{i}^{(n)}]\itqss_{1}\itqss_{k+1} +\itE_{i}^{(n)}\itpss_1\itrss_{1}\itqss_{k+1}-[\itqss_{1}, \itE_{i}^{(n)}]\itpss_{1}\itrss_{k+1}  \\
     &\qquad +\itE_{i}^{(n)}\itpss_{1}\itsss_1\itrss_{k+1}+[\itqss_{3}, \itE_{i}^{(n)}]\itqss_{k}-\itE_{i}^{(n)}\itpss_3\itrss_{k}-[\itqss_{1}, \itE_{i}^{(n)}]\itqss_{2}\itqss_{k}+\itE_{i}^{(n)}\itpss_1\itrss_{2}\itqss_{k}\\
     &\qquad -[\itqss_{1}, \itE_{i}^{(n)}]\itpss_{2}\itrss_{k}+\itE_{i}^{(n)}\itpss_{1}\itsss_2\itrss_{k}-\itE_{i}^{(n)}\itpss_{1}\itsss_{1}\itrss_{1}\itqss_{k}+[\itqss_{1},\itE_{i}^{(n)}]\itpss_{1}\itsss_{1}\itrss_{k}+\itE_{i}^{(n)}\itpss_{2}\itrss_{1}\itqss_{k}  \\
     &\qquad -[\itqss_{2},\itE_{i}^{(n)}]\itpss_{1}\itrss_{k}-\itE_{i}^{(n)}\itpss_{1}\itrss_{1}\itqss_{1}\itqss_{k} +[\itqss_{1},\itE_{i}^{(n)}]\itpss_{1}\itrss_{1}\itqss_{k}+[\itqss_{1},\itE_{i}^{(n)}]\itqss_{1}\itpss_{1}\itrss_{k}+\itE_{i}^{(n)}\itpss_{2}\itsss_{1}\itrss_{k}  \\
     &\qquad -\itE_{i}^{(n)}\itpss_{1}\itsss_{1}^2\itrss_{k}-\itE_{i}^{(n)}\itpss_{1}\itrss_{1}\itpss_{1}\itrss_{k}+[\itE_{i}^{(n)},\itqss_{2}]\itqss_{1}\itqss_{k}+[\itqss_{1},\itE_{i}^{(n)}]\itqss_{1}^2\itqss_{k} .
\end{split} \eeq 
Using \eqref{eq:zeqsum}, and \eqref{eq:zeqhighersum}, the right-hand side of \eqref{eq:zeqblock}, \eqref{eq:zeqotherblock}, and \eqref{eq:zeqotherblock2} can be expressed using only $\itqss_k$'s. We obtain  
\beq \label{eq:zeqblockq}\begin{split}
    &\partial_{x^{(i)}_1} \itqss_k = - [\itqss_{k+1}, \itE_{i}^{(n)}]+[\itqss_{1}, \itE_{i}^{(n)}]\itqss_k+\itE_{i}^{(n)}\partial\itqss_k , 
\end{split} \eeq 
\beq\label{eq:zeqotherblockq} \begin{split} 
    \partial_{x^{(i)}_2} \itqss_k 
    & = - [\itqss_{k+2}, \itE_{i}^{(n)}]+[\itqss_{1}, \itE_{i}^{(n)}]\itqss_{k+1}-[\itqss_{1}, \itE_{i}^{(n)}]\itqss_{1}\itqss_{k}+[\itqss_{2}, \itE_{i}^{(n)}]\itqss_{k} \\
    & \qquad +[\itqss_{1}, \itE_{i}^{(n)}]\partial\itqss_{k}+\itE_{i}^{(n)}\partial^2\itqss_{k}+2\itE_{i}^{(n)}\partial\itqss_{k+1}-2\itE_{i}^{(n)}(\partial\itqss_{1})\itqss_{k} , 
\end{split} \eeq 
and
\beq\label{eq:zeqotherblock2q} \begin{split} 
     \partial_{x^{(i)}_3}\itqss_k & = - [\itqss_{k+3},\itE_{i}^{(n)}]+[\itqss_{1}, \itE_{i}^{(n)}]\itqss_{k+2}+[\itqss_{2}, \itE_{i}^{(n)}]\itqss_{k+1}-[\itqss_{1}, \itE_{i}^{(n)}]\itqss_{1}\itqss_{k+1} \\
      &\qquad -[\itqss_{1}, \itE_{i}^{(n)}]\itqss_{2}\itqss_{k}+[\itqss_{3}, \itE_{i}^{(n)}]\itqss_{k}-[\itqss_{2},\itE_{i}^{(n)}]\itqss_{1}\itqss_{k}+[\itqss_{1},\itE_{i}^{(n)}]\itqss_{1}^2\itqss_{k}\\
     &\qquad +3\itE_{i}^{(n)}\partial\itqss_{k+2}+3\itE_{i}^{(n)}\partial^2\itqss_{k+1} +\itE_{i}^{(n)}\partial^3\itqss_{k}+2[\itqss_{1}, \itE_{i}^{(n)}]\partial\itqss_{k+1}  \\
     &\qquad +[\itqss_{1}, \itE_{i}^{(n)}]\partial^2\itqss_{k}+[\itqss_{2}, \itE_{i}^{(n)}]\partial\itqss_{k}-3\itE_{i}^{(n)}(\partial\itqss_{1})\itqss_{k+1}-3\itE_{i}^{(n)}(\partial\itqss_{1})(\partial\itqss_{k})-3\itE_{i}^{(n)}(\partial\itqss_{2})\itqss_{k}  \\
     &\qquad -3\itE_{i}^{(n)}(\partial^2\itqss_{1})\itqss_{k}-[\itqss_{1}, \itE_{i}^{(n)}]\itqss_{1}\partial\itqss_{k} -2[\itqss_{1}, \itE_{i}^{(n)}](\partial\itqss_{1})\itqss_{k}+3\itE_{i}^{(n)}(\partial\itqss_{1})\itqss_{1}\itqss_{k}.
\end{split} \eeq 

Now consider again \eqref{Sato}. 
In its the left-hand side, 
we use the definition $\itPo\deff \itI_n+\sum_{k=1}^{\infty}\itqss_k\partial^{-k}$ and take the derivative, $\dfrac{\partial\itPo}{\partial{x^{(i)}_\ell}}$. We then insert the formulas 
\eqref{eq:zeqblockq}, \eqref{eq:zeqotherblockq}, and \eqref{eq:zeqotherblock2q}.
For the right-hand side, we insert the definition of $\itPo$ and the formula \eqref{eq:Dj}. A tedious computation shows that the two sides are equal, and concluding the result.
\end{proof} 


Sato equations imply the following equations.

\begin{cor}[Zakharov-Shabat equations for multi-component KP hierarchy] \label{prop:KPhi}
Fix $1\le n\le m$. 
The $n\times n$ matrix-valued pseudo-differential operators $\laxKP_{\ell}^{(j)}$ from \eqref{eq:Dj} satisfy the following equations, known as the Zakharov-Shabat equations for the $n$-component KP hierarchy: 
\begin{equation}\label{eq:mkphier}
    	\dfrac{\partial \laxKP_\ell^{(i)}}{\partial x_{\ell'}^{(j)}}-\dfrac{\partial \laxKP_{\ell'}^{(j)}}{\partial x_\ell^{(i)}}=[\laxKP_{\ell'}^{(j)}, \laxKP_\ell^{(i)} ],
	\qquad i,  j \in \{ 1,\cdots, n\}, \quad \ell,  \ell' \in \{1,2, 3\}.     
\end{equation}
\end{cor}

\begin{proof}
Proposition 4.4 of \cite{kac_van_de_leur} proved that if $\itPo$ and $\laxKP_\ell^{(i)}$ satisfy the Sato equations \eqref{Sato}, then $\laxKP_\ell^{(i)}$ satisfy \eqref{eq:mkphier}. 

\end{proof}

Sato equations also give 
another proof of the matrix KP equation \eqref{eq:kpu}. 
The proof of the next result shows that 
Zakharov-Shabat equations \eqref{eq:mkphier} yield the matrix KP equation modulo an ``integration constant" $\theta$: see \eqref{eq:KPwiththeta} below. 
We can then show that $\theta=0$ using one of Sato equations.


\begin{cor}\label{cor:kp}
Set
\beqq
	\genu\deff -\partial_\ith \itqss \qquad \text{where} \quad \itqss\deff \itqss_1. 
\eeqq
Then, the matrix KP equation \eqref{eq:kpu} holds:
\beq \label{eq:KPnewder}
 	-4 \partial_{\itt} \genu  + \partial_{\ith}^3 \genu +  6 \partial_{\ith} (\genu^2) - 3 \partial_{\itl}^2 \itqss + 6[ \genu, \partial_{\itl} \itqss  ]=0. 
 \eeq 
\end{cor}

\begin{proof}
We use the notations $x^{(i)}_1=\ith_i$, $x^{(i)}_2=\itl_i$, $x^{(i)}_3=\itt_i$, 
and
\beqq
	\partial_\itt= \sum_{i=1}^n \partial_{\itt_i} , \qquad \partial_\itl= \sum_{i=1}^n \partial_{\itl_i} , \qquad
	\partial_\ith= \sum_{i=1}^n \partial_{\ith_i} = \partial .
\eeqq
Set $\laxKP_{\ell}\deff \sum_{j=1}^n \laxKP_{\ell}^{(j)}$ for $\ell\in \{1,2,3\}$.
From \eqref{eq:Dj}, we find that 
\beq \label{eq:Djsum}
\begin{split} 
    &\laxKP_{1}= \partial, \qquad \laxKP_{2} = \partial^2-2 (\partial \itqss_1),\qquad  \laxKP_{3} = \partial^3   -3 (\partial \itqss_1)\partial -3 (\partial \itqss_2) -3 (\partial^2 \itqss_1) +3 (\partial \itqss_1)\itqss_1 .
\end{split}
\eeq 
Setting 
\beqq 
	\alpha \deff   - (\partial \itqss_2) + (\partial \itqss_1)\itqss_1
\eeqq
and using $\genu\deff -\partial_\ith \itqss_1$, these expressions become
\beq \label{eq:Djsum2}
\begin{split} 
    &\laxKP_{1}= \partial, \qquad \laxKP_{2} = \partial^2 +2 \genu  \qquad  
    \laxKP_{3} = \partial^3   +3 \genu \partial + 3(\partial_{\ith} \genu +\alpha) . 
\end{split}
\eeq

Consider the $\ell'=3$ and $\ell=2$ case of the equations \eqref{eq:mkphier}.  Taking the sum over $i, j=1, \cdots, n$, it becomes 
\begin{equation}\label{eq:mkphiersum}
    	\partial_\itt \laxKP_2 - \partial_{\itl} \laxKP_3 =[\laxKP_3, \laxKP_2 ]   .
\end{equation}
Inserting \eqref{eq:Djsum2} and considering the terms of $\partial^1$ and $\partial^0$ we obtain two equations, namely
\beq \label{eq:hi1}
	\partial_{\itl} \genu = \partial_{\ith}^2 \genu + 2 \partial_{\ith}\alpha
\eeq
and
\beq \label{eq:hi2}
	2 \partial_{\itt} \genu - 3 \partial_{\itl} \partial_{\ith} \genu 
	+\partial_{\ith}^3 \genu - 6 \genu \partial_{\ith} \genu 
	= 3\partial_{\itl} \alpha - 3 \partial_{\ith}^2 \alpha + 6 [ \partial_{\ith} \genu+ \alpha, \genu] .
\eeq
Since $\genu\deff -\partial_\ith \itqss$, \eqref{eq:hi1} implies that there is a function $\theta$ such that 
\beq \label{eq:hi3}
	\alpha= - \frac12 \left( \partial_{\ith} \genu + \partial_{\itl} \itqss_1+ \theta \right) 	
	\qquad\text{and}\qquad\partial_{\ith} \theta =0.
\eeq
Inserting it into \eqref{eq:hi2}, we obtain 
 \beq \label{eq:KPwiththeta}
 	-4 \partial_{\itt} \genu  + \partial_{\ith}^3 \genu +  6 \partial_{\ith} (\genu^2) - 3( \partial_{\itl}^2 \itqss  + \partial_{\itl} \theta )+ 6[ \genu, \partial_{\itl} \itqss +\theta ]=0. 
 \eeq

Now, the case of $\ell=2$ of Sato equations  \eqref{Sato} is precisely \eqref{eq:zeqotherblockq}. 
The sum of this equation over $i$ implies that $- (\partial \itqss_2) + (\partial \itqss_1)\itqss_1= - \frac12 \left( \partial_{\ith} \genu + \partial_{\itl} \itqss_1 \right)$. 
This shows that $\theta=0$ in \eqref{eq:KPwiththeta}, completing the proof.
\end{proof}

\subsection{Matrix ODE system \`a la Tracy and Widom: proof of Theorem~\ref{thm:matrixODE} and Corollary~\ref{thm:matrixp2}}\label{sec:proofPII}

In this section, we assume that $\itH$ is a {\it strongly cubic integrable operator} as in Definition \ref{def:cubicintegrabl}. 
Thus, $\Omega$ is a finite union of simple contours without endpoints and, furthermore, the matrix $\itJ_0$ in \eqref{eq:rhpjupp} is constant on each connected component of $\Omega$. 
The last property allows us to derive an additional Lax equation for the derivative with respect to the variable $z$, which is often called the spectral variable. Because this equation will have coefficients which are polynomial in $z$, the corresponding Lax pairs are usually associated to isomonodromic deformations \cite{Fokas-Its-Kapaev06}, but we will not explore this fact explicitly.

\begin{proof}[Proof of Theorem~\ref{thm:matrixODE}]
Recall that 
\beqq   
    \partial= \sum_{j=1}^m \itt_j \partial_{\ith_j}
\eeqq
was defined in \eqref{def:diffoperPform}. 
Taking a weighted sum of the first equation of \eqref{eq:laxeq}, we find that 
\beqq
	 \partial\itW(z) = ( z \whkptm+ \laxD_1 ) \itW(z) \qquad \text{with} \quad \laxD_1\deff [\itZ_1, \whkptm].
\eeqq
where recall $\whkptm \deff \diag( \itt_1, \cdots, \itt_m, 0)$ from \eqref{eq:multibytyx}. 
Inserting $\itWsn(z) = \itZsn(z) \itconjb (z)$, using the asymptotic series \eqref{eq:expYk}, and comparing the terms of $\Boh(z^{-1})$, we obtain
\beq \label{eq:X1X2tauAX1}
	\partial\itZ_1 = - [\itZ_{2}, \whkptm] + [\itZ_1, \whkptm]\itZ_1 . 
\eeq
The lower-right block of this identity is exactly~\eqref{eq:parsinrmp}. 

On the other hand, consider the differentiation $\partial_z$ with respect to the variable $z$. We have 
\beq \label{eq:Deltapatialz}
	\partial_z\itconjb=\sum_{j=1}^m (3z^2\itt_j  + 2z\itl_j+\ith_j  )\itE_j \itconjb 
	= (3z^2\whkptm + 2z\whkpym+\whkpxm ) \itconjb.
\eeq 
Since $\itJ_0$ is independent of $z$ on each connected component of the contour, $\partial_z \itW$ satisfies the same jump condition as $\itW$, and hence, $(\partial_z \itW)\itW^{-1}$ is an entire function. 
Considering the large $z$ asymptotic formula again, we obtain 
\beq \label{Z1prime}
	\partial_z \itW(z) = (3z^2 \whkptm + z (3\laxD_1+2\whkpym) + (3\laxD_2+2\laxD_3+\whkpxm))\itW(z)
\eeq
with $\laxD_1\deff [\itZ_1, \whkptm]$ as before, and 
\begin{equation*}
	\laxD_2\deff  [\itZ_2, \whkptm] - [\itZ_1, \whkptm] \itZ_1 = - \partial\itZ_1,\quad \laxD_3\deff [\itZ_1, \whkpym],
\end{equation*}
where we used \eqref{eq:X1X2tauAX1} for the last equality. 

Thus, $\itW$ is a common solution to the $\partial$ and $\partial_z$ differential equations, and the compatibility of these two equations, namely $\partial \partial_z\itW = \partial_z\partial \itW$, implies the zero curvature equation   
\beqq
	z (3[\whkptm,\laxD_2]+2[\whkptm, \itC_3]+2[\laxD_1, \whkpym] - 3\partial\laxD_1)
	+( 3[\laxD_1, \laxD_2]+2[\laxD_1,\laxD_3]+[\laxD_1,\whkpxm]  - 3\partial\laxD_2-2\partial\laxD_3 ) =0 . 
\eeqq
Using the formulas for $\laxD_1,\laxD_2$ and $\laxD_3$ and the cyclicity of the commutator, we obtain that the coefficient of $z$ is trivially identically zero. 
The constant term, on the other hand, is non-trivial and gives the equation 
\beqq
	3[\laxD_1, \laxD_2]+2[\laxD_1,\laxD_3]+[\laxD_1,\whkpxm]  - 3\partial\laxD_2-2\partial\laxD_3 =0 .
\eeqq
Inserting the formula of $\laxD_1$, $\laxD_2$, $\laxD_3$, we obtain \eqref{eq:matrxp2ii}. 
\end{proof}

\begin{proof}[Proof of Corollary~\ref{thm:matrixp2}]
The equation \eqref{eq:thPformSystem} follows by inserting $\itZ_1 = \big( \begin{smallmatrix} \itqsn & \itpsn \\ \itrsn & \itssn \end{smallmatrix} \big)$ into \eqref{eq:matrxp2ii} after removing $\partial\itssn$ using 
\eqref{eq:parsinrmp}.
If $\itt_1=\cdots= \itt_m\revdeff\itt$, the $m\times m$ upper left corner of \eqref{eq:X1X2tauAX1} yields $\itDsnh \itqsn=-\itpsn\itrsn$. 
This was also already proved in Lemma \ref{lem:basic}. 
Using this equation, the bottom two equations of \eqref{eq:thPformSystem} become the bottom two equations of \eqref{eq:thPformSystemred}. 
\end{proof}

\subsection{Deformation formulas for cubic admissible functions: proof of Proposition~\ref{lem:basic2}}\label{RHP}

Assume that $\itH$ is cubic integrable and $\itid-\itH$ is invertible, and consider the Fredholm determinant $\genFunc=\det(\itid-  \itH )$. 
We compute the derivatives of $\log \genFunc$ with respect to the variables $\ith_i, \itl_i, \itt_i$. 

\begin{proof}[Proof of Proposition~\ref{lem:basic2}]
Recall that $\itE_{i}$ is the $(m+1)\times(m+1)$ matrix whose $(i,i)$ entry is $1$ and all other entries are $0$.
From \eqref{def:mjexpon} and \eqref{eq:itconjb}, we have 
$\itconjb(z)= \diag\left(\genm_1(z),\genm_2(z),\hdots,\genm_m(z),1\right)$ with 
$\genm_j(z)= e^{\gent_j z^3 + \geny_j z^2 + \genx_j z}$. Thus, 
\beqq
	\partial_{\genx_i} \itconjb(z) = z \itE_i \itconjb(z). 
\eeqq
Hence, the vector functions $\genf$ and $\geng$ satisfy  
\beqq
	\partial_{\ith_i} \genf(z)  = z  \itE_{i} \genf(z)+ \alpha(z) \genf(z), \qquad 
	\partial_{\ith_i} \geng(z) =- z  \itE_{i} \geng(z)- \alpha(z) \geng(z)
	\qquad \text{where $\alpha(z)\deff \frac{\partial_{\genx_i} \kpc(z)}{\kpc(z)}$.} 
\eeqq
Thus, for $z\neq w$, 
\beq \label{eq:Hxde}
	\partial_{\genx_i} \genH(z,w)= \frac{ (\partial_{\genx_i} \genf(z))^T \geng(w)+ \genf(z)^T \partial_{\genx_i} \geng(w)}{z-w}
	= \genf(z)^T \itE_{i} \geng(w) + \alpha(z) \genH(z,w) - \genH(z,w) \alpha(w). 
\eeq 
For $z=w$, since $\genH(z,z)=0$ by definition, we have $\partial_{\genx_i} \genH(z,z)=0$. 
On the other hand, noting that $\Delta$ is a diagonal matrix, the condition \eqref{eq:withEzero} implies that $\genf(z)^T \itE_{i} \geng(z)=0$. 
Thus, the right-hand side of \eqref{eq:Hxde} is zero when $z=w$. 
Therefore, the equation \eqref{eq:Hxde} is valid for all $z, w\in \Omega$. 

For vector functions $u$ and $v$, let $u\otimes v$ be the integral operator given by the kernel $u(z)^T v(w)$.  
Let $\alpha$ be the operator of multiplication by $\alpha(z)$. 
Then, the equation \eqref{eq:Hxde} can be written in the operator form $\partial_{\genx_i} \genH= \genf \otimes (\itE_{i} \geng)+[\itH,\alpha]$ with $[\cdot,\cdot]$ being the commutator of operators.
Thus, 
\begin{equation}\label{eq:partHtr1}
	\partial_{\genx_i} \log \det(\itid- \genH)  = - \Tr ((\itid-\genH)^{-1} \partial_{\genx_i} \genH) 
	= -\Tr ( (\itid-\genH)^{-1}\genf \otimes (\itE_{i} \geng))+ \Tr ( (\itid-\genH)^{-1} [\itH,\alpha]).
\end{equation}
Since $(\itid-\itH)^{-1} \itH =  \itH (\itid-\itH)^{-1}$, we use cyclicity of the trace and obtain 
\begin{equation}\label{eq:partHtr2}
	\Tr ((\itid-\genH)^{-1} [\itH,\alpha]) =  \Tr ((\itid-\itH)^{-1} \genH\alpha)) -\Tr (\alpha \genH (\itid-\genH)^{-1} )= 0.
\end{equation}
Therefore,  using \eqref{eq:ResolventKernel}, 
\begin{equation}\label{eq:partHtr3}
\begin{split}
	\partial_{\genx_i} \log \det(\itid- \genH) 
	=  -\Tr ( \itF \otimes \itE_{i} \geng) 
	& =  \int_{\Omega}  \genF(z)^T \itE_{i} \geng(z) \dd \mu(z) 
	=   - \int_{\Omega}  \Tr \left[ \genF(z)  \geng(z)^T \itE_{i} \right]  \dd\mu(z),
\end{split}
\end{equation}
where the last equality follows from the identity $v^Tw= \Tr(v^Tw)=\Tr(vw^T)$ which is valid for any two column vectors $v, w$.
The proof is now complete recalling the definition of $\itZ_1$ in \eqref{eq:pjinfg}.
\end{proof} 

In a similar spirit, we also obtain deformation formulas with respect to the other parameters. This result will be used in Section \ref{sec:symbolic}. 

\begin{prop}\label{prop:additionaldeformation}
Under the same conditions as Proposition~\ref{lem:basic2}, 
\beq\label{diffidenz} \begin{split} 
	&	\partial_{\itl_i} \log \det(\itid-\itH) = \Tr \left[ (\itZ_1^2 -2 \itZ_2 ) \itE_{i} \right], \\
	&	\partial_{\itt_i} \log \det(\itid-\itH) = \Tr \left[  ( - \itZ_1^3+2 \itZ_1\itZ_2 + \itZ_2\itZ_1 -3\itZ_3) \itE_{i} \right],
\end{split} \eeq 
where $\itZ_1, \itZ_2$, and $\itZ_3$ are defined in \eqref{eq:pjinfg}.
\end{prop}

\begin{proof}
Since $\partial_{\geny_i} \itconjb(z) = z^2 \itE_i \itconjb(z)$ and $\partial_{\gent_i} \itconjb(z) = z^3 \itE_i \itconjb(z)$, 
we find that 
\beqq \begin{aligned}
	&\partial_{\itl_i} \genf(z)  = {z^2}  \itE_{i} \genf(z)+ \beta(z) \genf(z), \qquad 
	&&\partial_{\itl_i} \geng(z) = - {z^2}  \itE_{i} \geng(z) - \beta(z) \geng(z),\\	
	&\partial_{\itt_i} \genf(z)  = {z^3} \itE_{i} \genf(z)+  \gamma(z) \genf(z), \qquad 
	&&\partial_{\itt_i} \geng(z) =-{z^3}{} \itE_{i} \geng(z)- \gamma(z) \geng(z), 
\end{aligned} \eeqq
where $\beta(z)\deff \frac{\partial_{\itl_i} \genc(z)}{\genc(z)}$ and $\gamma(z) \deff \frac{\partial_{\itt_i} \genc(z)}{\genc(z)}$. 
As in the proof of Proposition~\ref{lem:basic2}, we find that 
\beqq \begin{split}
	&\partial_{\itl_i} \itH(z,w) =   (z+w)  \genf(z)^T \itE_{i} \geng(w) + \beta(z)\itH(z,w)-\itH(z,w) \beta(w) , \\
	&\partial_{\itt_i} \itH(z,w) =   (z^2 + zw + w^2) \genf(z)^T \itE_{i} \geng(w) + \gamma(z)\itH(z,w)-\itH(z,w) \gamma(w), 
\end{split} \eeqq
for $z\neq w$. As in the last proposition, this identity also extends to $z=w$ because both sides are equal to $0$ in this case.
Let $\kpMz$ be the multiplication by $z$, $(\kpMz \kpf)(z)= z \kpf(z)$. 
In operator forms, the equations become
\beqq \begin{split}
	&\partial_{\itl_i} \itH = \kpMz  \itE_{i} \genf\otimes   \geng +  \genf\otimes \kpMz \itE_{i} \geng + [\beta, \itH],  \\
	&\partial_{\itt_i} \itH  =   \kpMz^2 \itE_{i}  \genf\otimes   \geng +   \kpMz \genf\otimes  \kpMz \itE_{i} \geng + \genf\otimes  \kpMz^2 \itE_{i}  \geng + [\gamma, \itH].
\end{split} \eeqq
Calculating as in \eqref{eq:partHtr1}--\eqref{eq:partHtr2}, we find
\beqq
\begin{split}
	\partial_{\itl_i} \log \det(\itid- \itH)  
	&=  - \Tr ((\itid-\itH)^{-1}  \kpMz  \itE_{i}  \genf \otimes  \geng)  - \Tr ((\itid-\itH)^{-1}  \genf \otimes  \kpMz \itE_{i} \geng) , \\
	\partial_{\itt_i} \log \det(\itid- \itH)  
	&=  -\Tr ( (\itid-\itH)^{-1} \kpMz^2 \itE_{i}  \genf \otimes  \geng)  - \Tr ((\itid-\itH)^{-1}  \kpMz \genf \otimes  \kpMz \itE_{i} \geng) 
	- \Tr ((\itid-\itH)^{-1}  \genf \otimes  \kpMz^2 \itE_{i} \geng). 
\end{split}\eeqq
We now relate each such term with the matrices $\itYb_{n}$ from \eqref{eq:pjinfg}.
First, proceeding similarly to \eqref{eq:partHtr3},
\beqq
	\Tr ( (\itid-\itH)^{-1} \genf \otimes \kpMz^n \itE_{i} \geng) 
	= \Tr ( \genF \otimes \kpMz^n \itE_{i} \geng)
	=  \int_{\Omega} \Tr \left[ z^n \genF(z) \geng(z)^T \itE_i \right] \dd\mu(z)
	= \Tr \left[ \itYb_{n+1} \itE_i \right]
\eeqq
for all $n\ge 1$.
The identity $\Tr ( Ku\otimes v)=\Tr ( u\otimes (K^Tv))$ is true
for any vectors $v$ and $v$ and trace class operator $K$. Using this identity and \eqref{eq:ResolventKernel} and \eqref{eq:pjinfg}, we find
\beqq 
\begin{split} 
	\Tr ( (\itid-\itH)^{-1} \kpMz^n \itE_{i} \genf \otimes   \geng) 
	&= \Tr ( \kpMz^n \itE_{i} \genf \otimes (\itid-\itH^T)^{-1}   \geng)
	=  \int_{\Omega} \Tr \left[ z^n \itE_i \genf(z) \genG(z)^T  \right] \dd\mu(z) \\
	&=\Tr\left[ \itE_i\int_{\Omega}  z^n \genf(z) \genG(z)^T \dd\mu(z) \right] = -\Tr \left[ \itE_i (\itYb^{-1})_{n+1}  \right] . 
\end{split} \eeqq
Finally, note that 
\begin{equation*}\begin{split}
	[\kpMz, (\itid-\itH)^{-1}] 
	&= (\itid-\itH)^{-1} [  \itid-\itH, \kpMz ] (\itid-\itH)^{-1}
	=  (\itid-\itH)^{-1} [\itH, \kpMz] (\itid-\itH)^{-1} \\
	&= (\itid-\itH)^{-1} \genf \otimes   \geng (\itid-\itH)^{-1}
	= (\itid-\itH)^{-1} \genf \otimes    (\itid-\itH^T)^{-1} \geng = \itF \otimes \itG . 
\end{split} \end{equation*}
This identity implies that $(\itid-\itH)^{-1} \kpMz = \kpMz (\itid-\itH)^{-1} - \itF \otimes \itG$, and thus acting on $\genf^T$,  we find 
$$
	(\itid-\itH)^{-1} \kpMz \genf^T
	= \kpMz \itF^T - \itF^T \left(  \int_{\Omega} \itG(z) \genf(z)^T \dd z \right) 
	= \kpMz \itF^T + \itF^T (\itYb^{-1})_1^T .
$$
Hence, 
\beqq 
\begin{split}
	\Tr ((\itid-\itH)^{-1}  \kpMz \genf \otimes  \kpMz \itE_{i} \geng) 
	&= \int_{\Omega} ( z \itF(z)^T + \itF(z)^T (\itYb^{-1})_1^T) z \itE_{i} \geng(z) \dd z \\
	&= \int_{\Omega} \Tr \left[  \left( z \itF(z)  + (\itYb^{-1})_1  \itF(z) \right) z  \geng(z)^T  \itE_{i} \right]  \dd z 
	= \Tr \left[ \itYb_3 \itE_{i} + (\itYb^{-1})_1 \itYb_2 \itE_{i} \right] .
\end{split} \eeqq

Combining the above calculations, 
we obtain 
\beqq\begin{split}
	&\partial_{\itl_i} \log \det(\itid- \itH)  = \Tr \left[ ( (\itYb^{-1})_2- \itYb_2 ) \itE_i \right], \qquad 
	\partial_{\itt_i} \log \det(\itid- \itH)  =   \Tr \left[ ( ( \itYb^{-1})_3-  2\itYb_3 - (\itYb^{-1})_1 \itYb_2)  \itE_{i}   \right] .
\end{split}\eeqq
Since $\itYb(z)\sim \itI + \frac{\itYb_1}{z}  + \frac{\itYb_2}{z^2} + \frac{\itYb_3}{z^3} + \cdots$ as $z\to \infty$, we see that 
\beqq \begin{split}
	\itYb(z)^{-1} = \itI - \frac{\itYb_1}{z}  + \frac{\itYb_1^2- \itYb_2}{z^2} - \frac{\itYb_1^3- \itYb_1\itYb_2 - \itYb_2\itYb_1+\itYb_3}{z^3}+\Boh(z^{-4}),
\end{split} \eeqq
which determines $(\itYb^{-1})_{j}$ for $j=1,2,3$ in terms of $\itYb_1, \itYb_2, \itYb_3$, and we conclude \eqref{diffidenz}. 
\end{proof}

\section{A class of Fredholm determinants}\label{sec:aspectsFreddet}

The formulas \eqref{eq:Fstepdefn} and \eqref{eq:PKPZFstepdefn}  of the multi-point distributions for the KPZ and the periodic KPZ fixed points  involve certain Fredholm determinants. 
The multi-point distributions of other models such as continuous-time and discrete-time totally asymmetric simple exclusion processes (TASEPs) on the line and the ring also have similar Fredholm determinants \cites{Baik-Liu19, Liu19, Liao21}. 
The operators for these determinants all have a common structure. 
In this section, we consider a class of operators with this common structure and show that their Fredholm determinants are equal to the Fredholm determinants of integrable operators with more transparent structure. The main result is Theorem \ref{thm:algebraicmain}. 
In the special case of the KPZ and the periodic KPZ fixed point, the new integrable operators will be cubic integrable, and we prove this statement in Section ~\ref{sec:cubicKPZ}. 
In a subsequent paper, we will consider the application of Theorem \ref{thm:genintoperred} to other models in the KPZ universality class, and derive differential equations for them. 


\subsection{A class of operators}

We introduce a class of operators which arise in the study of the multi-point distributions of various models in the KPZ universality class. 

Fix a positive integer $m$. Suppose that $\Omega_{1,1}, \cdots, \Omega_{1,m}, \Omega_{2,1},\cdots, \Omega_{2,m}$ are $2m$ pairwise disjoint subsets of $\C$. 
We assume that they are either all finite unions of contours or all discrete sets. Set
\begin{equation}\label{deff:genOmegaksets}
	\Omega_1\deff\Omega_{1,1}\cup \cdots \cup \Omega_{1,m}\qquad \text{and}\qquad   \Omega_2\deff \Omega_{2,1}\cup \cdots\cup \Omega_{2,m},
\end{equation}
and
\beq \begin{cases}
	\text{$\mu$ is the counting measure if $\Omega_{\ell, k}$ are discrete,}\\  
	\text{$\dd\mu=\dd z$ if $\Omega_{\ell, k}$ are contours.}  
\end{cases} \eeq
Let $\itA_1(z), \cdots, \itA_m(z)$ and $\itB_1(z), \cdots, \itB_m(z)$ be functions of $z\in \Omega_1\cup \Omega_2$, and define the $m\times m$ diagonal matrices
\begin{equation}
\label{eq:Kbasicfuntios}
	\Am(z) = \diag \big(\itA_{j}(z) \big)_{j=1}^m \qquad\text{and}\qquad
	\Bm(z) =  \diag \big( \itB_{j}(z) \big)_{j=1}^m.
\end{equation}
Finally, suppose that $\genm_1(z), \cdots, \genm_m(z)$ are non-vanishing functions of $z\in \Omega_1\cup \Omega_2$, and set 
\begin{equation}
\label{eq:defnofMm} 
	\Mm(z) =\diag(\itM_j(z))_{j=1}^m
	\deff \diag \bigg( \genm_1(z) , \frac{\genm_1(z)}{\genm_2(z)} , \frac{\genm_3(z)}{\genm_2(z)}, \frac{\genm_3(z)}{\genm_4(z)}, \cdots \bigg) .
\end{equation}

\begin{definition} \label{defn:whK12}
With the notions \eqref{deff:genOmegaksets}--\eqref{eq:defnofMm}, define integral operators
$$
	\wh\itK_1: L^2(\Omega_2, \dd\mu) \to L^2(\Omega_1, \dd\mu) \qquad \text{and}\qquad \wh\itK_2: L^2(\Omega_1, \dd\mu) \to L^2(\Omega_2, \dd\mu)
$$
by their kernels
$$
	\wh\itK_1 (u,v) = \left(\delta_i(j)+\delta_i(j+(-1)^i)\right) 
	\frac{\itM_i(u) \itA_i(u) \itB_j(v)}{u-v}, \qquad \text{if }u\in \Omega_{1,i},\; v\in \Omega_{2,j},
$$
and
$$
	\wh\itK_2 (u,v) = \left(\delta_i(j)+\delta_i(j-(-1)^i) \right) \frac{\itM_i(u)^{-1} \itA_{i}(u) \itB_j(v) }{u-v},   \qquad 
	\text{if }u\in \Omega_{2,i},\; v\in \Omega_{1,j},
$$
for $i,j \in\{ 1, \cdots, m\}$. 
\end{definition}

The $m$-point distributions of several models of the KPZ universality class are expressible in terms of the Fredholm determinant $\det(\itid -\wh \genK_1 \wh \genK_2)$ for some $\wh \genK_1$, $ \wh \genK_2$ of the above form. 
The main differences among the models are the function $\genm_i$. 
For the KPZ and the periodic KPZ fixed points, $\genm_i$ are cubic exponential functions $e^{\gent_j z^3 + \geny_j z^2 + \genx_j z}$ as in \eqref{def:mjexpon}; see Section~\ref{sec:cubicKPZ} for details.
For other models, the cubic exponential functions are changed to 
\begin{align*}
	& e^{t_iz} z^{k_i}(z+1)^{n_i} \quad \text{for continuous-time (periodic) TASEP,}\\
	& (1+pz)^{t_i-k_i}z^{k_i}(z+1)^{n_i}  \quad \text{for discrete-time (periodic) TASEP,}  \\
	& z^{k_i}e^{-b_iz+a_i/z} \quad \text{for (periodic) PNG model},
\end{align*}
where $a_i,b_i,k_i, n_i, t_i$ and $p$ 
are parameters of the models. 
These results are found in \cite[Section~3.4]{Baik-Liu19} and \cite[Section~2.1.1]{Liu19} for the continuous-time TASEP, in \cite[Sections~2.5 and 2.6]{Liao21} for the discrete-time TASEP \cite[Sections~2.5 and 2.6]{Liao21}, and in work in progress by Tejaswi Tripathi for the PNG model. 

In the next subsection, we show that $\det(\itid -\wh \genK_1 \wh \genK_2)$ is equal to the Fredholm determinant of an integrable operator with a more convenient structure. 

\subsection{A reformulation of a class of Fredholm determinants}

Set
\begin{equation}\label{eq:deffHilbSp}
	\Omega\deff\Omega_1\cup \Omega_2 . 
\end{equation}
We introduce $m\times m$ diagonal matrices of indicator functions 
\begin{equation}\label{defn:pkmchim}
	\pkmchim_\ell=\diag\left(\chi_{\Omega_{\ell, j}}\right)_{j=1}^m  \qquad \text{for } \ell=1,2,
\end{equation}
which, due to the pairwise disjointness of the sets $\Omega_{j,k}$, satisfy 
\begin{equation}\label{eq:chimeasy}
	\chim_{\ell}(z) \ite \ite^T \chim_{\ell'} (z)=  \chim_\ell(z) \delta_{\ell, \ell'}
\end{equation}
for $\ell, \ell'=1,2$. 
Define the $2\times 2$ matrix 
$$
	\itU\deff \begin{pmatrix} 1 & 1 \\ 0 & 0 \end{pmatrix}
$$
and introduce the $m\times m$ matrices
\begin{equation}\label{def:matrLambda}
	\itLambda_1\deff
\begin{dcases}
	\diag(\underbrace{ \itU,\cdots, \itU}_{m/2}), & m \mbox{ even}, \\
	\diag(\underbrace{ \itU,\cdots, \itU}_{(m-1)/2},1), & m \mbox{ odd},
\end{dcases}
\qquad
\mbox{and}
\qquad
	\itLambda_2\deff 
\begin{dcases}
	\diag(1,\underbrace{\itU,\cdots, \itU}_{(m-2)/2},1), & m \mbox{ even}, \\
	\diag(1,\underbrace{\itU,\cdots, \itU}_{(m-1)/2}), & m \mbox{ odd}.
\end{dcases}
\end{equation}
Furthermore, define an $m\times m$ matrix, an $(m+1)\times (2m)$ projection matrix, and a column vector given by 
\beq \label{eq:itE11m}
	\itE_{11}\deff \diag (1, 0, \cdots, 0), \qquad
	\itP=\begin{pmatrix} \itI_{m+1} & \itzero \end{pmatrix}, \qquad
	\ite = (1, 1, \cdots, 1)^T\in \R^m,
\eeq
respectively. Recall the $m\times m$ diagonal matrices $\Am(z)$ and $\Bm(z)$ from \eqref{eq:Kbasicfuntios}, and define the $(m+1)$-dimensional column vector-valued functions of $z\in \Omega$,
\begin{equation} \label{eq:ABcubicadm}
	\Fm(z) 
	\deff \itP \begin{pmatrix} \itLambda_1 & \itLambda_2-\itE_{11}  \\ \itzero &  \itE_{11} \end{pmatrix}  
	\begin{pmatrix}  \Am(z)\genchi_1(z)   \ite \\  \Am(z)\genchi_2(z) \ite \end{pmatrix}, 
	\qquad
	\Gm(z) 
	\deff \itP \begin{pmatrix} \itLambda_1 & \itLambda_2-\itE_{11}  \\ \itzero &  \itE_{11} \end{pmatrix}  
	\begin{pmatrix}  \Bm(z)\genchi_2(z)   \ite \\  \Bm(z)\genchi_1(z) \ite \end{pmatrix}.
\end{equation}
Finally, introduce the $(m+1)\times (m+1)$ diagonal matrix-valued function 
\begin{equation}\label{defn:itconjbgeneral}
	\itconjb(z) \deff \diag\left(\genm_1(z),\genm_2(z),\hdots,\genm_m(z),1\right),\qquad z\in \Omega, 
\end{equation}
where $\genm_j(z)$ are the functions in \eqref{eq:defnofMm}. 

\begin{definition} \label{defn:notationsforH}
Let $\genc(z)$ be a non-vanishing scalar function of $z\in \Omega$. Define the $(m+1)$-dimensional column vector-valued functions 
\begin{equation}\label{eq:fg22} \begin{split}
	&\genf(z) = \left(\genf_1(z), \hdots,\genf_{m+1}(z)\right)^T \deff \genc(z)\itconjb(z) \Fm(z),  \\
	&\geng(z) = \left(\geng_1(z), \hdots,\geng_{m+1}(z)\right)^T \deff \frac1{\genc(z)} \itconjb(z)^{-1} \Gm(z)
\end{split}\end{equation}
for $z\in \Omega$, and the integral operator $\genH: \genHilb\to \genHilb$ by its kernel
\begin{equation}\label{eq:intoperH}
	\genH(u,v) \deff \frac{\genf(u)^T \geng(v)}{u-v} \quad \text{for }u,v\in\Omega \text{ with }u\neq v, \quad \text{and}\quad \genH(u,u)\deff 0,\quad u \in \Omega.
	\end{equation}
\end{definition}

The following result is the main outcome of this section.

\begin{thm}\label{thm:algebraicmain}
Let $\wh\genK_1$ and $\wh\genK_2$ be operators in Definition~\ref{defn:whK12}, and $\genH$ in Definition~\ref{defn:notationsforH}. If the Fredholm determinant of  $\det(\itid -\wh \genK_1 \wh \genK_2)$ is well defined as a series expansion and $\genH$ is a trace class operator, then \beq \label{eq:intoperH001}
	\det(\itid -\wh \genK_1 \wh \genK_2)  = \det(\itid-\genH).   
\eeq
Furthermore, 
\beq \label{eq:genfgengdiag}
	\genf_i(z)  \geng_i(z)=0 \quad \text{for every $z\in \Omega$ and $i\in \{1, \cdots, m+1\}$.} 
\eeq
\end{thm}

The proof of Theorem~\ref{thm:algebraicmain} is given in Section~\ref{sec:prooftraceclassthm}. 

The scalar function $\genc(z)$ in \eqref{eq:fg22} is arbitrary, and changing it does not affect the series expansion of $\det(\itid-\genH)$. However, in Section \ref{sec:cubicKPZ}, we choose it appropriately in applications so that $\genH$ becomes a trace class operator.

The main difference of $\genH$ from $\wh\genK_1$ and $\wh\genK_2$ is that the matrix $\itconjb(z)$ does not involve the ratios of $\genm_i(z)$ unlike $\Mm(z)$. 
Thus, the above identity ``uncouples" $\genm_i(z)$s in the kernel formula. 
This is useful when we derive differential equations as in Section~\ref{sec:IIKSRHP}.

\bigskip

The formulation \eqref{eq:fg22} is useful in deriving differential equations. 
However, for other purposes, namely in proving that the operator is trace class and in obtaining asymptotic formulas, we find that alternative representations of $\genf$ and $\geng$ are more suitable, and we state them next.

\begin{lem} \label{lem:fgaltformula}
Define the $m\times m$ matrix 
\begin{equation}\label{def:Sm}
\begin{split}
	&\Sm(z)\deff \diag \bigg( \genm_1(z)^{1/2}, \frac{\genm_1(z)^{1/2}}{\genm_2(z)^{1/2}} , \frac{\genm_3(z)^{1/2}}{\genm_2(z)^{1/2}},  \frac{\genm_3(z)^{1/2}}{\genm_4(z)^{1/2}}, \frac{\genm_5(z)^{1/2}}{\genm_4(z)^{1/2}}, \cdots \bigg),\quad z\in \Omega,
\end{split} 
\end{equation}
for some fixed branch of the square-root. 
Choose the conjugating constant as 
\begin{equation}\label{eq:constantcL2} \begin{split}
	\genc(z) 
	=  \ite^T  (\genchi_1(z)  +\genchi_2(z)  ) 
	 	\diag \bigg( \frac1{\sqrt{\genm_1(z) }},  \frac1{\sqrt{ \genm_1(z)\genm_2(z)}} , 
	 \cdots, \frac1{\sqrt{ \genm_{m-1}(z)\genm_m(z)}}  \bigg) \ite . 
\end{split} \end{equation}
Then the functions $\genf$ and $\geng$ in \eqref{eq:fg22} can be written as 
\begin{equation}\label{eq:alternativefg}
\begin{split}
	& \genf(z)=  \itP 
	\begin{pmatrix} \itLambda_1 & \itLambda_2-\itE_{11}  \\ 0 &  \itE_{11} \end{pmatrix} 
	\begin{pmatrix} \Am(z)  \Sm(z) \genchi_1(z) & \itzero  \\ \itzero & \Am(z) \Sm(z)^{-1}  \genchi_2(z)  \end{pmatrix}
	\begin{pmatrix} \ite  \\ \ite  \end{pmatrix},  \\
	& \geng(z)= \itP 
	\begin{pmatrix} \itLambda_1 & \itLambda_2-\itE_{11}  \\ 0 &  \itE_{11} \end{pmatrix} 
	\begin{pmatrix} \Bm(z)  \Sm(z)^{-1} \genchi_2(z) & \itzero  \\ \itzero & \Bm(z) \Sm(z) \genchi_1(z)  \end{pmatrix}
	\begin{pmatrix} \ite  \\ \ite  \end{pmatrix}. 
\end{split}
\end{equation}
\end{lem}

\medskip

Before proving the lemma, we first  introduce the $m\times m$ diagonal matrices
\begin{equation}\label{def:itGoe}
\begin{split}
	&\Gmo(z)\deff \diag ( \genm_1(z), \genm_1(z), \genm_3(z), \genm_3(z), \genm_5(z), \cdots), \\
	&\Gme(z) \deff \diag ( 1, \genm_2(z), \genm_2(z), \genm_4(z), \genm_4(z), \cdots) , 
\end{split} 
\end{equation}
which we will also use later in this section. 

\begin{proof}[Proof of Lemma~\ref{lem:fgaltformula}]
Define the $m\times m$ matrix $\widetilde \itconjb \deff \diag\left(\genm_1,\genm_2,\hdots,\genm_m \right)$. 
Note that
$$
	\itconjb(z)^{\pm 1}  \itP= \begin{pmatrix}  \itconjb(z)^{\pm 1} & \itzero \end{pmatrix} 
	= \itP \begin{pmatrix}  \itconjb(z)^{\pm 1} & \itzero \\ \itzero & \itI_{m-1} \end{pmatrix} 
	= \itP \begin{pmatrix}  \widetilde \itconjb(z)^{\pm 1} & \itzero \\ \itzero & \itI_{m} \end{pmatrix}.
$$
It is straightforward to check that 
$$
\widetilde \itconjb(z)^{\pm 1} \itLambda_1  = \itLambda_1  \itGo(z)^{\pm 1}, \qquad 
\widetilde \itconjb(z)^{\pm 1} ( \itLambda_2 - \itE_{11})  = ( \itLambda_2 - \itE_{11})   \itGe(z)^{\pm 1}, \\
$$
and
$\itE_{11}= \itE_{11}\itGe(z) =\itE_{11}  \itGe(z)^{-1}$. 
Applying these identities to \eqref{eq:fg22}, we find that  
\begin{equation}\label{eq:fginterm1}
\begin{split}
	& \genf(z)=  \itP 
	\begin{pmatrix} \itLambda_1 & \itLambda_2-\itE_{11}  \\ 0 &  \itE_{11} \end{pmatrix} 
	\begin{pmatrix} \Gmo(z) \Am(z) \genchi_1(z)\ite  \\ \Gme(z) \Am(z) \genchi_2(z)\ite  \end{pmatrix} \genc(z),  \\
	& \geng(z)= \itP 
	\begin{pmatrix} \itLambda_1 & \itLambda_2-\itE_{11}  \\ 0 &  \itE_{11} \end{pmatrix} 
	\begin{pmatrix} \Gmo(z)^{-1} \Bm(z) \genchi_2(z)\ite  \\ \Gme(z)^{-1} \Bm(z) \genchi_1(z)\ite \end{pmatrix} \frac1{\genc(z)}  . 
\end{split} \end{equation}
In terms of \eqref{def:itGoe}, the conjugating constant chosen in \eqref{eq:constantcL2} can be written as  
\beqq \begin{split}
	\genc(z) &= \ite^T  (\genchi_1(z)  +\genchi_2(z)  ) \Gmo(z)^{-1/2} \Gme(z)^{-1/2} \ite, 
\end{split} \eeqq
and, due to the indicator functions,  it satisfies
\beqq
	\frac1{\genc(z)}= \ite^T  (\genchi_1(z)  +\genchi_2(z)  ) \Gmo(z)^{1/2} \Gme(z)^{1/2} \ite . 
\eeqq
Now using \eqref{eq:chimeasy}, we find that 
\beqq \begin{split}
	&\genchi_\ell (z)\ite  \genc(z)
	= \genchi_\ell (z)\ite  \ite^T  (\genchi_1(z)  +\genchi_2(z)  ) \Gmo(z)^{-1/2} \Gme(z)^{-1/2} \ite
	= \genchi_\ell (z)\Gmo(z)^{-1/2} \Gme(z)^{-1/2} \ite, \\
	&\genchi_\ell (z)\ite  \frac1{\genc(z)}
	= \genchi_\ell (z)\ite  \ite^T  (\genchi_1(z)  +\genchi_2(z)  ) \Gmo(z)^{1/2} \Gme(z)^{1/2} \ite
	= \genchi_\ell (z)\Gmo(z)^{1/2} \Gme(z)^{1/2} \ite
\end{split} \eeqq
for $\ell=1, 2$.
Thus, using $\Sm(z)= \Gmo(z)^{1/2} \Gme(z)^{-1/2}$, we find, for example, 
\beqq
	\Gmo(z) \Am(z) \genchi_1(z)\ite \genc(z)  
	= \Am(z) \genchi_1 (z)\Gmo(z)^{1/2} \Gme(z)^{-1/2} \ite
	= \Am(z)  \Sm(z) \genchi_1(z) \ite. 
\eeqq
The other three vectors in \eqref{eq:fginterm1} are computed in a similar way.
\end{proof}

\subsection{Proof of Theorem \ref{thm:algebraicmain}}\label{sec:prooftraceclassthm}

In this Subsection we prove \ref{thm:algebraicmain}. 
Replacing the kernels in Definition \ref{defn:whK12} by 
\beqq
	\chi_{\Omega_1}(u) \wh \itK_1(u,v)  \chi_{\Omega_2}(v) 
	\quad \text{and} \quad
	\chi_{\Omega_2}(u) \wh \itK_2(u,v)  \chi_{\Omega_1}(v) 
\eeqq
for $u, v\in \Omega$, we can extend $\wh \itK_1$ and $\wh \itK_2$ to operators from $\genHilb$ to $\genHilb$. We use the same notation 
$\wh \itK_1$ and $\wh \itK_2$ to denote the extended operators. The Fredholm determinant $\det(\itid -\wh \genK_1 \wh \genK_2)$ is unchanged by this modification. 

\begin{lem}\label{lem:kpzmatrixop}
For each $\ell=1, 2$, the kernel of $\wh\itK_\ell$ can be written as 
\beq \label{eq:whKkpztap}
	\wh\itK_\ell(u,v)=\frac{\wh \kpqq_\ell(u)^T \wh \kprr_\ell(v)}{u-v}, \qquad u, v\in \Omega, 
\eeq
with the $m\times 1$ column vectors 
\begin{equation*}
\begin{aligned}
	&\wh \gena_1(z)= \itLambda_1 \Am(z)\Mm(z)\genchi_{1} (z)\ite, 	 
	&&\wh \genb_1(z)= \itLambda_1 \Bm(z)\genchi_2(z) \ite, \\
	& \wh \gena_2(u)= \itLambda_2 \Am(u)\Mm(z)^{-1} \genchi_{2} (u)\ite,  \qquad 
	&& \wh \genb_2(z)= \itLambda_2 \Bm(z)\genchi_1(z) \ite.
\end{aligned}
\end{equation*}
\end{lem}

\begin{proof}
Set $\chi_{\Omega_{\ell,0}}=\chi_{\Omega_{\ell,m+1}}=0$ for $\ell=1,2$ and $\itB_0=\itB_{m+1}=0$. 
Due to indicator functions, the kernels in Definition \ref{defn:whK12} can be written as
\begin{align*}
	(u-v) \wh\itK_1(u,v)& =
\begin{multlined}[t]
	\sum_{\substack{k=1 \\ k\text{ odd} }}^m \chi_{\Omega_{1,k}}(u) \itM_k(u) \itA_k(u) \left(\itB_k(v)  \chi_{\Omega_{2,k}}(v) +\itB_{k+1}(v) \chi_{\Omega_{2,k+1}}(v)\right)  \\
	+ \sum_{\substack{k=1 \\ k\text{ even} }}^m \chi_{\Omega_{1,k}}(u) \itM_k(u) \itA_k(u)  \left(\itB_k(v) \chi_{\Omega_{2,k}}(v) +\itB_{k-1}(v) \chi_{\Omega_{2,k-1}}(v) \right)
\end{multlined}
\end{align*}
and
\begin{align*}
	(u-v) \wh\itK_2(u,v)& =
\begin{multlined}[t]
	\sum_{\substack{k=1 \\ k\text{ odd} }}^m \chi_{\Omega_{2,k}}(u) \itM_k(u)^{-1} \itA_k(u) \left(\itB_k(v) \chi_{\Omega_{1,k}}(v)  +\itB_{k-1}(v)  \chi_{\Omega_{1,k-1}}(v) \right)  \\
	+ \sum_{\substack{k=1 \\ k\text{ even} }}^m \chi_{\Omega_{2,k}}(u) \itM_k(u)^{-1} \itA_k(u) \left(\itB_k(v) \chi_{\Omega_{1,k}}(v) +\itB_{k+1}(v) \chi_{\Omega_{1,k+1}}(v) \right).
\end{multlined}
\end{align*}
It is direct to check that for any column vectors $\alpha=(\alpha_1,\hdots,\alpha_m)^T\in \R^m $ and $\beta=(\beta_1,\hdots,\beta_m)^T\in \R^m $, 
\beqq 
	\alpha^T \itLambda_1^T \itLambda_1 \beta
	=  \sum_{\substack{k=1 \\ k\text{ odd} }}^m \alpha_k (\beta_k+\beta_{k+1}) 
	+  \sum_{\substack{k=1 \\ k\text{ even} }}^m \alpha_k (\beta_k+\beta_{k-1})
\eeqq
and
\beqq 
	 \alpha^T \itLambda_2^T \itLambda_2 \beta
	=  \sum_{\substack{k=1 \\ k\text{ odd} }}^m \alpha_k (\beta_k+\beta_{k-1}) 
	+  \sum_{\substack{k=1 \\ k\text{ even} }}^m \alpha_k (\beta_k+\beta_{k+1}), 
\eeqq
where we set $\beta_0=\beta_{m+1}=0$. 
The result then follows.
\end{proof}

For the next lemma, recall the matrices $\itGo$ and $\itGe$ defined in \eqref{def:itGoe}. Also introduce $m\times m$ diagonal matrices 
\begin{equation}\label{eq:itconjdef}
	\genD_1(z) \deff \diag ( \genm_1(z) , 1 , \genm_3(z), 1, \genm_5(z), \cdots), 	
	\qquad
	\genD_2(z) \deff \diag ( 1, \genm_2(z) , 1 , \genm_4(z), 1, \cdots).
\end{equation}

\begin{lem}\label{lem:abformulaD}
Define the $m\times 1$ vector-valued functions
\begin{equation} \label{eq:kpzabvectors} 
\begin{split}
	 &\gena_1 (z) \deff  \itLambda_1 \Am(z) \Gmo(z)  \genchi_1 (z)  \ite =  \genD_1(z) \itLambda_1 \Am(z)\genchi_1(z)   \ite, \\
	 & \gena_2(z)  \deff \itLambda_2 \Am(z)\Gme(z) \genchi_{2} (z)\ite = \genD_2(z) \itLambda_2  \Am(z)\genchi_2(z) \ite,   \\
	 & \genb_1(z) \deff \itLambda_1 \Bm(z) \Gmo(z)^{-1}  \genchi_2(z) \ite =  \genD_1(z)^{-1} \itLambda_1 \Bm(z)\genchi_2(z)   \ite, \\
	 &\genb_2(z) \deff  \itLambda_2 \Bm(z) \Gme(z)^{-1} \genchi_1(z) \ite =  \genD_2(z)^{-1} \itLambda_2  \Bm(z)\genchi_1(z) \ite, 
\end{split} 
\end{equation}
for $z\in \Omega$. 
Then, there is a non-vanishing scalar function $\wh \genc$ on $\Omega$ such that for each $\ell=1,2$, 
\begin{equation*}
	\wh\gena_\ell(z)= \wh\genc(z) \gena_\ell(z)\quad \text{and}\quad \wh\genb_\ell(z)=\frac1{\wh\genc(z)} \genb_\ell(z)  .
\end{equation*}
\end{lem}

\begin{proof}
Define the scalar function
$$
	\wh \genc(z)= \ite^T ( \Gme(z)^{-1}  \genchi_1(z)  +  \Gmo(z)^{-1} \genchi_2(z)  ) \ite.
$$
Due to the indicator functions, we find that 
$$
	\frac1{\wh \genc(z)}= \ite^T ( \Gme(z)  \genchi_1(z)  +  \Gmo(z) \genchi_2(z)  ) \ite.
$$

Observe that \begin{equation}
\label{eq:Ueasy}
	\begin{pmatrix} 1 & 1 \\ 0 & 0 \end{pmatrix} \begin{pmatrix} a & 0 \\ 0 & a \end{pmatrix}
	= \begin{pmatrix} a & 0 \\ 0 & * \end{pmatrix} \begin{pmatrix} 1 & 1 \\ 0 & 0 \end{pmatrix}
\end{equation}
holds for any complex numbers $a$ and $*$. Using this identity, we find that  
$$
	\itLambda_1 \Gmo(z)^{\pm 1} = \genD_1(z)^{\pm 1}  \itLambda_1, \qquad
	\itLambda_2 \Gme(z)^{\pm 1} = \genD_2(z)^{\pm 1} \itLambda_2 .
$$
Since the matrices $\Am, \Bm ,\Gmo, \Gme$, and $\genchi_\ell$ are all diagonal, they commute. 
Since $\Mm=\Gmo\Gme^{-1}$, we find, using \eqref{eq:chimeasy}, that 
\beqq \begin{split} 
	\wh \gena_1(z) \frac1{\wh \genc(z)}  &= \itLambda_1 \Am(u) \Gmo(z)\Gme(z)^{-1}\genchi_{1} (u)\ite \ite^T ( \Gme(z)  \genchi_1(z)  +  \Gmo(z) \genchi_2(z)  ) \ite \\
	&= \itLambda_1 \Am(z) \Gmo(z)  \genchi_1 (z)  \ite = \genD_1(z) \itLambda_1 \Am(z) \genchi_1 (z) \ite = \gena_1(z) .
\end{split}  \eeqq 
The results for the other vector functions $\kprr_1$, $\kpqq_2$ and $\kprr_2$ follow similarly.
\end{proof}

For the next result, observe that \eqref{eq:kpzabvectors} and \eqref{eq:chimeasy} imply that  
\beq \label{eq:a1b122zero}
	\genb_1(z)\gena_1(z)^T=\genb_2(z)\gena_2(z)^T=0. 
\eeq

\begin{lem}\label{prop:genintoper}
Let $\gena_1, \gena_2, \genb_1, \genb_2$ be the vectors defined in \eqref{eq:kpzabvectors}.
Let $\genc(z)$ be an arbitrary non-vanishing scalar function on $\Omega$. 
For $\ell=1,2,$ define the operator $\genK_\ell: \genHilb\to \genHilb$ by its kernel
$$
	\genK_\ell(u,v)\deff   \frac{\genc(u) \gena_\ell(u)^T  \genb_\ell(v) \frac1{\genc(v)}}{u-v}  \quad 
	\text{for } u,v\in \Omega \text{ with }u\neq v, 
	\qquad \text{and}\qquad \genK_\ell(u,u)\deff 0,\quad u\in \Omega.
$$
If the series definition of $\det(\itid -\wh \genK_1 \wh \genK_2)$ is well-defined and $\genK_1$ and $\genK_2$ are trace class operators, then 
\beq \label{eq:intopersumK1K2}
	\det(\itid -\wh \genK_1 \wh \genK_2)  = \det(\itid-\genK_1-\genK_2).   
\eeq
\end{lem}

\begin{proof} 
Lemma~\ref{lem:kpzmatrixop} and~\ref{lem:abformulaD} show that for each $\ell=1,2$, 
\beqq
	\wh\genK_\ell(u,v)  =  \frac{\wh \genc(u)}{\genc(u)} \genK_\ell(u,v) \frac{\genc(v)}{\wh\genc(v)} .
\eeqq
Hence, from the series definition of Fredholm determinants, 
$\det(\itid -\wh \genK_1 \wh \genK_2) = \det(\itid - \genK_1 \genK_2)$. 

Equation~\eqref{eq:a1b122zero} implies that $\genK_\ell^2=0$ for $\ell=1,2$, and thus
$$
	\itid-\genK_1\genK_2 = (\itid + \genK_1) (\itid- \genK_1-\genK_2) (\itid+\genK_2).
$$
Now for $\ell=1,2$, we have $\Tr\genK_\ell=0$ because $\genK_\ell(u,u)=0$, and  $\Tr\genK_\ell^n=0$ for all $n\geq 2$ since $\genK_\ell^2=0$. 
Thus, the Plemelj-Smithies formula for Fredholm determinants, which is valid for trace class operators, implies that $\det(\itid+\genK_\ell)=1$.
Therefore, $\det ( \itid-\genK_1\genK_2)= \det (\itid- \genK_1-\genK_2)$, and we obtain the result. 
\end{proof}

We are ready to prove Theorem \ref{thm:algebraicmain}.
Introduce the $(2m)\times (2m)$ permutation matrix 
$$
	\itperm\deff \begin{pmatrix} \itperm_+ & \itperm_- \\ \itperm_- & \itperm_+ \end{pmatrix}\quad \text{where}\quad
	\itperm_+ \deff \diag ( 1,0,1,0, \cdots ),\quad 
	\itperm_- \deff \diag (0,1,0,1,  \cdots  ).
$$
The diagonal matrices $\itperm_\pm$ are of size $m\times m$. 
Recall the matrices $\itLambda_1$ and $\itLambda_2$ from  \eqref{def:matrLambda} and the matrix 
$\itE_{11}=\diag (1, 0, \cdots, 0)$ from \eqref{eq:itE11m}. 
It is easy to check that 
\beq \label{eq:permrelationPiLambda}
 	\itperm \begin{pmatrix}
	 \itLambda_1 &\itzero \\ \itzero &  \itLambda_2   \end{pmatrix}
 	= \begin{pmatrix}  \itLambda_1 & \itLambda_2-\itE_{11}  \\ \itzero & \itE_{11}  \end{pmatrix} .
\eeq

\begin{proof}[Proof of Theorem \ref{thm:algebraicmain}]
Let $\genK_1$ and $\genK_2$ be the operators in Lemma \ref{prop:genintoper}. 
We first show that 
\beqq
	\genK_1+\genK_2=\genH.
\eeqq
We need to show that  
\beqq
	\genc(u) \gena(u)^T \genb(v) \frac1{\genc(v)} = \genf(u)^T \geng(v),
\eeqq
where we denote $\gena(z)\deff (\gena_1(z), \gena_2(z))^T$ and $\genb(z)\deff (\genb_1(z), \genb_2(z))^T$. 
Since
$\itP^T\itP= \diag(\itI_m,\itE_{11}) $, 
we find from \eqref{eq:permrelationPiLambda} that 
\beqq
	\itperm\begin{pmatrix} \itLambda_1& \itzero \\ \itzero &\itLambda_2 \end{pmatrix}
	=\itP^T\itP \begin{pmatrix}  \itLambda_1 & \itLambda_2-\itE_{11}  \\ \itzero & \itE_{11}\end{pmatrix}
	=\itP^T\itP\itperm\begin{pmatrix} \itLambda_1&\itzero\\ \itzero&\itLambda_2 \end{pmatrix}.
\eeqq
Hence, from the first formulas of $\genb_1$ and $\genb_2$ in \eqref{eq:kpzabvectors}, $\itP^T\itP\itperm\genb(z) = \itperm \genb(z)$. 
Using $\itperm^T \itperm=\itI_{2m}$, we find
\beqq
	\gena(u)^T\genb(v) 
	= \gena(u)^T\itperm^T\itperm\genb(v)
	= \gena(u)^T\itperm^T\itP^T\itP\itperm\genb(v) 
	= (\itP \itperm \gena(u) )^T \itP\itperm\genb(v) .
\eeqq

We now insert the second formulas of $\gena_1$ and $\gena_2$ in \eqref{eq:kpzabvectors} to compute $\itP \itperm \gena(u)$. 
The relations
$$
\itperm_+\genD_1(z)=\itconjb(z)\itperm_+,\quad \itperm_-\genD_1(z)=\itperm_-, \quad \itperm_+\genD_2(z)=\itperm_+\quad \text{and}\quad \itperm_-\genD_2(z)=\itconjb(z)\itperm_-
$$
are straightforward, and yield
\beqq
	\itP \itperm\begin{pmatrix} \genD_1(z)&\itzero \\ \itzero &\genD_2(z) \end{pmatrix} \itperm
	= \itconjb(z) \itP \itperm. 
\eeqq
Thus, using \eqref{eq:permrelationPiLambda} again, we obtain 
\begin{equation}\label{eq:fcarelation}	
\itP\itperm \gena(z)
	=  \itconjb(z) \itP \itperm \begin{pmatrix} \itLambda_1 \Am(z)\genchi_1(z)   \ite \\ \itLambda_2  \Am(z)\genchi_2(z) \ite \end{pmatrix}
	= \itconjb(z) \itP \begin{pmatrix}  \itLambda_1 & \itLambda_2-\itE_{11}  \\ \itzero & \itE_{11}  \end{pmatrix}
	\begin{pmatrix}  \Am(z)\genchi_1(z)   \ite \\  \Am(z)\genchi_2(z) \ite \end{pmatrix} 
	= \frac1{\genc(z)} \genf(z) . 
\end{equation}
A similar computation yields 
\begin{equation}\label{eq:gcbrelation}
	\itP\itperm\genb(z) 
	=  \itconjb(z)^{-1} \itP \itperm \begin{pmatrix} \itLambda_1 \Bm(z)\genchi_2(z)   \ite \\ \itLambda_2  \Bm(z)\genchi_1(z) \ite \end{pmatrix} 
	= \itconjb(z)^{-1} \itP \begin{pmatrix}  \itLambda_1 & \itLambda_2-\itE_{11}  \\ \itzero & \itE_{11}  \end{pmatrix}
	\begin{pmatrix} \Bm(z)\genchi_2(z)   \ite \\  \Bm(z)\genchi_1(z) \ite \end{pmatrix} 
	= \genc(z) \geng(z). 
\end{equation}
Hence, $\genK_1+\genK_2=\genH$. 

Since
\beqq
	\chi_{\Omega_1}(u) \genH(u,v) \chi_{\Omega_2}(v) = \genK_1(u,v) , \qquad 
	\chi_{\Omega_2}(u) \genH(u,v) \chi_{\Omega_1}(v) = \genK_2(u,v) , 
\eeqq
if we assume that $\genH$ is a trace class operator, then $\genK_1$ and $\genK_2$ are trace class operators. 
Thus, Lemma \ref{prop:genintoper} applies and we obtain 
\beqq
	\det(\itid -\wh \genK_1 \wh \genK_2)  = \det(\itid-\genK_1-\genK_2) = \det(\itid-\genH),
\eeqq
and this proves \eqref{eq:intoperH001}.

Finally, we check that $\genf_i(z)\geng_i(z)=0$ for $i\in \{1, \cdots, m+1\}$, as claimed in \eqref{eq:genfgengdiag}. 
Since $\itE_i$ is the $(m+1)\times (m+1)$ diagonal matrix which has $1$ on the $i$-th diagonal entry and the other entries are all zero, 
we have $\genf_i(z)  \geng_i(z)= (\genf(z)^T \itE_i \geng(z))_{ii}$. 
Inserting the formula \eqref{eq:fg22} of $\genf$ and $\geng$ and noting that $\itconjb(z) \itE_i \itconjb(z)^{-1}= \itE_i$, we find that $\genf(z)^T \itE_i \geng(z)$ is equal to, with certain diagonal matrices $\alpha$, $\beta$, determined from $\itP^T \itE_i\itP$, 
\beqq
	\begin{pmatrix}  \genchi_1(z)  \itLambda_1^T \alpha \itLambda_1 \genchi_2(z) 
	& \genchi_1(z)  \itLambda_1^T \alpha ( \itLambda_2 - \itE_{11})  \genchi_1(z)    \\  
	\genchi_2(z) ( \itLambda_2^T - \itE_{11}) \alpha  \itLambda_1  \genchi_2(z)  
	& \genchi_2(z) ( \itLambda_2^T - \itE_{11}) \alpha ( \itLambda_2 - \itE_{11}) \genchi_1(z) + \genchi_2(z)  \itE_{11}^T \beta \itE_{11} \genchi_1(z)    \end{pmatrix}
\eeqq
multiplied by the row vector $((\Am(z) \ite)^T ,  ( \Am(z) \ite) ^T)$ on the left and by the column vector $(\Bm(z) \ite ,   \Bm(z) \ite) ^T$
on the right. 
It is direct to check that $ \itLambda_1^T \alpha ( \itLambda_2 - \itE_{11})= \itzero$ and $( \itLambda_2^T - \itE_{11}) \alpha ( \itLambda_2 - \itE_{11})=\itzero$ for every diagonal matrix $\alpha$. Thus, the off-diagonal blocks are zero. 
On the other hand, due to the characteristic functions, the diagonal blocks are also zero. 
Thus, $\genf_i(z)^T  \geng_i(z)=0$. 
\end{proof}

\subsection{A lemma on trace class operators} \label{sec:lemmatrace}

In order to apply Theorem~\ref{thm:algebraicmain}, we need to check that $\genH$ is a trace class operator. 
We will use the following lemma in the next section. We use $|v|$ to denote the usual Euclidean norm of a vector $v$.

\begin{lem} \label{lem:traceclasslemma}
Let $\genH$ be the operator in Definition~\ref{defn:notationsforH}. Suppose that each set $\Omega_{\ell}$ can be split into a finite union of disjoint subsets
\begin{equation}\label{eq:decompositionOmega}
\Omega_{\ell}=\bigcup_{j=1}^{N_{\ell}} \Sigma_{\ell,j}
\end{equation}
satisfying the following property: For every pair of indices $(j_1,j_2)$ with $j_\ell\in \{1,\hdots,N_\ell\}$, $\ell\in \{1,2\}$, there is a simple contour $C_{j_1,j_2}$ such that (a) it is either closed or extends to infinity, (b) it separates $\Sigma_{1,j_1}$ and $\Sigma_{2,j_2}$ in the sense that they are in different regions of $\C\setminus C_{j_1,j_2}$, and (c) the estimates 
\begin{equation} \label{eq:trcondition}
	\int_{\Sigma_{\ell, j_\ell}} \int_{C_{j_1,j_2}}  \frac{ |\genf(u)|^2}{ |s-u|^2 } |\dd s| |\dd \mu(u)| <\infty 
	\qquad \text{and}\qquad 
	\int_{\Sigma_{\ell, j_\ell}} \int_{C_{j_1,j_2}}   \frac{ |\geng(u)|^2}{ |s-u|^2 } |\dd s| |\dd \mu(u)|  <\infty
\end{equation}
hold for $\ell=1,2$. 
Then, $\genH$ is a trace class operator. 

Furthermore, there is a union $C$ of contours that is disjoint from $\Omega$ such that $\genH$ can be written as $\genH=\genT_2\genT_1$ for Hilbert-Schmidt operators
$\genT_1:\genHilb\to L^2(C,\dd z)$, $\genT_2: L^2(C,\dd z)\to \genHilb$ with Hilbert-Schmidt norms 
$$
\|\genT_1\|_2 \le \left(\int_\Omega \int_C \frac{ |\genf(u)|^2}{ |s-u|^2 } |\dd s| |\dd \mu(u)|\right)^{1/2}\qquad\text{and}\qquad \|\genT_2\|_2 \le \left(\int_\Omega \int_C \frac{ |\geng(u)|^2}{ |s-u|^2 } |\dd s| |\dd \mu(u)|\right)^{1/2} . 
$$
\end{lem}

The decomposition $\genH=\genT_2\genT_1$ will be used only in Lemma~\ref{lem:asymptfreddet}. 

\begin{proof}
We express the kernel $\genH$ in the form 
$$
	\genH(u,v)= \sum_{j_1,j_2} \chi_{\Sigma_{1,j_1}}(u) \genH(u,v) \chi_{\Sigma_{2,j_2}}(v)
	+  \sum_{j_1,j_2} \chi_{\Sigma_{2,j_2}}(u) \genH(u,v) \chi_{\Sigma_{1,j_1}}(v).
$$
In this decomposition terms of the form $\chi_{\Sigma_{\ell,j_\ell}}(u) \genH(u,v) \chi_{\Sigma_{\ell,j'_\ell}}(v)$, $\ell\in\{1, 2\}$, do not appear due to the structure of the kernel and the indicator functions in \eqref{eq:ABcubicadm} and \eqref{eq:intoperH}. 
To show that $\genH$ is trace class, it is enough to prove that each term in the sum defines a trace class operator. 
We only show it for the terms in the first sum since for the terms in the second sum the proof is analogous. 

For $u\in \Sigma_{1,j_1}$ and $v\in \Sigma_{2,j_2}$, the separation assumption and the Cauchy integral formula imply that
$$
	\frac1{u-v} = \frac{\pm 1}{2\pi \ii} \int_{C_{j_1,j_2}} \frac{\dd s}{(s-u)(s-v)}
$$
where the sign of the term $\pm 1$ depends on the relative location of $\Sigma_{1,j_1}$ to $C_{j_1,j_2}$. 
Thus,
$$
	 \chi_{\Sigma_{1,j_1}}(u) \genH(u,v) \chi_{\Sigma_{2,j_2}}(v)= 
	\frac{\pm 1}{2\pi \ii}  \int_{C_{j_1,j_2}} \frac{\chi_{\Sigma_{1,j_1}}(u) \genf(u)^T  \geng(v) \chi_{\Sigma_{2,j_2}}(v)}{(s-u)(s-v)} \dd s, 
	 \qquad u,v\in \Omega.
$$
This is the kernel of the product of two operators, one from $L^2(\Omega,\dd \mu)$ to $L^2(C_{j_1,j_2}, \dd z)$ and the other in the opposite way. 
Assumption~\eqref{eq:trcondition} implies that these two operators are Hilbert-Schmidt, and thus $\genH$ is a trace class operator. 

Furthermore, taking the $C$ to be the union of $C_{j_1, j_2}$, we find the decomposition $\genH=\genT_2\genT_1$, with $\genT_1$, $\genT_2$ being Hilbert-Schmidt operators with the claimed Hilbert-Schmidt norm bounds. 
\end{proof}


We note that from the formula \eqref{eq:alternativefg} and the definition of $\Omega$, the condition~\eqref{eq:trcondition} is satisfied if 
\begin{enumerate}[(i)]
\item $\Am, \Bm\in L^\infty (\Omega,\dd\mu)$, and 
\item for every $\ell_1, \ell_2 \in\{1, 2\}$ with $\ell_1\neq \ell_2$, and every $j_1,j_2$,
\beq \label{eq:condb1}
	\int_{\Omega_{\ell_k,j_{\ell_k}}} \int_{C_{j_1,j_2}}   \left|\frac{\sqrt{\genm_{j+1}(u)}}{\sqrt{\genm_j(u)}}\right|  \frac{ |\dd s| |\dd \mu(u)| }{  |s-u|^2 }  <\infty 
	 \quad \text{when $\ell_k-j_{\ell_k}$ is even,}
\eeq
and
\beq \label{eq:condb2}
	\int_{\Omega_{\ell_k,j_{\ell_k}}} 
	\int_{C_{j_1,j_2}} \left|\frac{\sqrt{\genm_{j}(u)}}{\sqrt{\genm_{j+1}(u) }}\right|  \frac{ |\dd s| |\dd \mu(u)| }{ |s-u|^2 }   <\infty  
	\quad \text{when $\ell_k-j_{\ell_k}$ is odd.} 
\eeq
\end{enumerate}
%

\section{Cubic integrable operators for the KPZ and the periodic KPZ fixed points}  \label{sec:cubicKPZ}

We state the Fredholm determinants \eqref{eq:Fstepdefn} and \eqref{eq:PKPZFstepdefn} obtained in \cite{Liu19} and \cite{Baik-Liu19} explicitly, and then apply Theorem~\ref{thm:algebraicmain} of the last section to express them in terms of new Fredholm determinants which we show to be cubic admissible determinants. 
We also prove a symmetry for the case of the KPZ fixed point, thus proving 
Theorem~\ref{thm:genintoperred}.

\subsection{KPZ fixed point - proof of Theorem~\ref{thm:genintoperred} (i)} \label{sec:Liuformula}

The paper \cite[Definition 2.23]{Liu19} showed that the function in \eqref{eq:Fstepdefn} is given by  
\begin{equation}\label{deff:itDBaikLiu}
	\itD (\kph,\kpga,\kpt\mid \oz) \deff \det( \itid - \wh\itK_1^\kpz \wh\itK_2^\kpz) 
\end{equation}
where the Fredholm determinant is defined through an absolutely convergent series definition and the operators $ \wh\itK_1^\kpz$ and $\wh\itK_2^\kpz$ are of the form in Definition \ref{defn:whK12} which we now describe. 
The above formula is valid for the parameters satisfying
\beq \label{eq:KPZparameco}
	0<\kpt_1\le \cdots \le \kpt_m \quad \text{where} \quad \text{$\kpga_i<\kpga_{i+1}$ if $\kpt_i=\kpt_{i+1}$,}
\eeq
which we can always assume by re-labeling the parameters. 
Note that the complex numbers $\zeta_1, \cdots, \zeta_{m-1}$ are fixed complex numbers of different moduli satisfying $0<|\zeta_i|<1$ for all $i$. 

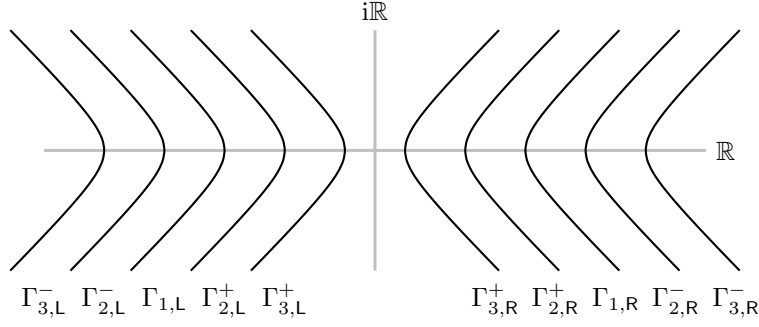
\begin{figure}\centering
\begin{tikzpicture}[scale=0.4]
\draw [line width=0.4mm,lightgray] (-11,0)--(11,0) node [pos=1,right,black] {$\R$};
\draw [line width=0.4mm,lightgray] (0,-4)--(0,4) node [pos=1,above,black] {$\ii\R$};
\draw[domain=-4:4,smooth,variable=\y, thick]  plot ({-(\y*\y+1)^(1/2)},{\y}) node at (-3,-5) {$\Gamma_{3, \kpL}^+$};
\draw[domain=-4:4,smooth,variable=\y, thick]  plot ({-2-(\y*\y+1)^(1/2)},{\y}) node at (-5,-5) {$\Gamma_{2, \kpL}^+$};
\draw[domain=-4:4,smooth,variable=\y, thick]  plot ({-4-(\y*\y+1)^(1/2)},{\y}) node at (-7,-5) {$\Gamma_{1, \kpL}$};
\draw[domain=-4:4,smooth,variable=\y, thick]  plot ({-6-(\y*\y+1)^(1/2)},{\y}) node at (-9,-5) {$\Gamma_{2, \kpL}^-$};
\draw[domain=-4:4,smooth,variable=\y, thick]  plot ({-8-(\y*\y+1)^(1/2)},{\y}) node at (-11,-5) {$\Gamma_{3, \kpL}^-$};
\draw[domain=-4:4,smooth,variable=\y, thick]  plot ({(\y*\y+1)^(1/2)},{\y}) node at (4,-5) {$\Gamma_{3, \kpR}^+$};
\draw[domain=-4:4,smooth,variable=\y, thick]  plot ({2+(\y*\y+1)^(1/2)},{\y}) node at (6,-5) {$\Gamma_{2, \kpR}^+$};
\draw[domain=-4:4,smooth,variable=\y, thick]  plot ({4+(\y*\y+1)^(1/2)},{\y}) node at (8,-5) {$\Gamma_{1, \kpR}$};
\draw[domain=-4:4,smooth,variable=\y, thick]  plot ({6+(\y*\y+1)^(1/2)},{\y}) node at (10,-5) {$\Gamma_{2, \kpR}^-$};
\draw[domain=-4:4,smooth,variable=\y, thick]  plot ({8+(\y*\y+1)^(1/2)},{\y}) node at (12,-5) {$\Gamma_{3, \kpR}^-$};
\end{tikzpicture}
	\caption{Plots of contours when $m=3$. The contours are oriented from the bottom to the top.}
	\label{fig:contourInfTASEP}
\end{figure}

Let $\Gamma_{m, \kpL}^+, \cdots, \Gamma_{2, \kpL}^+, \Gamma_{1, \kpL}, \Gamma_{2, \kpL}^-, \cdots, \Gamma_{m, \kpL}^-$ be pairwise disjoint contours that are horizontal translates of each other in the left half of the complex plane extending from $\infty e^{\ii\theta_1}$ to $\infty e^{\ii \theta_2}$ for some angles $\theta_1\in (7\pi/6, 5\pi/4)$ and $\theta_2\in (3\pi/4, 5\pi/6)$ and are arranged as indicated in Figure~\ref{fig:contourInfTASEP}. 
In addition, let $\Gamma_{m, \kpR}^+, \cdots, \Gamma_{2, \kpR}^+,$ $\Gamma_{1, \kpR},$ $\Gamma_{2, \kpR}^-, \cdots, \Gamma_{m, \kpR}^-$ be the reflections of the above contours about the imaginary axis: see Figure~\ref{fig:contourInfTASEP}. 

Set 
$$
	\Gamma_{i, \kpL}\deff  \Gamma^{+}_{i, \kpL} \cup \Gamma^{-}_{i, \kpL}\quad \text{and}\quad \Gamma_{i, \kpR}\deff \Gamma^{+}_{i, \kpR} \cup \Gamma^{-}_{i, \kpR}, \qquad i=2,\hdots,m.
$$
In terms of the notations in Definition \ref{defn:whK12}, define
\begin{equation}\label{deff:Omegakpz-1}
	\Omega_{1, j}^\kpz\deff \begin{cases} \Gamma_{j, \kpL}, \quad &\text{if $j$ is odd,}\\
	\Gamma_{j, \kpR},\quad &\text{if $j$ is even,} \end{cases} 
	\qquad
	\Omega_{2, j}^\kpz\deff \begin{cases} \Gamma_{j, \kpR}, \quad &\text{if $j$ is odd,}\\
	\Gamma_{j, \kpL}, \quad &\text{if $j$ is even,} \end{cases} 
\end{equation}
and also $\Omega_\ell^\kpz \deff \bigcup_{j=1}^m \Omega_{\ell,j}^\kpz$ for $\ell=1,2$. 
Note that $\Omega_2^\kpz$ is equal to $-\Omega_1^\kpz$ if we reverse the orientation. 
The Hilbert spaces are $L^2(\Omega_\ell^\kpz, \dd z)$ for $\ell=1,2$.

For each $j=1, \cdots, m$, let\footnote{In \cite[Definition~2.25]{Liu19}, the formulas $Q_1(j)$ and $Q_2(j)$ are used. Here, we use the function  $\itQ_{j}^\kpz (z)=-Q_2(j) \chi_{\Omega_{1,j}}(z)+ Q_1(j) \chi_{\Omega_{2,j}}(z)$.} 
\begin{equation}\label{def:itQkpz}
	\itQ_{j}^\kpz (z)\deff \begin{cases}
	-(1-\zeta_j) \chi_{\Omega_{1, j}}(z) + \big( 1-\frac{1}{\zeta_{j-1}}\big) \chi_{\Omega_{2,j}}(z), & \text{if } j \text{ is odd ,} \\
	-\big( 1-\frac{1}{\zeta_{j-1}}\big)  \chi_{\Omega_{1, j}}(z) + (1-\zeta_j) \chi_{\Omega_{2,j}}(z), & \text{if } j \text{ is even,}
\end{cases}
\end{equation}
with the convention that $\zeta_0=\infty$,  $\zeta_m=0$, 
and\footnote{The factor $\itP_j^\kpz$ is not present in the kernels in \cite[Equation~(5) and Definition~2.4]{Liu19}. This is because the $L^2$ space is weighted. Here we use $L^2(\dd z)$ space instead and include  $\itP_j^\kpz$  as a part of the kernel.} 
\begin{equation}\label{def:itPkpz}
	\itP_j^\kpz(z)\deff
\begin{cases}
	\frac{1}{2\pi \ii}, & z\in \Gamma_{1}, \\
	\frac{1}{2\pi \ii( 1-\zeta_{j-1})}, \qquad & z\in \Gamma_{j,\kpL}^+\cup \Gamma_{j,\kpR}^+,\quad j>1,\\
	\frac{-\zeta_{j-1}}{2\pi \ii (1-\zeta_{j-1}) }, & z\in \Gamma_{j,\kpL}^-\cup \Gamma_{j,\kpR}^-, \quad j>1. 
\end{cases}
\end{equation}
Define the $m\times m$ diagonal matrix functions $\Am^\kpz$ and $\Bm^\kpz$ in \eqref{eq:Kbasicfuntios} with the functions  
\begin{equation}\label{def:ABkpz}
	\itA_j^\kpz(z) \deff 1 \quad \text{and}
	\quad \itB^\kpz_j(z) \deff \itQ_j^\kpz (z) \itP^\kpz_j(z). 
\end{equation}

Finally, set
\begin{equation}\label{def:mjexpon101}
	\genm_j^\kpz(z)\deff e^{ \itt_i  z^3 + \itl_i z^2 + \ith_i z }
	\quad
	\text{where} \quad \itt_i = - \pkt_i/3, \quad \itl_i= \pkga_i \quad \ith_i= \pkh_i, 
\end{equation}
and define $\Mm^\kpz$ as \eqref{eq:defnofMm}. 
Here, $\kpt_j>0$, $\kpga_j\in \R$, $\kph_j\in \R$ are the time, position, and height variable of the KPZ fixed point, respectively. This notation is consistent with \eqref{def:xyttohgammataukpz}.

The operators $\wh\itK_1^\kpz$ and $\wh\itK_2^\kpz$ are given by the form in Definition \ref{defn:whK12} with the above functions and matrices. 

\bigskip

We now apply Theorem~\ref{thm:algebraicmain} to the above operators. 
Let $\genH^\kpz$ be the operator in Definition~\ref{defn:notationsforH} with the specific choice of the conjugating constant $\genc(z)$ given by \eqref{eq:constantcL2}. 
The operator acts on $L^2(\Omega^\kpz, \dd z)$ where $\Omega^\kpz=\Omega_1^\kpz\cup \Omega_2^\kpz$.  
From the formula of the kernel, in particular due to the formula \eqref{def:mjexpon101} of $\genm_j^\kpz$, we see that $\genH^\kpz$ is a cubic integrable operator. 
Furthermore, since $\itA_j^\kpz$ and $\itB_j^\kpz$ do not depend on the parameters $\itt,  \itl, \ith$, it is indeed strongly cubic integrable.

\begin{cor}\label{prop:cubicadmissibleKPZfp}
The operator $\genH^\kpz$ is trace class on $L^2(\Omega^\kpz, \dd z)$ and 
\begin{equation}\label{eq:identFredDetkpz}
	\itD (\kph,\kpga,\kpt\mid \oz)  =   \det(\itid-\genH^\kpz).  
\end{equation}
\end{cor}

\begin{proof} 
It only remains to check  that $\genH^\kpz$ is trace class. 
We omit the superscript $\kpz$ to lighten the notations.  
We use Lemma~\ref{lem:traceclasslemma}. 
Decompose $\Omega_\ell$ so that each $\Sigma_{\ell, j}$ is one of the connected contours of $\Omega$; see Figure~\ref{fig:contourInfTASEP}. 
If $\Sigma_{1,j_1}$ and $\Sigma_{2,j_2}$ are in different half planes (such as $\Gamma_{1, \kpL}$ and $\Gamma_{3, \kpR}^-$), we take the separating contour to be $C_{j_1, j_2}=\ii \R$. 
On the other hand, if $\Sigma_{1,j_1}$ and $\Sigma_{2,j_2}$ are in the same half plane, then we take $C_{j_1, j_2}$ to be a translation of them that sits between them. 
To check condition (c), it is enough to check (i) and (ii) in the discussion after Lemma~\ref{lem:traceclasslemma}. 
Since $\Am$ and $\Bm$ are constants on each connected component, they are in $L^\infty(\Omega)$ and (i) holds. For (ii), due to the conditions \eqref{eq:KPZparameco} on the parameters, we see that 
\begin{equation}\label{eq:ckpz}
	\frac{\sqrt{\genm_{i+1}(z)}}{\sqrt{\genm_i(z)}}
	= e^{ \frac12 ( \gent_i - \gent_{i+1}) z^3 + \frac12 ( \geny_i - \geny_{i+1})  z^2 +\frac12  (\genx_i - \genx_{i+1})  z} 
\end{equation}
decays super-exponentially fast as $|z|\to \infty$ for $z$ in the sectors $\arg(z)\in (\pi/6, \pi/4)\cup(7\pi/4, 11\pi/6)$, 
and similarly, its reciprocal 
decays super-exponentially fast in the sectors $\arg(z)\in (3\pi/4, 5\pi/6) \cup (7\pi/6, 5\pi/4)$. 
From the geometry of the contours and this exponential decay, we find that the decay condition (ii) is satisfied, and we obtain the result.
\end{proof}

From the above result and the next lemma, Theorem~\ref{thm:genintoperred} (i) is proved.

\begin{lem}\label{prop:cubicadmissibleKPZfp22}
The symmetry \eqref{eq:symmetryprkpz} holds true.
\end{lem}

\begin{proof} 
For simplicity, we omit the superscript $\kpz$ in the notations. 
From \eqref{eq:expYk} and \eqref{eq:pjinfg}, we have $\itYb_1 =\lim_{z\to \infty} z (\itYb(z)-\itI)$. 
Thus, a symmetry of $\itYb_1$ follows from a symmetry of $\itYb(z)$, the solution of the Riemann-Hilbert problem in Subsection~\ref{sec:IIKSdiscussion}. 
This further follows from a symmetry of the jump matrix. 

From \eqref{eq:Jcon} and \eqref{eq:fg22}, the jump matrix is 
\beq \label{eq:jmptep}
	\itJ(z) =\itI_{m+1}-2\pi \ii \genf(z) \geng(z)^T = \itI_{m+1}-2\pi \ii \itconjb(z) \Fm(z)  \Gm(z)^T\itconjb(z)^{-1}
\eeq
where the $(m+1)$-dimensional vector functions $\Fm$ and $\Gm$ are given in \eqref{eq:ABcubicadm} 
with $\Am$ and $\Bm$ given in \eqref{def:ABkpz}. 
Since  $\Am=\itI$, we have
$$
	 \Fm(z)\Gm(z)^T= \itP
\begin{pmatrix} \itLambda_1 & \itLambda_2-\itE_{11}  \\ \itzero &  \itE_{11} \end{pmatrix}  
\begin{pmatrix}
\itzero &  \Bm(z)\genchi_1(z) \\  \Bm(z)\genchi_2(z) & \itzero
\end{pmatrix}
\begin{pmatrix} \itLambda_1 & \itLambda_2-\itE_{11}  \\ \itzero &  \itE_{11} \end{pmatrix}^{T}  	
\itP^T. 
$$
It is tedious but straightforward from the definition \eqref{def:ABkpz}, also noting that $\genchi_2(-z)=\genchi_1(z)$, to check that
\begin{equation}\label{eq:symm1} \begin{split}
	 \itB_j(-z) \chi_{\Omega_{2, j}}(-z) = 
	 \begin{dcases}
	 - \frac{1- \frac{1}{\oz_{j-1}}}{1-\oz_j} \itB_j(z) \chi_{\Omega_{1, j}}(z) \quad &\text{if $j$ is odd},\\
	 - \frac{1-\oz_j}{1-\frac{1}{\oz_{j-1}}} \itB_j(z) \chi_{\Omega_{1, j}}(z) \quad &\text{if $j$ is even},
	 \end{dcases}
\end{split} \end{equation}
for $j=1, \cdots, m$ and $z\in \Omega$.

For arbitrary non-zero constants $L_1, \cdots, L_m$, it is straightforward to check that 
$$
	\diag( L_1,\cdots, L_m )^{\pm 1}\itLambda_1=\itLambda_1 \genL_{\rm o}^{\pm 1}, 
	\qquad  
	\diag( L_1,\cdots, L_m )^{\pm 1}(\itLambda_2-\itE_{11})= (\itLambda_2 - \itE_{11})  \genL_{\rm e}^{\pm 1},
$$
with $m\times m$ matrices $\genL_{\rm o} \deff\diag(L_{1}, L_1 ,L_3, L_3,\cdots)$ 
and $\genL_{\rm e}\deff\diag(1, L_2, L_2, L_4, L_4, \cdots)$. 
Thus, with the $(m+1)\times (m+1)$ matrix  
\beqq
	\genL \deff \diag\left(L_1 ,\cdots, L_m,1\right), 
\eeqq
we find, noting $\itE_{11} \genL_{\rm e}^{\pm} =\itE_{11}$, that  
$$
	\genL \Fm(z)\Gm(z)^T \genL^{-1}= \itP
\begin{pmatrix} \itLambda_1 & \itLambda_2-\itE_{11}  \\ \itzero &  \itE_{11} \end{pmatrix}  
\begin{pmatrix}
\itzero & \genL_{\rm o} \genL_{\rm e}^{-1} \Bm(z)\genchi_1(z) \\ \genL_{\rm o}^{-1} \genL_{\rm e} \Bm(z)\genchi_2(z) & \itzero
\end{pmatrix}
\begin{pmatrix} \itLambda_1 & \itLambda_2-\itE_{11}  \\ \itzero &  \itE_{11} \end{pmatrix}^{T}  	
\itP^T . 
$$
Now for $L_j$ given in \eqref{eq:symmetryLkpz}, we have
\beqq
    \genL_{\rm o} \genL_{\rm e}^{-1}
    = \diag\left( L_1, \frac{L_1}{L_2}, \frac{L_3}{L_2}, \frac{L_3}{L_4}, \cdots\right)
    = \diag \left( - (1-\zeta_1), - \frac{1-\frac1{\zeta_1}}{1-\zeta_2}, 
    - \frac{1-\zeta_3}{1-\frac1{\zeta_2}}, - \frac{1-\frac1{\zeta_3}}{1-\zeta_4}, \cdots \right).
\eeqq
Thus, from \eqref{eq:symm1} we find that
\beqq
    \genL_{\rm o}^{-1} \genL_{\rm e} \Bm(-z)\genchi_2(-z) = \Bm(z)\genchi_1(z), \qquad 
    \genL_{\rm o} \genL_{\rm e}^{-1} \Bm(-z)\genchi_1(-z) = \Bm(z)\genchi_2(z)
\eeqq
and hence, 
\beq \label{eq:LFGL}
	\genL \Fm(-z)\Gm(-z)^T \genL^{-1}=( \Fm(z)\Gm(z)^T )^T.
\eeq
Noting that $\genm_j(-z\mid \ith,-\itl,\itt)^{-1}=\genm_j(z\mid \ith,\itl,\itt)$, we have $\itconjb(-z\mid \ith, -\itl,\itt)^{-1} = \itconjb(z\mid \ith, \itl,\itt)$. 
Thus, the jump matrix \eqref{eq:jmptep} satisfies 
$$
	\left(\genL\itJ (-z\mid \ith,-\itl,\itt)\genL\right)^{T}=\itJ(z\mid \ith,\itl,\itt).
$$
A symmetry of the jump matrix implies a symmetry of the solution of the RHP, in our case here it reads
$$
	\left(\genL\itYb(-z\mid \ith,-\itl,\itt)\genL^{-1}\right)^{-T}=\itYb(z\mid \ith,\itl,\itt).
$$ 
Finally, this implies the symmetry \eqref{eq:symmetryprkpz} for $\itYb_1$. 
\end{proof}

\subsection{Periodic KPZ fixed point - proof of Theorem~\ref{thm:genintoperred} (ii)}
 \label{sec:BLformula}
 
For the periodic KPZ fixed point, \cite[Section~2.2.3]{Baik-Liu19}  show that the formula in \eqref{eq:PKPZFstepdefn} is given by 
\begin{equation}\label{defn:pkDdisc}
	\pkDdisc (\kph,\kpga,\kpt\mid \oz)  \deff \det ( \itid - \wh \genK_1^\per \wh\genK^\per_2 ),
\end{equation}
for operators $\wh\genK_1^\per$ and $\wh\genK_2^\per$ of the form in Definition~\ref{defn:whK12} which we now describe. 
The series formula of the Fredholm determinant was shown to be convergent and the above formula is valid for the parameters satisfying
\beq \label{eq:pKPZparameco}
	0<\kpt_1\le \cdots \le \kpt_m \quad \text{and} \quad \text{$\kph_i<\kph_{i+1}$ if $\kpt_i=\kpt_{i+1}$.}
\eeq
The fixed constants $\zeta_1, \cdots, \zeta_m$ are arbitrary complex numbers satisfying $0<|\zeta_1|<\cdots< |\zeta_m|<1$.

\begin{figure}[t] \centering
\begin{tikzpicture}[scale=0.5]
\draw [line width=0.4mm,lightgray] (-4,0)--(4,0) node [pos=1,right,black] {$\R$};
\draw [line width=0.4mm,lightgray] (0,-3)--(0,3) node [pos=1,above,black] {$\ii\R$};
\draw[domain=-4:4,smooth,variable=\y,black]  plot ({(\y*\y+1)^(1/2)},{\y});
\draw[domain=-4:4,smooth,variable=\y,black]  plot ({-(\y*\y+1)^(1/2)},{\y});
\draw[domain=-4:4,smooth,variable=\y,black]  plot ({(\y*\y+5)^(1/2)},{\y});
\draw[domain=-4:4,smooth,variable=\y,black]  plot ({-(\y*\y+5)^(1/2)},{\y});
\draw[domain=-4:4,smooth,variable=\y,black]  plot ({(\y*\y+10)^(1/2)},{\y});
\draw[domain=-4:4,smooth,variable=\y,black]  plot ({-(\y*\y+10)^(1/2)},{\y});
\filldraw [fill=black!20!white, draw=black] (1.06573, 0.368479) circle[radius=3.5pt] node [above,shift={(0pt,0pt)}] {};   
\filldraw [fill=black!20!white, draw=black] (2.68227,2.48889) circle[radius=3.5pt] node [above,shift={(0pt,0pt)}] {};
\filldraw [fill=black!20!white, draw=black] (3.6699, 3.5311) circle[radius=3.5pt] node [above,shift={(0pt,0pt)}] {};
\filldraw [fill=black!20!white, draw=black] (2.532128,-2.326300) circle[radius=3.5pt] node [above,shift={(0pt,0pt)}] {}; 
\filldraw [fill=black!20!white, draw=black] (3.5615, -3.4182) circle[radius=3.5pt] node [above,shift={(0pt,0pt)}] {};
\fill (-1.06573,-0.368479) circle[radius=3.5pt] node [above,shift={(0pt,0pt)}] {};   
\fill (-2.68227,-2.48889) circle[radius=3.5pt] node [above,shift={(0pt,0pt)}] {};
\fill (-3.6699, -3.5311) circle[radius=3.5pt] node [above,shift={(0pt,0pt)}] {};
\fill (-2.532128,2.326300) circle[radius=3.5pt] node [above,shift={(0pt,0pt)}] {}; 
\fill (-3.5615, 3.4182) circle[radius=3.5pt] node [above,shift={(0pt,0pt)}] {};
\fill (2.242905, 0.1749963) circle[radius=3.5pt] node [above,shift={(0pt,0pt)}] {};   
\fill (3.102974, 2.151383) circle[radius=3.5pt] node [above,shift={(0pt,0pt)}] {};
\fill (3.9620469, 3.2707515) circle[radius=3.5pt] node [above,shift={(0pt,0pt)}] {};
\fill (2.983158,-1.974647) circle[radius=3.5pt] node [above,shift={(0pt,0pt)}] {};
\fill (3.8636663, -3.1508598) circle[radius=3.5pt] node [above,shift={(0pt,0pt)}] {}; 
\filldraw [fill=black!20!white, draw=black] (-2.242905, -0.1749963) circle[radius=3.5pt] node [above,shift={(0pt,0pt)}] {};   
\filldraw [fill=black!20!white, draw=black] (-3.102974, -2.151383) circle[radius=3.5pt] node [above,shift={(0pt,0pt)}] {};
\filldraw [fill=black!20!white, draw=black] (-3.9620469, -3.2707515) circle[radius=3.5pt] node [above,shift={(0pt,0pt)}] {};
\filldraw [fill=black!20!white, draw=black] (-2.983158, 1.974647) circle[radius=3.5pt] node [above,shift={(0pt,0pt)}] {};
\filldraw [fill=black!20!white, draw=black] (-3.8636663, 3.1508598) circle[radius=3.5pt] node [above,shift={(0pt,0pt)}] {}; 
\filldraw [fill=black!20!white, draw=black] (3.164709, 0.1240240) circle[radius=3.5pt] node [above,shift={(0pt,0pt)}] {};   
\filldraw [fill=black!20!white, draw=black] (3.652472, 1.827717) circle[radius=3.5pt] node [above,shift={(0pt,0pt)}] {};
\filldraw [fill=black!20!white, draw=black] (4.3462637, 2.9816117) circle[radius=3.5pt] node [above,shift={(0pt,0pt)}] {};
\filldraw [fill=black!20!white, draw=black] (3.567434,-1.651239) circle[radius=3.5pt] node [above,shift={(0pt,0pt)}] {};
\filldraw [fill=black!20!white, draw=black] (4.2615333, -2.8566879) circle[radius=3.5pt] node [above,shift={(0pt,0pt)}] {}; 
\fill (-3.164709, -0.1240240) circle[radius=3.5pt] node [above,shift={(0pt,0pt)}] {};    
\fill (-3.652472, -1.827717) circle[radius=3.5pt] node [above,shift={(0pt,0pt)}] {};
\fill (-4.3462637, -2.9816117) circle[radius=3.5pt] node [above,shift={(0pt,0pt)}] {};
\fill (-3.567434,1.651239) circle[radius=3.5pt] node [above,shift={(0pt,0pt)}] {};
\fill (-4.2615333, 2.8566879) circle[radius=3.5pt] node [above,shift={(0pt,0pt)}] {}; 
\end{tikzpicture}
	\caption{The nodes are  $\pkS_{i, \pkL}, \pkS_{i, \pkR}$ and the solid curves are the hyperbolas $\mathcal P_{i,\kpL}, \mathcal P_{i,\kpR}$ for $i=1, 2, 3$ for certain choices of $\zeta, \zeta_2, \zeta_3$. The black-filled nodes represent $\Omega_{1}^\per$ while the gray-filled nodes represent $\Omega_2^\per$. }
	\label{fig:PTASEP}
\end{figure}
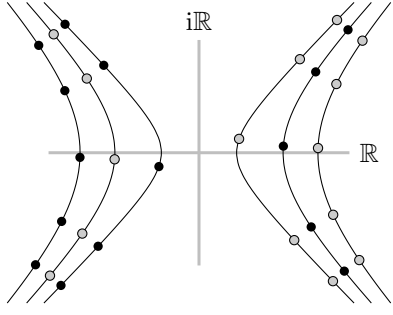

For each $i=1, \cdots, m$, consider the roots of the equation $e^{-s^2/2}= \zeta_i$, which are called  {\it Bethe roots}, and 
define the discrete sets  
\begin{equation}\label{deff:pkS}
 	\pkS_{i, \pkL}\deff\{ s\in \C\mid e^{-s^2/2}=\zeta_i, \; \re s<0 \},\qquad 
 	\pkS_{i, \pkR}\deff\{ s\in \C\mid e^{-s^2/2}=\zeta_i, \; \re s>0 \}.
\end{equation} 
Note that $\pkS_{i, \pkL}=-\pkS_{i, \pkR}$, and
\begin{equation}\label{eq:betheparabolas}
	\pkS_{i,\kpR}\subset \mathcal P_{i,\kpR}\quad \text{and}\quad \pkS_{i,\kpL}\subset \mathcal P_{i,\kpL},
\end{equation}
where 
\begin{equation}\label{deff:parabolasbetheroots}
	\mathcal P_{i,\kpR}\deff  \{s\in \C\mid |e^{-s^2/2}|=|\zeta_i|, \; \re s>0  \}=\left\{s\in \C\mid  (\re s)^2 - (\im s)^2=-2\log|\zeta_i|,    \; \re s>0 \right\} 
\end{equation}
and $\mathcal P_{i,\kpL}\deff -\mathcal P_{i,\kpR}$ are hyperbolas. 
Thus, the sets $\pkS_{i,\kpR}$ extend to $\infty$ along angles $7\pi/4$ and $\pi/4$ and the sets $\pkS_{i,\kpL}$ extend to $\infty$ with angles $5\pi/4$ and $3\pi/4$. See Figure \ref{fig:PTASEP}. 
In terms of the notation of \eqref{deff:genOmegaksets}, define the discrete sets 
\begin{equation}\label{def:pkS12}
	\Omega_{1, i}^\per\deff \begin{cases} \pkS_{i, \kpL}, \quad &\text{if $i$ is odd,}\\
	\pkS_{i, \kpR},\quad &\text{if $i$ is even,} \end{cases} 
	\qquad
	\Omega_{2, i}^\per\deff \begin{cases} \pkS_{i, \kpR}, \quad &\text{if $i$ is odd,}\\
	\pkS_{i, \kpL}, \quad &\text{if $i$ is even,} \end{cases} 
\end{equation} 
and 
\beqq
	\Omega_\ell^\per\deff\bigcup_{i=1}^m \Omega_{\ell,i}^\per
	\quad 
\eeqq
for $\ell=1,2$. Note that $\Omega_1^\per=-\Omega_2^\per$. 
The Hilbert spaces are $\ell^2(\Omega_\ell^\per)$ for $\ell=1,2$ with the counting measure.

For $i=1, \cdots, m$, set
 \begin{equation}\label{def:pQ}
 	\pkQ_{1,i}^\per\deff 1- \frac{\zeta_{i-(-1)^i}}{\zeta_i}, \qquad
 	\pkQ_{2,i}^\per\deff 1- \frac{\zeta_{i+(-1)^i}}{\zeta_i} 
 \end{equation}
with the convention that $\zeta_0=\zeta_{m+1}=0$. Define
 \begin{equation}\label{deff:pkH}
 	\pkH_i(z)\deff 
 	\begin{dcases}  \exp \bigg( - \int_{-\ii\infty}^{\ii\infty} \frac{\log(1-\zeta_i e^{w^2/2} ) }{w-z} \dd w \bigg)  
 	\quad &\text{for $\re(z)<0$,}\\
 	\exp \bigg( - \int_{-\ii\infty}^{\ii\infty} \frac{\log(1-\zeta_i e^{w^2/2} ) }{w+z} \dd w \bigg) 
 	\quad & \text{for $\re(z)>0$.} \end{dcases}
 \end{equation}
The function $\pkH_i(z)$ does not depend on the parameters, it is analytic and non-vanishing on $\C\setminus \ii\R$, and satisfies $\pkH_i(-z)=\pkH_i(z)$ and $\pkH_i(z)\to 1$ as $z\to \infty$. 
Now,  setting $\pkH_0=\pkH_{m+1}\equiv 1$,  define
$$
	\itA_j^\per(z)\deff \frac{1}{z\pkH_j(z)}\times \begin{dcases} \pkH_{j+1}(z) \quad & \text{for } z\in \pkS_{j, \pkL}, \\
	 \pkH_{j-1}(z) \quad & \text{for } z\in \pkS_{j, \pkR},
	\end{dcases}
$$
and
$$
	\pkB^\per_j(z)\deff \left(\itQ_{1,j}^\per \chi_{\pkS_1}(z)+ \itQ_{2,j}^\per \chi_{\pkS_2}(z)\right)\times\begin{cases} \pkH_{j+1}(z)  \quad & \text{for }  z\in \pkS_{j, \pkL}, \\
	\pkH_{j-1}(z)  \quad & \text{for }  z\in \pkS_{j, \pkR},
	\end{cases} 
$$
for $j=1,\hdots,m$, and define the matrices $\Am=\Am^\per$ and $\Bm=\Bm^\per$ by \eqref{eq:Kbasicfuntios}. 
Observe that $\pkA_j^\per,\pkB^\per_j$ are bounded on $\Omega^\per$.

Finally, define the functions
\beq
	 \genm_j^\per(z) \deff e^{ \itt_i  z^3 + \itl_i z^2 + \ith_i z } 
	\quad
	\text{where} \quad \itt_i = - \pkt_i/3, \qquad \itl_i= \pkga_i/2 \qquad \ith_i= \pkh_i,
\eeq 
and define the matrix $\genM^\per(z)$ as \eqref{eq:defnofMm}. 
Recall that $\kpt_j>0$, $\kpga_j\in [0,1)$, $\kph_j\in \R$ correspond to the time, position, and height variables, and which is consistent with \eqref{def:xyttohgammatau}. 
With these formulas, we define the operators $\wh \genK_1^\per$ and $\wh\genK^\per_2$ by the formula of Definition~\ref{defn:whK12}.

\bigskip

Let $\Omega^\per\deff \Omega_1^\per \cup \Omega_2^\per $ and define the operator $\genH^\per$ on $\ell^2(\Omega^\per_\ell)$ as Definition~\ref{defn:notationsforH} with the specific choice of the conjugating constant $\genc(z)$ given by \eqref{eq:constantcL2}. 
Then, from  the formula of $\genm_j^\per$, we find that $\genH^\per$ is a cubic admissible operator. 
The next result proves Theorem~\ref{thm:genintoperred} (ii).

\begin{cor}\label{prop:cubicadmissiblePKPZfp}
The operator $\genH^\per$ above is trace class on $\ell^2(\Omega^\per_\ell)$ and  
$$
	\pkDdisc (\kph,\kpga,\kpt\mid \oz)  =   \det(\itid-\genH^\per).  
$$
\end{cor}

\begin{proof} 
It is enough to show that $\genH^\per$ is trace class. 
We omit the superscript $\per$. 
We use Lemma~\ref{lem:traceclasslemma} and the discussion that follows. 
We already observed that $\Am$ and $\Bm$ are bounded on $\Omega$.

We take the decomposition \eqref{eq:decompositionOmega} of $\Omega_\ell$ using  $\Sigma_{\ell,j}=\Omega_{\ell,j}^\per$ for $\ell=1, 2$ and $j=1, \cdots, m$. 
If the sets $\Omega^\per_{1,j_1}$ and $\Omega^\per_{2,j_2}$ are in different half planes, then we take $C_{j_1,j_2}=\ii\R$ for the separating contour, and this case is easy, and we skip the detail. 
If the sets $\Omega^\per_{1,j_1}$ and $\Omega^\per_{2,j_2}$ are in the same half plane, then from the definition of the sets, we see that $j_1\neq j_2$. 
Thus, $|\zeta_{j_1}|\neq |\zeta_{j_2}|$. We take the separating contour as 
$$
	C_{j_1,j_2}= \left\{s\in \C\mid |e^{-s^2/2}|=(|\zeta_{j_1}|+|\zeta_{j_2}|)/2\right\}
$$
It is enough to check that 
\eqref{eq:condb1} and \eqref{eq:condb2} in the discussion after Lemma~\ref{lem:traceclasslemma} are satisfied. 
The hyperbola $C_{j_1,j_2}$ does not intersect the hyperbolas $\{s\in \C\mid |e^{-s^2/2}|= |\zeta_{j_1}| \}$ and $\{s\in \C\mid |e^{-s^2/2}|= |\zeta_{j_2}| \}$ (see \eqref{deff:parabolasbetheroots}) and all three hyperbolas have the same asymptotes. 
From the geometry we can check that 
\begin{equation}\label{eq:Cauchytr}
	\int_{C_{j_1,j_2}} \frac{ |\dd s| }{  |s-u|^2 }  
\end{equation}
is finite and grows at most like a polynomial as $u\to \infty$ on $\Omega_{\ell, j_1}^\per$ or $\Omega_{\ell, j_2}^\per$. 
We have 
\beq \label{eq:mratiodis}
	\frac{\sqrt{\genm_{j+1}(z)}}{\sqrt{\genm_j(z)}}
	= e^{ \frac12 ( \gent_j- \gent_{j+1}) z^3 + \frac12 ( \geny_j - \geny_{j+1})  z^2 +\frac12  (\genx_j- \genx_{j+1})  z} .
\eeq 
Note for $z\in \Omega^\per$, we have $e^{-z^2/2}= \zeta_i$ for some $i$, and thus, 
\beqq
	| e^{ \frac12 ( \geny_j - \geny_{j+1})  z^2} | = |\zeta_i|^{\geny_{j+1} - \geny_{j}} 
\eeqq
does not depend on $z$. 
Therefore, since the set $\Omega^\per$ has asymptotes to angles $\pi/4, 3\pi/4, 5\pi/4, 7\pi/4$, the conditions on \eqref{eq:pKPZparameco} on parameters imply that \eqref{eq:mratiodis} decays to zero exponentially fast on the set $\Omega_{\ell, j}^\per$ if it is on the right half plane. 
Similarly, the reciprocal of \eqref{eq:mratiodis} decays exponentially on $\Omega_{\ell, j}^\per$ if it is on the left half plane.
Due to exponential decay of $\frac{\sqrt{\genm_{j+1}(u)}}{\sqrt{\genm_j(u)}}$ or $\frac{\sqrt{\genm_{j}(u)}}{\sqrt{\genm_{j+1}(u)}}$ and the at most polynomial growth of \eqref{eq:Cauchytr}, we find that \eqref{eq:condb1} and \eqref{eq:condb2} are satisfied, and we obtain the result. 
\end{proof}

\section{Large height limits and the proof of Proposition~\ref{thm:kpformula}}\label{sec:TWformulas}

It is interesting to consider the asymptotic properties of the solutions to the differential equations in Subsection~\ref{sec:DEsection} 
as some of the parameters tend to infinity. 
In this paper, we only consider one simple case and leave other cases for a future work. 
The next result shows the asymptotic behavior of the Fredholm determinant and $\itYb_1$ for the (periodic) KPZ fixed points as the height parameters $\genx_i= \kph_i$ tend to positive infinity in a certain way. 
Since the result for both the KPZ and periodic KPZ fixed points are the same, we use the variables $\genx, \geny, \gent$ instead of the physical variables $\kph,\kpga,\kpt$: see \eqref{def:xyttohgammataukpz} and \eqref{def:xyttohgammatau} for the correspondence. 

\begin{lem}\label{lem:asymptfreddet}
Let either $\gendetD=\gendetD^\kpz$ and $\itYb_1=\itYb_1^\kpz$, or $\gendetD=\gendetD^\per$ and $\itYb_1=\itYb_1^\per$. 
Fix parameters $\genx,\geny,\gent$.
Let $b_1, \cdots, b_m>0$ be constants and
\beqq
	\veca\deff (a_1, a_2, \cdots, a_m), \qquad a_i\deff b_1+\cdots+ b_i.
\eeqq
Set
\beqq
	\genx_\xi\deff \genx+ \xi \veca= (\genx_1+ \xi \veca_1, \cdots, \genx_m+\xi \veca_m), \qquad \xi\in \R. 
\eeqq
Then, there exist constants $c>0$ and $\xi_0>0$, independent of $\genx,\geny,\gent$, for which the inequalities
$$
	|\gendetD(\genx_\xi,\geny,\gent)-1|\leq e^{-c \xi} \qquad \text{and}\qquad \max_{i,j=1,\cdots,m}|(\itYb_1(\genx_\xi,\geny,\gent))_{i,j}|\leq e^{-c\xi}
$$
hold true for all $\xi\geq \xi_0$.
\end{lem}

\begin{proof}
We omit the superscripts $\kpz$ and $\per$.
We also indicate the variable $z$ and the parameter $\xi$, but omit other parameters. 
In this proof, we use the alternative formula \eqref{eq:alternativefg} of the vectors $\genf$ and $\geng$. 
This formula involves $\Sm(z) \genchi_1(z)$ and $\Sm(z)^{-1} \genchi_2(z)$, which we first consider. 

From the definition in \eqref{def:mjexpon}, $\genm_i(z\mid \genx_\xi)= \genm_i(z\mid\genx) e^{\xi b_i z}$ holds. Thus, the matrix $\Sm$ in \eqref{def:Sm} satisfies  
\begin{equation} \label{eq:Goec0}
	\Sm(z\mid\genx_\xi) = \Sm(z\mid \genx) \bm R(z\mid  \xi)\qquad \text{with}\quad
	\bm R(z \mid \xi) \deff \diag ( e^{\xi b_1 z}, e^{-\xi b_2 z}, e^{\xi b_3 z}, \cdots). 
\end{equation}
Using the identity \eqref{eq:chimeasy} for indicator functions and noting that $\bm R$ and $\genchi_\ell$ are diagonal matrices, we find that 
\begin{equation}\label{eq:Goec2}
\begin{split}
	 &\genchi_1(z) \ite \ite^T( \genchi_1(z) \bm R(z \mid \xi) + \genchi_2(z) \bm R(z \mid \xi)^{-1})\ite
	 = \genchi_1(z) \bm R(z\mid \xi)   \ite = \bm R(z\mid \xi)  \genchi_1(z) \ite  , \\
	 &\genchi_2(z) \ite \ite^T( \genchi_1(z) \bm R(z \mid \xi) + \genchi_2(z) \bm R(z \mid \xi)^{-1} )\ite
	 = \genchi_2(z) \bm R(z\mid\xi)^{-1}   \ite =  \bm R(z\mid\xi)^{-1}  \genchi_2(z) \ite .
\end{split}
\end{equation}
The equations \eqref{eq:Goec0} and \eqref{eq:Goec2} imply that 
\begin{equation}\label{eq:Goec1}
	\Sm(z\mid\genx_\xi) \genchi_1(z) \ite =  \Sm(z\mid\genx) \genchi_1(z) \ite \phi(z\mid\xi) \quad \text{and} 
	\quad 
	\Sm(z\mid\genx_\xi)^{-1}  \genchi_2(z) \ite  =  \Sm(z\mid\genx)^{-1} \genchi_2(z) \ite \phi(z\mid\xi), 
\end{equation}
where
\beq
	\phi(z\mid\xi)  \deff  \ite^T( \genchi_1(z) \bm R(z \mid \xi) + \genchi_2(z) \bm R(z \mid \xi )^{-1})\ite 
\eeq
is a scalar function. 
Since $\Am$ and $\Bm$ do not depend on $\genx$, the formula \eqref{eq:alternativefg} implies that 
\begin{equation}\label{eq:fgR1R2}
	\genf(z\mid \genx_\xi)=\genf(z\mid \ith)\phi(z\mid\xi)\quad \text{and}\quad \geng(z\mid \genx_\xi) =\geng(z \mid \ith) \phi(z\mid \xi). 
\end{equation}

Now, $\phi(z\mid\xi)$ is a sum of $e^{ \pm b_j \xi z} \chi_{\Omega_{\ell, j}}(z)$ for $\ell=1, 2$ and $j=1, \cdots, m$ with the certain choice of the sign $\pm$. Being careful with the sign and noting that $\Omega_{\ell, j}$ is in either of the two halves of the complex plane, we find that 
\begin{equation}\label{eq:phixiineq}
	|\phi(z\mid \xi)| \le e^{- b\xi  |\re z |}, \qquad b\deff\min\{b_1, \cdots, b_m\},
\end{equation}
for every $z\in \Omega$ and every $\xi>0$.  
Since the sets in $\Omega$ are away from the imaginary axis by a positive distance, say $d$,  (see Figure \ref{fig:contourInfTASEP} and \ref{fig:PTASEP}) , this implies that  $|\phi(z\mid \xi)|\leq e^{-bd\xi}$ for every $z\in \Omega$ and every $\xi>0$.

Let $\genT_1=\genT_1(\ith_\xi)$ and $\genT_2=\genT_2(\ith_\xi)$ be the Hilbert-Schmidt operators from Lemma~\ref{lem:traceclasslemma}. From \eqref{eq:fgR1R2}--\eqref{eq:phixiineq} we find
$$
	\|\genH(\ith_\xi)\| \le \|\genH(\ith_\xi)\|_1\leq e^{-bd\xi}\|\genT_1(\ith)\|_2\|\genT_2(\ith)\|_2
$$
where the norms are the operator norm, the trace class norm, and the Hilbert-Schmidt norm, respectively. 
The norms $\|\genT_1(\ith)\|_2$ and $\|\genT_2(\ith)\|_2$ are finite and independent of $\xi$. 
Thus, from the inequality $\|\det(\itid-\genH)-1\|_1\leq \|\genH\|_1e^{1+\|\genH \|_1}$ for trace class operators (see \cite[Theorem~3.4]{Simon05}), 
the claim on $\gendetD=\det(\itid-\genH) $ then follows. 
On the other hand, from \eqref{defn:itYb1} and \eqref{eq:fgR1R2}--\eqref{eq:phixiineq}, we find 
$$
	\max_{i,j=1,\cdots, m}|(\itYb_1(\ith_\xi))_{i,j}|\leq  e^{-2bd\xi}\|(\itid-\genH(\ith_\xi))^{-1} \| \|\genf(\cdot\mid \ith)\|_{L^2(\Omega,\mu)}\|\geng(\cdot\mid \ith)\|_{L^2(\Omega,\mu)}.
$$
and the result for $\itYb_1$  follows. 
\end{proof}

Thus, in particular, the matrix functions $\genq, \genp, \genr, \gens$ for the (periodic) KPZ fixed points, which solve various differential equations, decay exponentially to zero as $\xi\to \infty$. 
Using the above lemma, we prove Proposition~\ref{thm:kpformula}. 

\begin{proof}[Proof of Proposition~\ref{thm:kpformula}]
Let $\gendetD$ be a cubic admissible determinant and $\genq_i$ the matrix-valued function in \eqref{eq:L1O1not} corresponding to the choice of subset $\itss=\{1, 2,\cdots, i\}$. 
Then, from Proposition~\ref{lem:basic2}, 
\beqq
	(\partial_{\ith_1}+ \cdots + \partial_{\ith_i}) \log\gendetD  =- \Tr\left( \genq_i \right) . 
\eeqq
Thus, for setting $\genx_\xi\deff \genx+ \xi \veca$ with $a=(1,\cdots,m)$ as in \eqref{deff:ithR} and $\xi\in \R$, we have 
\beq \label{eq:defiddDqkpz}
	\frac{\dd }{\dd \xi} \log\gendetD(\ith_\xi) 
	=\sum_{k=1}^m  k \partial_{\ith_k}\log\gendetD(\ith_\xi)
	=\sum_{i=1}^m  (\partial_{\ith_1}+ \cdots + \partial_{\ith_i}) \log\gendetD(\ith_\xi)
	=-\sum_{i=1}^m  \Tr  \genq_i (\genx_\xi) .
\eeq
By Lemma~\ref{lem:basic}, the above term is also equal to $\sum_{i=1}^m  \Tr\left( \gens_i (\genx_\xi)\right)$. 
Now, when $\gendetD$ is either $\gendetD^\kpz$ or $\gendetD^\per$, the last lemma shows that 
$\gendetD(\ith_\xi)\to 1$ and $\Tr \genq_i (\genx_\xi)\to 0$ exponentially fast as $\xi\to \infty$. 
Thus, integrating \eqref{eq:defiddDqkpz} we obtain Proposition~\ref{thm:kpformula}. 
\end{proof}


\section{Adler and van Moerbeke PDE and a non-equal time extension}\label{sec:symbolic}

Let  $\genH$ be a \underline{strongly} cubic integrable operator which is also trace class. Assume that $\itid-\genH$ is invertible. 
Set $M= \log \det(\itid-\itH)$ for the case when $m=2$ and with the change of variables given in \eqref{eq:twoptpara}, 
\beq \label{eq:AVcv}
	\itt_2= -\frac{\aitau}{3},\qquad 
	\itl_2 =\ail, \qquad
	\ith_1=\frac{\aiE+\aiW}{2},\qquad\ith_2=\frac{\aiE-\aiW}{2}-\frac{\ail^2}{\aitau}    
\eeq
with $\itt_1=-\frac{1}{3}$ and $\itl_1=0$.
Then, $M$ is a function of 4 variables, $\aitau, \aiE, \aiW$, and $\ail$. 
We now discuss how we can obtain partial differential equations \eqref{eq:maple5} and \eqref{eq:maple6} for $M$. 
As we discussed in Subsection~\ref{sec:AVPDE}, the special cases of these equations when $\aitau=1$ yield the equal-time differential equations for the Airy$_2$ process obtained by Adler and van Moerbeke \cite{Adler-van_Moerbeke05} and Quastel and Remenik \cite{Quastel-Remenik19b}. 

In the derivation of the differential equations for cubic admissible operators, we used Lax equations obtained from a Riemann-Hilbert problem. 
For the NLS, mKdV, and KP derivations, we used the Lax equations with respect to the parameters $\itt, \ith, \itl$, which in terms of \eqref{eq:AVcv} become Lax equations with respect to $\aitau, E, W$, and $y$. 
For the derivation of Tracy-Widom type ODE, we also considered the Lax equation with respect to the spectral variable $z$ for the case of strongly cubic admissible operators, and used the compatibility of $z$-Lax equation and $\ith$-Lax equation. 
In this section, we use the compatibility of the $z$-Lax equation with each of the $\aitau$-, $y$-, $E$-, and $W$-Lax equations. 

We start with an analysis of the $z$-Lax equation. 
The solution $\itYb(z)$ of the RHP in Subsection~\ref{sec:IIKSdiscussion} has the asymptotic series $\itYb(z)\sim \itI+ \frac{\itYb_1}{z} + \frac{\itYb_2}{z^2}+\cdots$ as $z\to \infty$. 
The asymptotic series can be written uniquely as 
\beq \label{eq:phioddef}
	\itYb(z)\sim\left( \itI+\sum_{k=1}^\infty \frac{\itYb^{(o)}_k}{z^k}\right)\exp\left( \sum_{k=1}^\infty \frac{\itYb^{(d)}_k}{z^k}\right)
\eeq
for  diagonal matrices  $\itYb^{(d)}_k$ and  off-diagonal matrices $\itYb^{(o)}_k$. 
Considering the $O(z^{-1})$ and $O(z^{-2})$ terms, we see that 
\beq \label{eq:phi12phiod}
	\itYb_1= \itYb^{(o)}_1 + \itYb^{(d)}_1 \quad \text{and} \quad \itYb_2= \itYb^{(o)}_2+ \itYb^{(o)}_1 \itYb^{(d)}_1+ \itYb^{(d)}_2+ \frac12 (\itYb^{(d)}_1)^2 .
\eeq

Consider the Lax equation \eqref{Z1prime} with respect to the spectral parameter $z$, where recall \eqref{eq:itconjb} that $\itW(z) \deff\itZ(z) \itconjb(z)$. We can write it as 
\beq \label{eq:Wzdertvs}
	\partial_z \itW(z) = ( z^2 \itQ_2 + z \itQ_1+ \itQ_0) \itW(z) 
\eeq
where, thanks to \eqref{eq:phi12phiod}, the coefficient matrices are
\beqq
	\itQ_2=3\whkptm,\quad 
	\itQ_1=2\whkpym + 3[\itYb^{(o)}_{1},\whkptm],\quad 
	\itQ_0=\whkpxm + 2 [ \itYb^{(o)}_{1}, \whkpym  ] +3[\itYb^{(o)}_{2},\whkptm] - 3[\itYb^{(o)}_{1},\whkptm] \itYb^{(o)}_{1},
\eeqq
and where we also recall \eqref{eq:multibytyx} for the definition of $\whkptm$, $\whkpym$, and $\whkpxm$. 
Inserting $\itW(z) =\itZ(z) \itconjb(z)$, the equation becomes, using \eqref{eq:Deltapatialz}, 
\beqq
	\partial_z \itZ(z) = ( z^2 \itQ_2 + z \itQ_1+ \itQ_0) \itZ(z) - \itZ(z)  (3z^2\whkptm + 2z\whkpym+\whkpxm ) .
\eeqq
We insert the asymptotic series \eqref{eq:phioddef} and for each $k=1,2, \cdots$ collect the coefficients of $z^{-k}$ to obtain an equation for $\itYb^{(o)}_i$ and $\itYb^{(d)}_i$. 
It is useful to note that $\itYb^{(d)}_i$, $\whkptm$, $\whkpym$, and $\whkpxm$ are diagonal matrices and commute. 
We find, after inserting the formula of $\itQ_0, \itQ_1, \itQ_2$, that for $k\ge 1$, 
\beq \label{recursion} \begin{split}
	&3 [ \itYb^{(o)}_{k+2} , \whkptm ] + 2 [\itYb^{(o)}_{k+1} , \whkpym ] 
	- 3[ \itYb^{(o)}_{k} , \whkptm ]  \itYb^{(o)}_{k+1}  - 3[\itYb^{(o)}_{k+1},\whkptm]  \itYb^{(o)}_{k}  + [ \itYb^{(o)}_{k} , \whkpxm ] \\
	&\quad - ( 2 [ \itYb^{(o)}_{k}, \whkpym  ] - 3[\itYb^{(o)}_{k},\whkptm] \itYb^{(o)}_{k}  ) \itYb^{(o)}_{k} 
	- (k-1)\itYb^{(o)}_{k-1}
	=   (k-1)\itYb^{(d)}_{k-1}  
	+ \sum_{j=2}^{k-1}(j-1)\itYb^{(o)}_{k-j} \itYb^{(d)}_{j-1} . 
\end{split} \eeq
Note that all matrices are of size $3\times 3$ since $m=2$. 
These equations are recursive equations. 
Note that $\itYb^{(o)}_{k+2}$ is the matrix of the highest index for the off-diagonal matrices and $\itYb^{(o)}_{k-1}$ is the matrix of the highest index for the diagonal matrices. 
It turns out that the recursion determines all matrices if we know $ \itYb^{(o)}_{1}$ and some entries of $\itYb^{(o)}_{2}$.

\begin{lem} \label{lem:recurdfflti}
Suppose that $a_1, \cdots, a_{10}$, are 10 functions of $\aitau, \aiE, \aiW, \ail$. 
Let 
\beq \label{eq:Yoff12}
	\itYb^{(o)}_{1}=\begin{pmatrix} 0 &a_1&a_3 \\ a_2 & 0&a_5\\ a_4&a_6&0\end{pmatrix}, \qquad 
	\itYb^{(o)}_{2}=\begin{pmatrix} 0 &a_7&c_1  \\ a_8 & 0&c_2\\ a_9&a_{10}&0\end{pmatrix},
\eeq
where $c_1,\,c_2$ are functions of $a_1, \cdots, a_{10}$ given by 
\beq \label{eq:c1c2Y10t1t2} \begin{split}
	&c_1=-\frac{a_3a_9}{a_4}-{\aitau a_1a_5}+\frac{\aitau a_1a_8}{a_4}+\frac{\aitau a_2a_7}{a_4}-\frac{2\ail a_1 a_2}{a_4}+a_1a_5+\frac{a_2a_3a_6}{a_4}-\frac{a_1a_8}{a_4}-\frac{a_2a_7}{a_4},\\& c_2=\frac{a_1a_4a_5}{a_6}+a_2a_3-\frac{a_1a_8}{a_6}-\frac{a_2a_7}{a_6}-\frac{a_5a_{10}}{a_6}+\frac{2\ail a_1 a_2}{\aitau a_6}+\frac{2\ail a_5 }{\aitau}-\frac{a_2a_3}{\aitau} +\frac{a_1a_8}{a_6\aitau}+\frac{a_2a_7}{a_6\aitau}.
\end{split} \eeq
Then, $\itYb^{(o)}_{n}$, $n\ge 3$, and $\itYb^{(d)}_{n}$, $n\ge 1$, are determined uniquely from the recursion \eqref{recursion}. 
Furthermore, all entries of $\itYb^{(o)}_n$ and $\itYb^{(d)}_{n}$ are Laurent polynomials of $a_1, \cdots, a_{10}$. 
\end{lem}

\begin{proof} 
Since $\whkptm$ is a diagonal matrix, the commutator of any matrix $W$ with $\whkptm$ is of the form 
\beq \label{eq:Mcom5}
	[ W , \whkptm ] = \begin{pmatrix} 0 &*& * \\ * & 0&* \\ * & * & 0\end{pmatrix} 
\eeq
and each $*$ is a constant multiple of the entry of $W$ at the same spot. 

Consider the case $k=1$ of \eqref{recursion}. The right-hand side of the equation is zero, so it becomes an equation of 
$\itYb^{(o)}_{1}$, $\itYb^{(o)}_{2}$ and $\itYb^{(o)}_{3}$ only,
$$
	3 [ \itYb^{(o)}_{3} , \whkptm ] + 2 [\itYb^{(o)}_{2} , \whkpym ] 
	- 3[ \itYb^{(o)}_{1} , \whkptm ]  \itYb^{(o)}_{2}  - 3[\itYb^{(o)}_{2},\whkptm]  \itYb^{(o)}_{1}  + [ \itYb^{(o)}_{1} , \whkpxm ]
	- ( 2 [ \itYb^{(o)}_{1}, \whkpym  ] - 3[\itYb^{(o)}_{1},\whkptm] \itYb^{(o)}_{1}  ) \itYb^{(o)}_{1} 
	= 0 . 
$$
Since $\itYb^{(o)}_{3}$ appears only in the term $[\itYb^{(o)}_{3} , \whkptm ]$, the diagonal entries of this (matrix) equation yield 3 scalar equations involving only the entries of $\itYb^{(o)}_{1}$ and $\itYb^{(o)}_{2}$. 
Two of them are \eqref{eq:c1c2Y10t1t2} and they imply the third equation. 
Now if we consider the off-diagonal entries of the equation, from \eqref{eq:Mcom5} we see that all 6 off-diagonal entries of $\itYb^{(o)}_{3}$ are Laurent polynomials of $a_1, \cdots, a_{10}$. Thus, $\itYb^{(o)}_{3}$, which is an off-diagonal matrix, is determined. 

We now prove the lemma using an induction on $k$. 
For $k\ge 2$, suppose that $\itYb^{(o)}_1, \cdots, \itYb^{(o)}_{k+1}$, and $\itYb^{(d)}_1, \cdots, \itYb^{(d)}_{k-2}$ are known and their entries are Laurent polynomials of $a_1, \cdots, a_{10}$. 
It is enough to show that we can determine $\itYb^{(o)}_{k+2}$ and $\itYb^{(d)}_{k-1}$ from the induction hypothesis, and their entries are Laurent polynomials of the quantities of the hypothesis. 
These two matrices appear in the equation \eqref{recursion} only in $3 [ \itYb^{(o)}_{k+2} , \whkptm ]$ and $(k-1)\itYb^{(d)}_{k-1}$, 
Thus, considering the diagonal entries of the equation, we can solve $\itYb^{(o)}_{k+2}$ in terms of the terms of the induction hypothesis, and the off-diagonal entries of the equation allows us to solve $ [ \itYb^{(o)}_{k+2} , \whkptm ]$, and thus $\itYb^{(o)}_{k+2}$ from  \eqref{eq:Mcom5}, in terms of the induction hypothesis. 
From the solution, we find that they are polynomials of the terms of the induction hypothesis, and thus the lemma is proved.  
\end{proof} 

Consider the Lax equations  \eqref{eq:tyxeqforZ} with respect to $\itt$, $\itl$,  and $\ith$.  With the  choice \eqref{eq:AVcv} of the parameters, they become four equations, 
\beqq \begin{split}
	 &
	 \partial_{\aiE} \itWsn(z) = \left( \frac{1}{2}z \itE_1+\frac{1}{2}z \itE_2+  \frac{1}{2}\itAsn_{1,1} +\frac{1}{2}\itAsn_{1,2}\right) \itWsn(z), \\
	 &\partial_{\aiW} \itWsn(z) = \left( \frac{1}{2}z \itE_1-\frac{1}{2}z \itE_2+   \frac{1}{2}\itAsn_{1,1}-\frac{1}{2}\itAsn_{1,2}  \right) \itWsn(z), \\
	 & 
	 \partial_{\ail} \itWsn(z) =\left( z^2 \itE_2+  z \itAsn_{1,2} -z \frac{2\ail}{\aitau}\itE_2+\itAsn_{2,2}- \frac{2\ail}{\aitau}\itAsn_{1,2} \right) \itWsn(z) ,\\&
	 \partial_{\aitau} \itWsn(z) = \left( -\frac{1}{3}z^3 \itE_2-\frac{1}{3}  z^2 \itAsn_{1,2} - \frac{1}{3}z \itAsn_{2,2}+z\frac{\ail^2}{\aitau^2}\itE_2 - \frac{1}{3} \itAsn_{3,2}+\frac{\ail^2}{\aitau^2}\itAsn_{1,2}\right) \itWsn(z),
\end{split}
\eeqq
where $\itAsn_{i,1}$ and $\itAsn_{i,2}$ are $\itAsn_{i}$ in \eqref{eq:itAssdef1} with $\itEn$ replaced by $\itE_1$ and $\itE_2$, respectively. 
For each equation, consider the compatibility condition with the $z$-Lax equation \eqref{eq:Wzdertvs}. 
If we consider the entries of these 4 equations (called zero-curvature equations; cf. the derivation of \eqref{eq:3eqnsfromzc}) yield 40 equations of the form 
\beqq \begin{split}
	&\partial_{{\aiE}}a_j=P^{(1)}_{j}(a),\quad \partial_{{\aiW}}a_j=P^{(2)}_{j}(a),\quad \partial_{{\ail}}a_j=P^{(3)}_{j}(a),\quad \partial_{{\aitau}}a_j=P^{(4)}_{j}(a), 
\end{split} \eeqq 
for $j=1, \cdots 10$, where $P^{(k)}_j(a)$ are Laurent polynomials. 
On the other hand, from Proposition \ref{lem:basic2} and Proposition \ref{prop:additionaldeformation}, the partial derivatives of $M= \log \det(\itid-\itH)$ are  polynomials of $\itZ_1$ and $\itZ_2$, and hence, by Lemma~\ref{lem:recurdfflti}, they and their higher partial derivatives are Laurent polynomials of $a_1, \cdots, a_{10}$. 

Motivated by the formula of the PDE of Adler and van Moerbeke, and of Quastel and Remenik, we consider the 27 partial derivatives  
\beqq\begin{split} 
&\partial_{{\itHe_a}}M,\quad 	\partial_{{\itHe_b}}\partial_{{\itHe_c}}M,\quad 
	\partial_{{\itHe_d}}\partial_{{\itHe_e}}\partial_{{\itHe_f}}M, \quad 
\partial_{{\itHe_h}}\partial_{{\itHe_h}}\partial_{{\itHe_i}}\partial_{{\itHe_j}}M,\quad 
\partial_{{\itHe_k}}\partial_{{\itHe_l}}M\partial_{{\itHe_m}}\partial_{{\itHe_n}}M, \\
&	\partial_{{\ail}}M,\quad	\partial_{{\ail}}\partial_{{\itHe_o}}M,\quad  
	\partial_{{\ail}}^2M, \quad 
	\partial_{{\aitau}}M,\quad \partial_{{\aitau}}\partial_{{\itHe_p}}M
\end{split}
\eeqq 
with the indices $a, b,c, \cdots$ being either $1$ and $2$, and $\itHe_1=E$, $\itHe_2=W$. 
These are Laurent polynomials of $a_1, \cdots, a_{10}$ involving 163 Laurent monomials. 
We then asked if there are linear relations between these 27 Laurent polynomials. 
We regarded each of 163 Laurent monomials as independent vectors and considered the vector space of Laurent polynomials spanned by them. 
Symbolic computations using Maple showed that the above 27 Laurent polynomials are not linearly independent and furthermore, the rank of their coefficient matrix is 23. As a result we obtain 4 linear relations of these Laurent polynomials. 
Two of them are \eqref{eq:maple5} and \eqref{eq:maple6}.
The other two are more complicated and we do not present them here.

\bibliographystyle{abbrv}  
\bibliography{bibliography}

\end{document}